\documentclass[amsfonts,amsmath,12pt]{article}
\usepackage{amsmath}
\usepackage{amsfonts}
\usepackage{microtype}
\newtheorem{theorem}{Theorem}[section]
\newtheorem{lem}[theorem]{Lemma}

\newtheorem{cor}[theorem]{Corollary}
\newtheorem{defn}[theorem]{Definition}
\newtheorem{res}[theorem]{Result}
\newenvironment{defn-new}{\begin{defn} \em}{\end{defn}}
\newtheorem{rem}[theorem]{Remark}
\newenvironment{rem-new}{\begin{rem} \em}{\end{rem}}
\newtheorem{ex}[theorem]{Example}
\newenvironment{ex-new}{\begin{ex} \em}{\end{ex}}

\newenvironment{notation-new}{\begin{rem} \em}{\end{rem}}

\newenvironment{agr-new}{\begin{rem} \em}{\end{rem}}

\makeatletter \@addtoreset{equation}{section} \makeatother

\makeatletter \@addtoreset{figure}{section} \makeatother

\begin{document}

\begin{center}
{\bf {\Large On H-Conformal Semi-invariant Submersion}}\\[0pt]

{\bf Punam Gupta}\footnote{
School of Mathematics, Devi Ahilya Vishwavidyalaya, Indore-452 001, M.P. 
INDIA\newline
Email: punam2101@gmail.com} and {\bf Kirti Gupta}\footnote{
Department of Mathematics \& Statistics, Dr. Harisingh Gour Vishwavidyalaya,
Sagar-470 003, M.P. INDIA \newline
Email: guptakirti905@gmail.com} 
\end{center}

\noindent {\bf Abstract:} We explore h-conformal semi-invariant submersions and almost h-conformal semi-invariant submersions originating from quaternionic K\"ahler manifolds to Riemannian manifolds. Our investigation focuses on the geometric characteristics of these submersions, including the integrability of distributions and the geometry of foliations. Additionally, we establish the necessary and sufficient conditions for such submersions to be totally geodesic. We also examine the equivalent conditions for the total manifold of the submersion to be twisted product manifold. Finally, we present a series of examples illustrating quaternionic K\"ahler manifolds and h-conformal semi-invariant submersions from quaternionic K\"ahler manifolds to Riemannian manifolds. \newline

\section{Introduction}

Riemannian submersion was first introduced by B. O`Neill \cite{O'Neill} in 1966, where he formulated the essential equations governing Riemannian submersions for Riemannian manifolds. This concept is crucial in the field of differential geometry and has numerous applications, including supergravity, statistical machine learning, statistical analysis on manifolds, and image processing. The geometric characteristics of Riemannian submersions have been explored by Besse \cite{Besse}, Escobales Jr. \cite{Eco}, and various other researchers.

B. Watson \cite{Watson} introduced and examined almost Hermitian submersions between almost complex manifolds, demonstrating that the horizontal and vertical distributions remain invariant with respect to the almost complex structure of the total space. In 2010, B. Sahin \cite{Sahin-10} investigated anti-invariant Riemannian submersions originating from almost Hermitian manifolds to a Riemannian manifold, positing that the fibers of these anti-invariant submersions are anti-invariant under the influence of the almost complex structure of the total space, indicating that the horizontal distribution does not maintain invariance under this action.

Conformal submersion serves as a natural extension of Riemannian submersion. The foundational theory of conformal submersions between Riemannian manifolds was initiated by Ornea and Romani\cite{Ornea}, and has since been explored by numerous researchers\cite{Akyol-Sahin,
Akyol}. The concept of horizontally conformal submersion, also a natural extension of Riemannian submersion, was independently introduced by Fuglede \cite{Fuglede} and Ishihara \cite{Ishi}. In 2020, Park \cite{Park-20} examined almost $h$-conformal semi-invariant submersions originating from almost quaternionic hermitian manifolds, yielding several noteworthy findings.

This manuscript examines the geometric characteristics of $h$-conformal semi-invariant submersions originating from quaternionic K\"ahler manifolds and mapping onto Riemannian manifolds. Section $2$ presents foundational concepts related to Riemannian submersions, conformal submersions, semi-invariant submersions and quaternionic K\"ahler manifolds. In Section $3$, we study the $h$-conformal semi-invariant submersion, $h$-homothetic and analyze the integrability of distributions, the geometry of foliations, and the criteria for such submersions to be totally geodesic, along with the necessary and sufficient conditions for the integrability of distributions. Furthermore, we explore the geometry of the leaves associated with horizontal and vertical distributions. In section $4$, we offer numerous examples to improve the understanding of these submersions.

\section{Preliminaries}
This section revisits fundamental definitions and presents established results related to Riemannian submersions and horizontally conformal submersions, which will be utilized in the paper.

\subsection{Submersion}

We will now present a set of definitions and results for future reference.

\begin{defn-new}
{\rm \cite{Baird}} A $C^{\infty }$-map $F:\left(
M,g_{M}\right)  \rightarrow \left( N,g_{N}\right) $ be  between Riemannian  manifolds $\left( M,g_{M}\right) $ and $\left( N,g_{N}\right) $.is said to be  $C^{\infty }$-submersion if

\begin{description}
\item[(i)] $F$ is surjective. 

\item[(ii)] The differential $\left(F_*\right)_p$ has maximal rank for any $
p \in M$.\newline
\end{description}
\end{defn-new}

\begin{defn-new}
{\rm {\cite{Falcitelli-Ianus and Pastore, O'Neill}}}
Let $(M,g_{M})$ and $
(N,g_{N})$ be Riemannian manifolds with $\dim (M)=m>n=\dim (N)$. A smooth surjective map $F :M\rightarrow N$ is said to be Riemannian submersion if $F$ satisfies the
following axioms:

\begin{description}
\item[(i)] $F$ has maximal rank.

\item[(ii)] The differential $F_{\ast }|_{(\operatorname{ker} F
_{\ast })^{\perp }}$ is a linear isometry.
\end{description}
\end{defn-new}

\noindent{\bf Remark:} For each $q\in N$, $F ^{-1}(q)$ is an $(m-n)$-dimensional submanifold of $%
M $. The submanifolds $F ^{-1}(q)$, $q\in N$, are called fibers. A vector
field on $M$ is called vertical if it is always tangent to the fibers. A vector
field on $M$ is called horizontal if it is always orthogonal to the fibers. A
vector field $X$ on $M$ is called basic if $X$ is horizontal and $F $%
-related to a vector field $X^{\prime }$ on $N$, that is, $F _{\ast
}X_{p}=X_{F _{\ast }(p)}^{\prime }$ for all $p\in M.$ We denote the
projection morphisms on the distributions $\operatorname{ker}F _{\ast }$ and $(\operatorname{ker} F
_{\ast })^{\perp }$ by ${\cal V}$ and ${\cal H}$, respectively. The sections
of ${\cal V}$ and ${\cal H}$ are called vertical vector fields and
horizontal vector fields, respectively. So
\[
{\cal V}_{p}=T_{p}\left( F ^{-1}(q)\right) ,\qquad {\cal H}%
_{p}=T_{p}\left( F ^{-1}(q)\right) ^{\perp }.
\]

\medskip

\begin{defn-new}
{\rm {\cite{Ornea}}}\ Let $(M,g_{M})$ and $(N,g_{N})$ be Riemannian 
manifolds with $\dim (M)=m>n=\dim (N)$. A horizontally conformal submersion $F:(M,g_{M})\rightarrow (N,g_{N})$ is a map of $M$ onto
$N$ satisfying the following conditions: 

\begin{description}
\item[(i)] $F$ has maximal rank. 

\item[(ii)] The angle between the horizontal vectors is preserved by the 
differential $F_{\ast }$,  
\[
g_{N}(F_{\ast }V,F_{\ast }W)=\lambda ^{2}(p)g_{M}(V,W),\quad V,W\in \Gamma\left( {%
 \operatorname{ker}}\left( F_{\ast }\right) _{p}\right) ^{\perp},  
\]
\noindent where $\lambda (p)$ is said to be dilation of $F$ at $p$. 
\end{description}
\end{defn-new}
From now on, we use the $h$-conformal submersion instead of the horizontally conformal submersion.
\begin{defn-new}
{\rm {\cite{Baird}}}\ An $h$-conformal submersion $F:(M,g_{M})\rightarrow 
(N,g_{N})$ is called $h$-homothetic if the slope of its dilation $\lambda $ 
is vertical, that is,  
\[
{\cal H}(\nabla \lambda )=0  
\]
at $p\in M$, where ${\cal H}$ is the projection on the horizontal space $
\left( \operatorname{ker} F_{\ast }\right) ^{\perp }$. 
\end{defn-new}

\begin{defn-new}
{\rm \cite{Baird}} \ Let $(M,g_{M})$ and $(N,g_{N})$ be Riemannian 
manifolds with $\dim (M)=m>n=\dim (N)$. A map $F:(M,g_{M})\rightarrow (N,g_{N})$ is called horizontally weakly conformal or 
semi-conformal at $p\in M$ if it satisfies either {\rm (i)} $\left( F_{\ast
}\right) _{p}=0$ or {\rm (ii)} $\left( F_{\ast }\right) _{p}$ is surjective and 
there exists a number $\lambda (p)>0$ such that  
\[
g_{N}\left( \left( F_{\ast }\right) _{p}V,\left( F_{\ast }\right) 
_{p}W\right) =\lambda ^{2}g_{M}(V,W)\quad V,W\in \Gamma\left( {\operatorname{ker}}\left( F_{\ast
}\right) _{p}\right) ^{\perp }  
\]
and the $\lambda (p)$ is said to be dilation of $F$ at $p$. 
\end{defn-new}

The point $p\in M$ is said to be a critical point if it satisfies the type (i)
and  the point $p\in M$ is said to be a regular point if it satisfies the
condition (ii). The map $F$  is called horizontally weakly conformal if it
is horizontally weakly  conformal at any point of $M$. It is said to be a
conformal submersion if $F$  has no critical point.

Let $F:\left( M,g_{M}\right) \rightarrow \left( N,g_{N}\right) $ 
be an $h$-conformal submersion. Given any vector field \\$U\in \Gamma (TM)$, 
we have  
\[
U={\cal V}U+{\cal H}U,  
\]
where ${\cal V}U\in \Gamma \left( {\operatorname{ker}}F_{\ast }\right)$ and
${\cal H}U\in \Gamma \left( \left( {\operatorname{ker}}F_{\ast }\right) ^{\perp }\right) $.

The second fundamental tensors of all fibers $F ^{-1}(q),\ q\in N$ give
rise to configuration tensors $T$ and $A$ in $M$ defined by O'Neill 
\cite{O'Neill}.  
\begin{eqnarray}
{\cal A}_{E}F &=&{\cal H}\nabla _{{\cal H}E}{\cal V}F+{\cal V}\nabla _{{\cal %
H}E}{\cal H}F, \\
{\cal T}_{E}F &=&{\cal H}\nabla _{{\cal V}E}{\cal V}F+{\cal V}\nabla _{{\cal %
V}E}{\cal H}F
\end{eqnarray}
for $E,F\in \Gamma (TM)$, where $\nabla $ is the Levi-Civita connection of $
g_{M}$.

From the above equations, we have  
\begin{equation}
\nabla _{X}Y={\cal T}_{X}Y+\widehat{\nabla }_{X}Y,
\end{equation}
\vspace{-0.7cm}
\begin{equation}
\nabla _{X}V={\cal H}\nabla _{X}V+{\cal T}_{X}V ,  \label{eq-aa}
\end{equation}
\begin{equation}
\nabla _{V}X={\cal A}_{V}X+{\cal V}\nabla _{V}X,
\end{equation}
\begin{equation}
\nabla _{V}W={\cal H}\nabla _{V}W+{\cal A}_{V}W,
\end{equation}
for $X,Y\in \Gamma  \left( {\operatorname{ker}}F_{\ast }\right)$ and $V,W\in
\Gamma \left( \left( {\operatorname{ker}}F_{\ast }\right) ^{\perp }\right) $, where $\widehat{
\nabla}_{X}Y={\cal V}\nabla _{X}Y$. If $V$ is basic, then ${\mathcal{A}}_{V}X={\cal H} 
\nabla _{X}V.$

It is easy to find that for $p\in M,$ $X\in {\cal V}_{p}$ and $U\in {\cal H}
_{p}$, the linear operators  
\[
{\mathcal{T}}_{X},{\mathcal{A}}_{U}:T_{p}M\rightarrow T_{p}M  
\]
are skew-symmetric.

Let $F:\left( M,g_{M}\right) \rightarrow \left( N,g_{N}\right) $ be a $
C^{\infty }$-map between Riemannian manifolds. Then the differential $
F_{\ast }$ of $F$ can be observed as a section of the bundle $
Hom(TM,F^{-1}TN)\rightarrow M$, where $F^{-1}TN$ is the bundle which has 
fibres $\left( F^{-1}TN\right) _{p}=T_{f(p)}N$, has a connection $\nabla $ 
induced from the Riemannian connection $\nabla ^{M}$ and the pullback 
connection $\nabla ^{N}$. Then the second fundamental form of $F$ is given by  
\begin{equation}
\left( \nabla F_{\ast }\right) (X,Y):=\nabla _{X}^{N}F_{\ast }Y-F_{\ast
}\left( \nabla _{X}^{M}Y\right) \text{ for }X,Y\in \Gamma (TM).
\label{eq-bb}
\end{equation}

\begin{rem}
  It is easy to check that 
 \begin{itemize}
 \item The second fundamental form $\nabla F_{\ast }$ is 
symmetric. 
\item If $X\in \Gamma \left( {\operatorname{ker}}F_{\ast }\right) $ and $V\in \Gamma 
\left( \left( {\operatorname{ker}}F_{\ast }\right) ^{\perp }\right)$, then $\lbrack 
X,V]\in \Gamma \left( {\operatorname{ker}}F_{\ast }\right).$
 \end{itemize}   
\end{rem}

\begin{defn-new}
{\rm {\cite{Baird}}}\ A map $F$ is said to be harmonic if the tension field 
$\tau(F)={trace}\left(\nabla F_*\right)=0$ and $F$ is said to be a totally 
geodesic map if $(\nabla F_{*})(X, Y)=0$ for $X, Y \in \Gamma(T M)$. 
\end{defn-new}

Now, we state the result given by researchers for further use. 

\begin{res}
    \label{2.10}
{\rm \cite{Urakawa}} Let $F$ : $\left( M,g_{M}\right) \rightarrow \left( 
N,g_{N}\right) $ be a $C^{\infty }$-map between Riemannian manifolds. Then
for $X,Y\in \Gamma (TM)$  
\[
\nabla _{X}^{N}F_{\ast }Y-\nabla _{Y}^{N}F_{\ast }X-F_{\ast }([X,Y])=0.  
\]
\end{res}


\begin{res}
{\rm \cite{Baird}} Let $F:\left( M,g_{M}\right) \rightarrow \left( 
N,g_{N}\right) $ be an $h$-conformal submersion with dilation $\lambda $. 
Then  
\begin{equation}
\left( \nabla F_{\ast }\right) (V,W)=V(\ln \lambda )F_{\ast }W+W(\ln \lambda
)F_{\ast }V-g_{M}(V,W)F_{\ast }(grad \ln \lambda )  \label{eq-cc}
\end{equation}
for $V,W\in \Gamma \left( \left( {\operatorname{ker}}F_{\ast }\right) ^{\perp }\right) $. 
\end{res}

\subsection{Quaternionic K\"ahler manifold}

Let $\left( M,g_{M},J\right) $ be an almost Hermitian manifold, where $J$ is
an almost complex structure on $M$, that is,  
\[
J^{2}=-id,\hspace*{0.15cm}g_{M}(JX,JY)=g_{M}(X,Y)\hspace*{0.15cm}for\hspace*{%
		0.1cm}X,Y\in \Gamma (TM).  
\]
Let $M$ be a $4m$-dimensional $C^{\infty }$-manifold and let $E$ be a rank 3
subbundle of End $(TM)$ such that for any point $p\in M$ with a neighborhood
$U$, there exists a local basis $\left\{ J_{1},J_{2},J_{3}\right\} $ of 
sections of $E$ on $U$ satisfying  
\[
J_{a}^{2}=-id,\quad J_{\alpha }J_{\alpha +1}=-J_{\alpha +1}J_{\alpha
}=J_{\alpha +2}  
\]
for all $\alpha \in \{1,2,3\}$,  where the indices are taken from $\{1,2,3\}$
modulo $3$. Then $E$ is said to  be an almost quaternionic structure on $M$
and $(M,E)$ an almost  quaternionic manifold \cite{Alekseevsky-Marchiafava}.

Moreover, let $g_M$ be a Riemannian metric on $M$ defined by  
\[
g_M\left( J_{\alpha }X,J_{\alpha }Y\right) =g_M(X,Y)  
\]
for all vector fields $X,Y\in \Gamma (TM)$, where the indices are taken from
$\{1,2,3\}$ modulo $3$. Then $(M,E,g_M)$ is said to be an almost quaternionic Hermitian manifold \cite%
{Ianus-Mazzocco-Vilcu} and the basis $\left\{ J_{1},J_{2},J_{3}\right\} $ is
said to be a quaternionic Hermitian basis.

An almost quaternionic Hermitian manifold $(M,E,g_M)$ is said to be a 
quaternionic K\"{a}hler manifold \cite{Ianus-Mazzocco-Vilcu} if there exist 
locally defined $1$-forms $\omega _{1},\omega _{2},\omega _{3}$ such that 
for $\alpha \in \{1,2,3\}$ modulo $3$  
\begin{equation}
\nabla _{X}J_{\alpha }=\omega _{\alpha +2}(X)J_{\alpha+1}-\omega _{\alpha
+1}(X)J_{\alpha +2}  \label{eq-qkm}
\end{equation}
for $X\in \Gamma (TM)$, where the indices are taken from $\{1,2,3\}$ modulo $
3$.

If there exists a global parallel quaternionic Hermitian basis $\left\{ 
J_{1},J_{2},J_{3}\right\} $ of sections of $E$ on $M$ (i.e. $\nabla J_{a}=0$
for $\alpha \in \{1,2,3\}$, where $\nabla $ is the Levi-Civita connection of
the metric $g_M$), then $(M,E,g_M)$ is said to be a hyperk\"ahler manifold \cite%
{Besse}.
\begin{defn-new}
{\rm \cite{Sahin}} Let $(M,g_{M},J)$ be an almost Hermitian manifold and $(N,g_{N})$ be a Riemannian manifold. A Riemannian submersion $F:(M,g_{M})\rightarrow (N,g_{N})$ is called semi-invariant if there is a distribution ${\mathcal{D}}_1 \subseteq {\operatorname{ker}}F_{\ast}$ such that
$$
\operatorname{ker} F_*={\mathcal{D}}_1 \oplus {\mathcal{D}}_2
$$
and
$$
J\left({\mathcal{D}}_1\right)={\mathcal{D}}_1,\quad  \quad J\left({\mathcal{D}}_2\right) \subseteq\left(\operatorname {ker} F_{\ast}\right)^{\perp},
$$
where ${\mathcal{D}}_2$ is orthogonal complementary to ${\mathcal{D}}_1$ in $ \operatorname{ker}F_{\ast}$.
\end{defn-new}

\begin{rem}
Every anti-invariant Riemannian submersion from an almost Hermitian manifold onto a Riemannian manifold is a semi-invariant Riemannian submersion with ${\mathcal{D}}_1=\{0\}. $    
\end{rem}


\begin{defn-new} 
{\rm\cite{Akyol-Sahin-17}} An $h$-conformal submersion $F: (M,g_{M},J) \mapsto (N,g_{N})$ is called a conformal semi-invariant submersion if there is a distribution ${\mathcal{D}}_{1} \subset \operatorname{ker}F_{*}$ such that 
$$
\operatorname{ker} F_*={\mathcal{D}}_1 \oplus {\mathcal{D}}_2, 
J\left({\mathcal{D}}_1\right)={\mathcal{D}}_1, J\left({\mathcal{D}}_2\right) \subset\left(\operatorname{ker} F_*\right)^{\perp} \text {, }
$$
where ${\mathcal{D}}_2$ is the orthogonal complement of ${\mathcal{D}}_1$ in $\operatorname{ker}F_*$. 
\end{defn-new}

\begin{defn-new}
\cite{Park-12} Let $\left(M, E, g_M\right)$ be an almost quaternionic Hermitian manifold and $\left(N, g_N\right)$ a Riemannian manifold. A Riemannian submersion $F:\left(M, E, g_M\right) \mapsto\left(N, g_N\right)$ is called an $h$-semi-invariant submersion if given a point $p \in M$ with a neighborhood $U$, there exists a quaternionic Hermitian basis $\{J_{1}, J_{2}, J_{3}\}$ of sections of $E$ on $U$ such that for any $ \alpha \in\{1, 2, 3\}$, there is a distribution ${\mathcal{D}}_1 \subset \operatorname{ker} F_*$ on $U$ such that
$$
\operatorname{ker} F_*={\mathcal{D}}_1 \oplus {\mathcal{D}}_2, J_{\alpha}\left({\mathcal{D}}_1\right)={\mathcal{D}}_1, J_{\alpha}\left({\mathcal{D}}_2\right) \subset\left(\operatorname{ker} F_*\right)^{\perp} \text {, }
$$
where ${\mathcal{D}}_2$ is the orthogonal complement of ${\mathcal{D}}_1$ in $\operatorname{ker}F_*$.
We call such a basis $\{J_{1}, J_{2}, J_{3}\}$ an $h$-semi-invariant basis. 
\end{defn-new}

\begin{defn-new}
\cite{Park-12} A Riemannian submersion $F:\left( M,E,g_{M}\right) \rightarrow \left( 
N,g_{N}\right) $ is called an almost $h$-semi-invariant submersion if 
given a point $p\in M$ with a neighborhood $U$, there exists a quaternionic 
Hermitian basis $\{J_{1},J_{2},J_{3}\}$ of sections of $E$ on $U$ such that for each $ \alpha \in\{1, 2, 3\}$, there is a distribution ${{\mathcal{D}}_1}^{J_{\alpha}} \subset \operatorname{ker} F_*$ on $U$ such that
$$
\operatorname{ker} F_*={{\mathcal{D}}_1}^{J_{\alpha}} \oplus {{\mathcal{D}}_2}^{J_{\alpha}}, J_{\alpha}\left({{\mathcal{D}}_1}^{J_{\alpha}}\right)={{\mathcal{D}}_1}^{J_{\alpha}}, J_{\alpha}\left({{\mathcal{D}}_2}^{J_{\alpha}}\right) \subset\left(\operatorname{ker} F_*\right)^{\perp} \text {, }
$$
where ${{\mathcal{D}}_2}^{J_{\alpha}}$ is the orthogonal complement of ${{\mathcal{D}}_1}^{J_{\alpha}}$ in $\operatorname{ker} F_*$.
We call such a basis $\{J_{1}, J_{2}, J_{3}\}$ an almost $h$-semi-invariant basis.  
 \end{defn-new}
 
\begin{defn-new}
Let $F$ be an almost $h$-semi-invariant submersion from an almost quaternionic Hermitian manifold $(M,E,g_{M})$ onto a Riemannian manifold $(N,g_{N}).$ If $ J_{\alpha}({\mathcal{D}}_{2}^{J_{\alpha}}) = \left(\operatorname{ker} F_*\right)^{\perp}$ for ${\alpha} \in \{1,3\}$ and $J_{2}\left(\operatorname{ker} F_*\right) = \operatorname{ker} F_*$  (i.e., $ {\mathcal{D}}_{2}^{J_{2}}=\{0\}$), then we call the map $F$ an anti-holomorphic almost $h$-semi-invariant submersion. We call such a basis $\{J_{1},J_{2},J_{3}\}$ an anti-holomorphic almost $h$-semi-invariant basis.
\end{defn-new}

\begin{defn-new}
Let $F$ be an almost $h$-semi-invariant submersion from a quaternionic K\"{a}hler manifold $(M,J_{1},J_{2},J_{3},g_{M})$ onto a Riemannian manifold $(N,g_{N})$ such that $\{J_{1}, J_{2}, J_{3}\}$ is an almost $h$-semi-invariant basis. We call the distribution ${\mathcal{D}}_{2}^{J_{\alpha}}$ for $\alpha \in \{1,2,3\}$ parallel along $ \left( \operatorname{ker}{F_{\ast}}\right)$ if $ \nabla_{X}V \in \Gamma \left( {\mathcal{D}}_{2}^{J_{\alpha}} \right)$ for $X \in \Gamma \left(\operatorname{ker} F_*\right)^{\perp}$ and $V \in \Gamma \left( {\mathcal{D}}_{2}^{J_{\alpha}} \right).$
\end{defn-new}

\begin{defn-new} Let $F$ be an almost h-conformal semi-invariant submersion from a quaternionic K\"ahler manifold $(M,J_{1},J_{2},J_{3},g_{M})$ onto a Riemannian manifold $(N,g_{N})$ such that $\{J_{1}, J_{2}, J_{3}\}$ is an almost h-conformal semi-invariant basis. Given $R \in \{J_{1},J_{2},J_{3} \},$ we call the distribution ${\cal{D}}_{2}^{R}$ parallel along $\left( \operatorname{ker} F_{\ast}\right)^{\perp}$ if $ \nabla_{X}V \in \Gamma({\cal{D}}_{2}^{R})$ for $X \in \Gamma \left( \left( \operatorname{ker} F_{\ast}\right)^{\perp}\right)$ and $ V \in \Gamma({\cal{D}}_{2}^{R}). $
\end{defn-new}

\begin{defn-new}
 Let $F$ be an almost $h$-conformal semi-invariant submersion from a quaterni-onic-K\"ahler manifold $(M,J_{1},J_{2},J_{3},g_{M})$ onto a Riemannian manifold $(N,g_{N})$ such that $\{J_{1}, J_{2}, J_{3}\}$ is an almost $h$-conformal semi-invariant basis. Then given $R \in \{J_{1},J_{2},J_{3} \},$ we call the distribution $\mu^{R}$ parallel along $\operatorname{ker} F_{\ast}$ if $\nabla_{U}X \in \Gamma(\mu^{R})$ for $U \in \Gamma \left( \operatorname{ker} F_{\ast}\right)$ and $X \in \Gamma (\mu^{R}).$  
\end{defn-new}
\begin{defn-new}\cite{Zawa}
    Let $F: (M,g_{M}) \rightarrow (N,g_{N})$ be an $h$-conformal submersion. The map $F$ is called an $h$-conformal submersion with totally umbilical fibres if 
    $$ {\cal{T}}_{X}Y = g_{M}(X,Y)H  \text{  for }  X,Y \in \Gamma \left( \operatorname{ker} F_{\ast} \right),$$
    where $H$ is the mean curvature vector field of the distribution $\operatorname{ker} F_{\ast}.$
\end{defn-new}
\begin{defn-new}
Let $F$ be an almost $h$-conformal semi-invariant submersion from a quaternionic K%
\"{a}hler manifold $\left( M,J_{1},J_{2},J_{3},g_{M}\right) $ onto a 
Riemannian manifold $\left( N,g_{N}\right) $ such that $\{J_{1},J_{2},J_{3} 
\} $ is an almost $h$-conformal semi-invariant basis. For $R\in \left\{ 
J_{1},J_{2},J_{3}\right\} $, the map $F$ is said to be $\left( R{\cal{D}}_{2}^{R},\mu ^{R}\right) $-totally geodesic if $\left( \nabla F_{\ast }\right) 
(RV,X)=0$ for $V\in \Gamma \left( {\cal{D}}_{2}^{R}\right) $ and $X\in \Gamma 
\left( \mu ^{R}\right) $. 
\end{defn-new}
\section{$h$-conformal semi-invariant submersion}

In this section, we study the $h$-conformal semi-invariant submersions and almost $
h $-conformal semi-invariant submersions from almost quaternionic K\"ahler 
manifolds onto Riemannian manifolds and investigate their geometric 
properties.

\begin{defn-new}
Let $(M,E,g_{M})$ be an almost quaternionic Hermitian manifold and $(N,g_{N})$ a Riemannian manifold. An $h$-conformal submersion $F: (M,E,g_{M})\mapsto (N,g_{N})$ is called an $h$-conformal semi-invariant submersion if given a point $p \in M$ with a neighborhood $U,$ there exists a quaternionic Hermitian basis $\{J_{1},J_{2},J_{3}\}$ of sections of $E$ on $U$ such that for any $R \in \{J_{1},J_{2},J_{3}\},$ there is a distribution ${\mathcal{D}}_{1} \subset (\operatorname{ker}F_{*})$ on $U$ such that 
$$ \operatorname{ker}F_{*} = {\mathcal{D}}_{1} \oplus{\mathcal{D}}_{2}, 
R({\mathcal{D}}_{1})= {\mathcal{D}}_{1},
R({\mathcal{D}}_{2}) \subset (\operatorname{ker}F_{*})^{\perp},$$
where ${\mathcal{D}}_{2}$ is orthogonal complement of ${\mathcal{D}}_{1}$ in $\operatorname{ker}F_{*}.$ We call such a basis $\{J_{1},J_{2},J_{3}\}$ an  $h$-conformal semi-invariant basis.
\end{defn-new}

\begin{defn-new} Let $(M,E,g_{M})$ be an almost quaternionic Hermitian manifold and $(N,g_{N})$ be a Riemannian manifold. An $h$-conformal submersion $F: (M,E,g_{M})\mapsto (N,g_{N})$ is called an almost $h$-conformal semi-invariant submersion if given a point $p \in M$ with a neighborhood $U,$ there exists a quaternionic Hermitian basis $\{J_{1},J_{2},J_{3}\}$ of sections of $E$ on $U$ such that for any $R \in \{J_{1},J_{2},J_{3}\},$ there is a distribution ${\mathcal{D}}^{R}_{1} \subset (\operatorname{ker}F_{*})$ on $U$ such that 
\[ \operatorname{ker}F_{*} = {\mathcal{D}}^{R}_{1} \oplus{\mathcal{D}}^{R}_{2}, \vspace{0.4cm}  R({\mathcal{D}}^{R}_{1})= {\mathcal{D}}^{R}_{1},\vspace{0.4cm}  R({\mathcal{D}}^{R}_{2}) \subset (\operatorname{ker}F_{*})^{\perp},\]  
where ${\mathcal{D}}^{R}_{2}$ is orthogonal complement of ${\mathcal{D}}^{R}_{1}$ in $\operatorname{ker}F_{*}.$ We call such a basis $\{J_{1},J_{2},J_{3}\}$ an almost $h$-conformal semi-invariant basis.
\end{defn-new}
Let $F$ be an almost $h$-conformal semi-invariant submersion with an almost $h$-conformal semi-invariant basis $
\{J_{1},J_{2},J_{3}\}.$ 

Denote the orthogonal complement of $J_{\alpha}{{\mathcal{D}}_{2}}^{J_{\alpha}} \in \Gamma
 \left( {\operatorname{ker}}F_{\ast }\right) ^{\perp } $ by $\mu^{J_{\alpha}}$ for $\alpha \in \{1,2,3\}.$ So
$$
\left( { \operatorname{ker} } F_*\right)^{\perp}=J_{\alpha} {{\mathcal{D}}_2}^{J_{\alpha}} \oplus \mu^{J_{\alpha}} \quad \text { for } \alpha \in\{1, 2, 3\}
$$ Clearly, $\mu^{J_{\alpha}}$ is $J_{\alpha}$-invariant for $\alpha \in \{1,2,3\}.$
 
 \noindent Given $X \in \Gamma\left(\operatorname{ker} F_*\right)$, we write
$$
J_{\alpha} X=\phi_{J_{\alpha}} X+\omega_{J_{\alpha}} X,
$$
where $\phi_{J_{\alpha}} X \in \Gamma\left({{\mathcal{D}}_1}^{J_{\alpha}}\right)$ and $\omega_{J_{\alpha}} X \in \Gamma\left(J_{\alpha} {{\mathcal{D}}_2}^{J_{\alpha}}\right)$ for $\alpha \in\{1, 2, 3\}$.

\noindent Given $Z \in \Gamma\left(\left(\operatorname{ker} F_*\right)^{\perp}\right)$, we get
$$
J_{\alpha} Z=B_{J_{\alpha}} Z+C_{J_{\alpha}} Z,
$$
where $B_{J_{\alpha}} Z \in \Gamma\left({{\mathcal{D}}_2}^{J_{\alpha}}\right)$ and $C_{J_{\alpha}} Z \in \Gamma\left(\mu^{J_{\alpha}}\right)$ for $\alpha \in\{1, 2, 3\}$.
Therefore 
$$
g_M\left(C_{J_{\alpha}} X, J_{\alpha} V\right)=0
$$
for $X \in \Gamma\left(\left(\operatorname{ker} F_*\right)^{\perp}\right)$ and $V \in \Gamma\left({{\mathcal{D}}_2}^{J_{\alpha}}\right)$.

We can easily find that 
\begin{eqnarray}
    \left(\nabla_X \phi_{J_{\alpha}}\right) Y:=\widehat{\nabla}_X \phi_{J_{\alpha}} Y-\phi_{J_{\alpha}} \widehat{\nabla}_X Y         
\label{3.5}
\end{eqnarray}
and
\begin{eqnarray}
\left(\nabla_X \omega_{J_{\alpha}}\right) Y:={\mathcal{H}} \nabla_X \omega_{J_{\alpha}} Y-\omega_{J_{\alpha}} \widehat{\nabla}_X Y
\label{3.6}\end{eqnarray}
for $X, Y \in \Gamma\left(\operatorname{ker} F_{\ast}\right)$ and $\alpha \in\{1,2,3\}$.

Now, we have the following result:

\begin{lem}\label{Lemma 1}
Let $F$ be an almost h-conformal semi-invariant submersion from a quaternionic K\"ahler manifold $\left(M ,J_{1}, J_{2}, J_{3}, g_M\right)$ onto a Riemannian manifold $\left(N, g_N\right)$ such that $\{J_{1}, J_{2}, J_{3}\}$ is an almost h-conformal semi-invariant basis. Then 
\begin{itemize}
\item[\rm(1)]    For $ X, Y \in \Gamma\left(\operatorname{ker} F_*\right)$ 
$$
\begin{aligned}
& \widehat{\nabla}_X \phi_{J_{1}} Y+{\mathcal{T}}_{X} \omega_{J_{1}}Y=\phi_{J_{1}} \widehat{\nabla}_X Y+B_{J_{1}} {\mathcal{T}}_X Y + \omega_{3}(X) \phi_{J_{2}}Y- \omega_{2}(X) \phi_{J_{3}}Y,\\
& {\mathcal{T}}_X \phi_{J_{1}} Y+{\mathcal{H}} \nabla_X \omega_{J_{1}} Y=\omega_{J_{1}} \widehat{\nabla}_X Y+C_{J_{1}} {\mathcal{T}}_X Y + \omega_{3}(X) \omega _{J_{2}}Y- \omega_{2}(X) \omega_{J_{3}}Y.
\end{aligned}
$$

$$
\begin{aligned}
& \widehat{\nabla}_X \phi_{J_{2}} Y+{\mathcal{T}}_{X} \omega_{J_{2}}Y=\phi_{J_{2}} \widehat{\nabla}_X Y+B_{J_{2}} {\mathcal{T}}_X Y + \omega_{1}(X) \phi_{J_{3}}Y- \omega_{3}(X) \phi_{J_{1}}Y,\\
& {\mathcal{T}}_X \phi_{J_{2}} Y+{\mathcal{H}} \nabla_X \omega_{J_{2}} Y=\omega_{J_{2}} \widehat{\nabla}_X Y+C_{J_{2}} {\mathcal{T}}_X Y + \omega_{1}(X) \omega _{J_{3}}Y- \omega_{3}(X) \omega_{J_{1}}Y.
\end{aligned}
$$

$$
\begin{aligned}
& \widehat{\nabla}_X \phi_{J_{3}} Y+{\mathcal{T}}_{X} \omega_{J_{3}}Y=\phi_{J_{3}} \widehat{\nabla}_X Y+B_{J_{3}} {\mathcal{T}}_X Y + \omega_{2}(X) \phi_{J_{1}}Y- \omega_{1}(X) \phi_{J_{2}}Y,\\
& {\mathcal{T}}_X \phi_{J_{3}} Y+{\mathcal{H}} \nabla_X \omega_{J_{3}} Y=\omega_{J_{3}} \widehat{\nabla}_X Y+C_{J_{3}} {\mathcal{T}}_X Y + \omega_{2}(X) \omega _{J_{1}}Y- \omega_{1}(X) \omega_{J_{2}}Y.
\end{aligned}
$$

\item[\rm (2)]  For $ Z, W \in \Gamma\left(\left(\operatorname{ker} F_*\right)^{\perp}\right)$ 
$$
\begin{aligned}
& {\mathcal{V}}{\nabla}_Z B_{J_{1}} W+{\mathcal{A}}_Z C_{J_{1}} W=\phi_{J_{1}} {\mathcal{A}}_Z W+B_{J_{1}} {\mathcal{H}} \nabla_Z W + \omega_{3}(Z)B_{J_{2}}W- \omega_{2}(Z)B_{J_{3}}W, \\
& {\mathcal{A}}_Z B_{J_{1}} W+{\mathcal{H}} \nabla_Z C_{J_{1}} W=\omega_{J_{1}} {\mathcal{A}}_Z W+C_{J_{1}} {\mathcal{H}} \nabla_Z W + \omega_{3}(Z)C_{J_{2}}W- \omega_{2}(Z)C_{J_{3}}W.
\end{aligned}
$$

$$
\begin{aligned}
& {\mathcal{V}}{\nabla}_Z B_{J_{2}} W+{\mathcal{A}}_Z C_{J_{2}} W=\phi_{J_{2}} {\mathcal{A}}_Z W+B_{J_{2}} {\mathcal{H}} \nabla_Z W + \omega_{1}(Z)B_{J_{3}}W- \omega_{3}(Z)B_{J_{1}}W, \\
& {\mathcal{A}}_Z B_{J_{2}} W+{\mathcal{H}} \nabla_Z C_{J_{2}} W=\omega_{J_{2}} {\mathcal{A}}_Z W+C_{J_{2}} {\mathcal{H}} \nabla_Z W + \omega_{1}(Z)C_{J_{3}}W- \omega_{3}(Z)C_{J_{1}}W.
\end{aligned}
$$

$$
\begin{aligned}
& {\mathcal{V}}{\nabla}_Z B_{J_{3}} W+{\mathcal{A}}_Z C_{J_{3}} W=\phi_{J_{3}} {\mathcal{A}}_Z W+B_{J_{3}} {\mathcal{H}} \nabla_Z W + \omega_{2}(Z)B_{J_{1}}W- \omega_{1}(Z)B_{J_{2}}W, \\
& {\mathcal{A}}_Z B_{J_{3}} W+{\mathcal{H}} \nabla_Z C_{J_{3}} W=\omega_{J_{3}} {\mathcal{A}}_Z W+C_{J_{3}} {\mathcal{H}} \nabla_Z W + \omega_{2}(Z)C_{J_{1}}W- \omega_{1}(Z)C_{J_{2}}W.
\end{aligned}
$$

\item[\rm(3)] For $X \in \Gamma\left(\operatorname{ker} F_*\right), Z \in \Gamma\left(\left(\operatorname{ker} F_*\right)^{\perp}\right)$
$$
\begin{aligned}
& \widehat{\nabla}_X B_{J_{1}} Z+{\mathcal{T}}_X C_{J_{1}} Z=\phi_{J_{1}} {\mathcal{T}}_X Z+B_{J_{1}} {\mathcal{H}} \nabla_X Z + \omega_{3}(X)B_{J_{2}}Z- \omega_{2}(X)B_{J_{3}}Z, \\
& {\mathcal{T}}_X B_{J_{1}} Z+{\mathcal{H}} \nabla_X C_{J_{1}} Z=\omega_{J_{1}} {\mathcal{T}}_X Z+C_{J_{1}} {\mathcal{H}} \nabla_X Z + \omega_{3}(X) C_{J_{2}}Z - \omega_{2}(X)C_{J_{3}}Z.
\end{aligned}
$$

$$
\begin{aligned}
& \widehat{\nabla}_X B_{J_{2}} Z+{\mathcal{T}}_X C_{J_{2}} Z=\phi_{J_{2}} {\mathcal{T}}_X Z+B_{J_{2}} {\mathcal{H}} \nabla_X Z + \omega_{1}(X)B_{J_{3}}Z- \omega_{3}(X)B_{J_{1}}Z, \\
& {\mathcal{T}}_X B_{J_{2}} Z+{\mathcal{H}} \nabla_X C_{J_{2}} Z=\omega_{J_{2}} {\mathcal{T}}_X Z+C_{J_{2}} {\mathcal{H}} \nabla_X Z + \omega_{1}(X) C_{J_{3}}Z - \omega_{3}(X)C_{J_{1}}Z.
\end{aligned}
$$

$$
\begin{aligned}
& \widehat{\nabla}_X B_{J_{3}} Z+{\mathcal{T}}_X C_{J_{3}} Z=\phi_{J_{3}} {\mathcal{T}}_X Z+B_{J_{3}} {\mathcal{H}} \nabla_X Z + \omega_{2}(X)B_{J_{1}}Z- \omega_{1}(X)B_{J_{2}}Z, \\
& {\mathcal{T}}_X B_{J_{3}} Z+{\mathcal{H}} \nabla_X C_{J_{3}} Z=\omega_{J_{3}} {\mathcal{T}}_X Z+C_{J_{3}} {\mathcal{H}} \nabla_X Z + \omega_{2}(X) C_{J_{1}}Z - \omega_{1}(X)C_{J_{2}}Z.
\end{aligned}
$$
\end{itemize}
\end{lem}

\textbf{Proof:} (1) For $ X,Y \in \Gamma\left(\operatorname{ker} F_*\right) $, we have  
\begin{eqnarray*}
    \nabla_{X}J_{\alpha}Y &=& \nabla_{X}\phi_{J_{\alpha}}Y + \nabla_{X}\omega_{J_{\alpha}}Y, \\
    \left(\nabla_{X}J_{\alpha}\right)Y + J_{\alpha}\left(\nabla_{X}Y\right) 
    &=& {\mathcal{T}}_{X} \phi_{J_{\alpha}}Y + \widehat {\nabla}_{X}\phi_{J_{\alpha}}Y \\
    && + {\mathcal{H}}\nabla_{X}\omega_{J_{\alpha}}Y + {\mathcal{T}}_{X}\omega_{J_{\alpha}}Y, \\
    \left(\nabla_{X}J_{\alpha}\right)Y + B_{J_{\alpha}}{\mathcal{T}}_{X}Y + C_{J_{\alpha}}{\mathcal{T}}_{X}Y 
    &+& \phi_{J_{\alpha}}\widehat{\nabla}_{X}Y + \omega_{J_{\alpha}}\widehat{\nabla}_{X}Y \\
    &=& {\mathcal{T}}_{X}\phi_{J_{\alpha}}Y + \widehat {\nabla}_{X}\phi_{J_{\alpha}}Y \\
    && + {\mathcal{H}}\nabla_{X}\omega_{J_{\alpha}}Y + {\mathcal{T}}_{X}\omega_{J_{\alpha}}Y.
\end{eqnarray*}

For $\alpha ={1}$
\begin{eqnarray*}
\left(\nabla_{X}J_{1}\right)Y + B_{J_{1}}{\mathcal{T}}_{X}Y + C_{J_{1}}{\mathcal{T}}_{X}Y 
    &+& \phi_{J_{1}}\widehat{\nabla}_{X}Y + \omega_{J_{1}}\widehat{\nabla}_{X}Y \\
    &=& {\mathcal{T}}_{X}\phi_{J_{1}}Y + \widehat {\nabla}_{X}\phi_{J_{1}}Y \\
    && + {\mathcal{H}}\nabla_{X}\omega_{J_{1}}Y + {\mathcal{T}}_{X}\omega_{J_{1}}Y, \\
    \omega_{3}(X)J_{2}Y - \omega_{2}(X)J_{3}Y + B_{J_{1}}{\mathcal{T}}_{X}Y 
    &+& C_{J_{1}}{\mathcal{T}}_{X}Y + \phi_{J_{1}}\widehat{\nabla}_{X}Y + \omega_{J_{1}}\widehat{\nabla}_{X}Y \\
    &=& {\mathcal{T}}_{X}\phi_{J_{1}}Y + \widehat {\nabla}_{X}\phi_{J_{1}}Y \\
    && + {\mathcal{H}}\nabla_{X}\omega_{J_{1}}Y + {\mathcal{T}}_{X}\omega_{J_{1}}Y, \\
    \omega_{3}(X)\left(\phi_{J_{2}}Y+\omega_{J_{2}}Y\right) - \omega_{2}(X)\left(\phi_{J_{3}}Y + \omega_{J_{3}}Y\right) 
    &+& B_{J_{1}}{\mathcal{T}}_{X}Y + C_{J_{1}}{\mathcal{T}}_{X}Y \\
    && + \phi_{J_{1}}\widehat{\nabla}_{X}Y + \omega_{J_{1}}\widehat{\nabla}_{X}Y \\
    &=& {\mathcal{T}}_{X}\phi_{J_{1}}Y + \widehat {\nabla}_{X}\phi_{J_{1}}Y \\
    && + {\mathcal{H}}\nabla_{X}\omega_{J_{1}}Y + {\mathcal{T}}_{X}\omega_{J_{1}}Y.
    \end{eqnarray*}

Therefore we have,
$$ 
\begin{aligned}
& \widehat{\nabla}_X \phi_{J_{1}} Y+{\mathcal{T}}_{X} \omega_{J_{1}}Y=\phi_{J_{1}} \widehat{\nabla}_X Y+B_{J_{1}} {\mathcal{T}}_X Y + \omega_{3}(X) \phi_{J_{2}}Y- \omega_{2}(X) \phi_{J_{3}}Y,\\
& {\mathcal{T}}_X \phi_{J_{1}} Y+{\mathcal{H}} \nabla_X \omega_{J_{1}} Y=\omega_{J_{1}} \widehat{\nabla}_X Y+C_{J_{1}} {\mathcal{T}}_X Y + \omega_{3}(X) \omega _{J_{2}}Y- \omega_{2}(X) \omega_{J_{3}}Y.
\end{aligned}
$$
Similarly, we can find others.\\

(2) For $ Z,W \in \Gamma\left(\left(\operatorname{ker} F_*\right)^{\perp}\right) $
\begin{eqnarray*}
    \nabla_{Z}J_{\alpha}W &=& \nabla_{Z}\left(B_{J_{\alpha}}W+C_{J_{\alpha}}W\right),\\
    &=& \nabla_{Z}B_{J_{\alpha}}W + \nabla_{Z}C_{J_{\alpha}}W,\\
    \left( \nabla_{Z}J_{\alpha} \right)W + J_{\alpha} \left(\nabla_{Z}W\right) &=& {\mathcal{A}}_{Z}B_{J_{\alpha}}W + {\mathcal{V}}{\nabla}_{Z}B_{J_{\alpha}}W + {\mathcal{H}}\nabla_{Z}C_{J_{\alpha}}W + {\mathcal{A}}_{Z}C_{J_{\alpha}}W.
\end{eqnarray*}
For $\alpha= {1}$, we have
\begin{align*}
    \left( \nabla_{Z}J_{1} \right)W + J_{1} \left(\nabla_{Z}W\right) &= {\mathcal{A}}_{Z}B_{J_{1}}W + {\mathcal{V}}{\nabla}_{Z}B_{J_{1}}W + {\mathcal{H}}\nabla_{Z}C_{J_{1}}W + {\mathcal{A}}_{Z}C_{J_{1}}W,  
\end{align*}

\begin{align*}
    \omega_{3}(Z)J_{2}W - \omega_{2}(Z)J_{3}W 
    &+ B_{J_{1}}{\mathcal{H}}\nabla_{Z}W + C_{J_{1}}{\mathcal{H}}\nabla_{Z}W \\
    &+ \phi_{J_{1}}{\mathcal{A}}_{Z}W + \omega_{J_{1}}{\mathcal{A}}_{Z}W \\
    &= {\mathcal{A}}_{Z}B_{J_{1}}W + {\mathcal{V}}{\nabla}_{Z}B_{J_{1}}W \\
    &+ {\mathcal{H}}\nabla_{Z}C_{J_{1}}W + {\mathcal{A}}_{Z}C_{J_{1}}W.
\end{align*}

\begin{align*}
    \omega_{3}(Z)\left( B_{J_{2}}W + C_{J_{2}}W\right) 
    & - \omega_{2}(Z)\left(B_{J_{3}}W + C_{J_{3}}W\right) \\
    & + B_{J_{1}}{\mathcal{H}}\nabla_{Z}W + C_{J_{1}}{\mathcal{H}}\nabla_{Z}W \\
    & + \phi_{J_{1}}{\mathcal{A}}_{Z}W + \omega_{J_{1}}{\mathcal{A}}_{Z}W \\
    &= {\mathcal{A}}_{Z}B_{J_{1}}W + {\mathcal{V}}{\nabla}_{Z}B_{J_{1}}W \\
    & + {\mathcal{H}}\nabla_{Z}C_{J_{1}}W + {\mathcal{A}}_{Z}C_{J_{1}}W.
\end{align*}

  Therefore we have,
  $$
\begin{aligned}
& {\mathcal{V}}{\nabla}_Z B_{J_{1}} W+{\mathcal{A}}_Z C_{J_{1}} W=\phi_{J_{1}} {\mathcal{A}}_Z W+B_{J_{1}} {\mathcal{H}} \nabla_Z W + \omega_{3}(Z)B_{J_{2}}W- \omega_{2}(Z)B_{J_{3}}W,\\
& {\mathcal{A}}_Z B_{J_{1}} W+{\mathcal{H}} \nabla_Z C_{J_{1}} W=\omega_{J_{1}} {\mathcal{A}}_Z W+C_{J_{1}} {\mathcal{H}} \nabla_Z W + \omega_{3}(Z)C_{J_{2}}W- \omega_{2}(Z)C_{J_{3}}W.
\end{aligned}
$$
Similarly, we can find the rest two conditions.\\

(3) For $ X \in \Gamma\left(\operatorname{ker} F_*\right) $ and $ Z \in \Gamma\left(\left(\operatorname{ker} F_*\right)^{\perp}\right) $ 

\begin{eqnarray*}
    \nabla_{X}J_{\alpha}Z&=& \nabla_{X}B_{J_{\alpha}}Z + \nabla_{X}C_{J_{\alpha}}Z,\\
    \left(\nabla_{X}J_{\alpha}\right)Z + J_{\alpha} \left(\nabla_{X}Z\right) &=& \widehat{\nabla}_{X}B_{J_{\alpha}}Z + {\mathcal{T}}_{X}B_{J_{\alpha}}Z+ {\mathcal{H}}\nabla_{X}C_{J_{\alpha}}Z + {\mathcal{T}}_{X}C_{J_{\alpha}}Z,\\
    \left(\nabla_{X}J_{\alpha}\right)Z + J_{\alpha} \left({\mathcal{H}}\nabla_{X}Z+ {\mathcal{T}}_{X}Z\right) &=& \widehat{\nabla}_{X}B_{J_{\alpha}}Z + {\mathcal{T}}_{X}B_{J_{\alpha}}Z+ {\mathcal{H}}\nabla_{X}C_{J_{\alpha}}Z + {\mathcal{T}}_{X}C_{J_{\alpha}}Z.
\end{eqnarray*}
For $\alpha={1}$, we get
\begin{align*}
    \left(\nabla_{X}J_{1}\right)Z 
    + J_{1} \left({\mathcal{H}}\nabla_{X}Z + {\mathcal{T}}_{X}Z\right)
    &= \widehat{\nabla}_{X}B_{J_{1}}Z 
    + {\mathcal{T}}_{X}B_{J_{1}}Z 
    + {\mathcal{H}}\nabla_{X}C_{J_{1}}Z 
    + {\mathcal{T}}_{X}C_{J_{1}}Z.
\end{align*}

\begin{align*}
    \left( \nabla_{X}J_{1} \right) Z
    &+ \omega_{J_{1}}{\mathcal{T}}_{X}Z + \phi_{J_{1}}{\mathcal{T}}_{X}Z 
    + B_{J_{1}}{\mathcal{H}}\nabla_{X}Z + C_{J_{1}}{\mathcal{H}}\nabla_{X}Z \\
    &= \widehat{\nabla}_{X}B_{J_{1}}Z + {\mathcal{T}}_{X}B_{J_{1}}Z 
    + {\mathcal{H}}\nabla_{X}C_{J_{1}}Z + {\mathcal{T}}_{X}C_{J_{1}}Z.
\end{align*}

\begin{align*}
    \omega_{3}(X)J_{2}Z - \omega_{2}(X)J_{3}Z 
    &\quad + \omega_{J_{1}}{\mathcal{T}}_{X}Z + \phi_{J_{1}}{\mathcal{T}}_{X}Z \\
    &\quad + B_{J_{1}}{\mathcal{H}}\nabla_{X}Z + C_{J_{1}}{\mathcal{H}}\nabla_{X}Z \\
    &= \widehat{\nabla}_{X}B_{J_{1}}Z + {\mathcal{T}}_{X}B_{J_{1}}Z + {\mathcal{H}}\nabla_{X}C_{J_{1}}Z + {\mathcal{T}}_{X}C_{J_{1}}Z.
\end{align*}

\begin{align*}
    \omega_{3}(X)B_{J_{2}}Z + \omega_{3}(X)C_{J_{2}}Z 
    &- \omega_{2}(X)B_{J_{3}}Z - \omega_{2}(X)C_{J_{3}}Z \\
    &+ \omega_{J_{1}}{\mathcal{T}}_{X}Z + \phi_{J_{1}}{\mathcal{T}}_{X}Z \\
    &+ B_{J_{1}}{\mathcal{H}}\nabla_{X}Z + C_{J_{1}}{\mathcal{H}}\nabla_{X}Z \\
    &= \widehat{\nabla}_{X}B_{J_{1}}Z + {\mathcal{T}}_{X}B_{J_{1}}Z \\
    &\quad + {\mathcal{H}}\nabla_{X}C_{J_{1}}Z + {\mathcal{T}}_{X}C_{J_{1}}Z.
\end{align*}
Therefore we have,
$$
\begin{aligned}
& \widehat{\nabla}_X B_{J_{1}} Z+{\mathcal{T}}_X C_{J_{1}} Z=\phi_{J_{1}} {\mathcal{T}}_X Z+B_{J_{1}} {\mathcal{H}} \nabla_X Z + \omega_{3}(X)B_{J_{2}}Z- \omega_{2}(X)B_{J_{3}}Z, \\
& {\mathcal{T}}_X B_{J_{1}} Z+{\mathcal{H}} \nabla_X C_{J_{1}} Z=\omega_{J_{1}} {\mathcal{T}}_X Z+C_{J_{1}} {\mathcal{H}} \nabla_X Z + \omega_{3}(X) C_{J_{2}}Z - \omega_{2}(X)C_{J_{3}}Z.
\end{aligned}
$$
Similarly, we can find others.

\begin{rem} By {\rm (\ref{3.5}), (\ref{3.6}),} and Lemma{\rm~\ref{Lemma 1}(1)} , we have\\
For $J_{1}$
$$
\begin{aligned}
& \left(\nabla_X \phi_{J_{1}}\right) Y=B_{J_{1}} {\mathcal{T}}_X Y-{\mathcal{T}}_X \omega_{J_{1}} Y + \omega_{3}(X) \phi_{J_{2}}Y - \omega_{2}(X) \phi_{J_{3}}Y, \\
& \left(\nabla_X \omega_{J_{1}}\right) Y=C_{J_{1}} {\mathcal{T}}_X Y-{\mathcal{T}}_X \phi_{J_{1}} Y + \omega_{3}(X) \omega_{J_{2}}Y- \omega_{2}(X)\omega_{J_{3}}Y.
\end{aligned}
$$
For $J_{2}$
$$
\begin{aligned}
& \left(\nabla_X \phi_{J_{2}}\right) Y=B_{J_{2}} {\mathcal{T}}_X Y-{\mathcal{T}}_X \omega_{J_{2}} Y + \omega_{1}(X) \phi_{J_{3}}Y - \omega_{3}(X) \phi_{J_{1}}Y, \\
& \left(\nabla_X \omega_{J_{2}}\right) Y=C_{J_{2}} {\mathcal{T}}_X Y-{\mathcal{T}}_X \phi_{J_{2}} Y + \omega_{1}(X) \omega_{J_{3}}Y- \omega_{3}(X)\omega_{J_{1}}Y.
\end{aligned}
$$
For $J_{3}$
$$
\begin{aligned}
& \left(\nabla_X \phi_{J_{3}}\right) Y=B_{J_{3}} {\mathcal{T}}_X Y-{\mathcal{T}}_X \omega_{J_{3}} Y + \omega_{2}(X) \phi_{J_{1}}Y - \omega_{1}(X) \phi_{J_{2}}Y, \\
& \left(\nabla_X \omega_{J_{3}}\right) Y=C_{J_{3}} {\mathcal{T}}_X Y-{\mathcal{T}}_X \phi_{J_{3}} Y + \omega_{2}(X) \omega_{J_{1}}Y- \omega_{1}(X)\omega_{J_{2}}Y.
\end{aligned}
$$
for $X, Y \in \Gamma\left(\right. \operatorname{ker} \left.F_{\ast}\right)$.

\end{rem}
Now, we investigate the integrability of some distributions.

\begin{lem}
 Let $F$ be an h-conformal semi-invariant submersion from a quaternionic K\"ahler manifold $\left(M, I, J, K, g_M\right)$ onto a Riemanian manifold $\left(N, g_N\right)$ such that $\{J_{1}, J_{2}, J_{3}\}$ is an $h$-confornal semi-invariant basis. Then we have the followings:\newline
(i) The distribution ${\mathcal{D}}_2$ is always integrable.\newline
(ii) The following conditions are equivalent:
\vspace{-0.2cm}
\begin{itemize}
\item[(a)]The distribution  ${\mathcal{D}}_1$ is integrable.
\vspace{-0.2cm}
\item[(b)]$\left(\nabla F_*\right)(Z, J_{1} Y)-\left(\nabla F_*\right)(Y, J_{1} Z) \in \Gamma\left(F_* \mu^{J_{1}}\right)$  for  $Y, Z \in \Gamma\left({\mathcal{D}}_1\right)$.
\vspace{-0.2cm}
\item[(c)] $\left(\nabla F_*\right)(Z, J_{2} Y)-\left(\nabla F_*\right)(Y, J_{2} Z) \in \Gamma\left(F_* \mu^{J_{2}}\right)$  for  $ Y, Z \in \Gamma\left({\mathcal{D}}_1\right)$.
\vspace{-0.2cm}
\item[ (d)] $\left(\nabla F_*\right)(Z, J_{3} Y)-\left(\nabla F_*\right)(Y, J_{3} Z) \in \Gamma\left(F_* \mu^{J_{3}}\right)$  for  $Y, Z \in \Gamma\left({\mathcal{D}}_1\right)$.
\end{itemize}

\end{lem}
 \noindent {\bf Proof:} By result \ref{2.10}, we have $[Y, Z] \in \Gamma\left(\operatorname{ker} F_*\right)$ for $Y, Z \in \Gamma\left(\operatorname{ker} F_*\right)$.\newline
We claim that ${\mathcal{T}}_Y J_{\alpha} Z={\mathcal{T}}_Z J_{\alpha} Y$ for $Y, Z \in \Gamma\left({\mathcal{D}}_2\right)$ and $ \alpha\in\{1, 2, 3\}$.\newline
Let  $X \in \Gamma\left(\operatorname{ker} F_*\right)$ and $\alpha=1$. Using the skew-symmetry property of ${\mathcal T}_Y$, we get
$$
\begin{aligned}
g_M\left({\mathcal{T}}_Y J_{1} Z, X\right) & =-g_M\left(J_{1} Z, {\mathcal{T}}_{Y}X\right)=-g_M\left(J_{1} Z, \nabla_{Y} X - \widehat{\nabla}_{Y}X\right)= -g_{M}\left( J_{1}Z, \nabla_{Y}X\right)\\ &=-g_M\left( J_{1} Z, \nabla_{X}Y\right) 
 =g_M\left(\nabla_{X}J_{1}Z,Y\right)=g_M\left(J_{1} \left(\nabla_{X}Z\right) + \left( \nabla_{X}J_{1} \right)Z,  Y\right)\\ &=g_M \left( J_{1}\left(\nabla_{X}Z\right) , Y\right) +g_M \left( \left( \nabla_{X}J_{1}\right) Z,Y \right) \\ &= -g_M \left( \nabla_{X} Z, J_{1}Y\right) + g_M \left( \omega_{3}(X)J_{2}Z- \omega_{2}(X)J_{3}Z,Y   \right)\\ & = 
  -g_M \left( \nabla_{Z} X, J_{1}Y\right) + g_M \left( \omega_{3}(X)J_{2}Z- \omega_{2}(X)J_{3}Z,Y   \right)\\ & = g_M \left( X, \nabla_{Z} J_{1}Y\right) + g_M \left( \omega_{3}(X)\left( \phi_{J_{2}}Z+ \omega_{J_{2}}Z\right)- \omega_{2}(X) \left( \phi_{J_{3}}Z + \omega_{J_{3}}Z\right),Y   \right)\\ & =
  g_M \left( X, \nabla_{Z} J_{1}Y\right) + g_M \left( \omega_{3}(X)\phi_{J_{2}}Z + \omega_{3}(X)\omega_{J_{2}}Z- \omega_{2}(X)  \phi_{J_{3}}Z - \omega_{2}(X) \omega_{J_{3}}Z,Y \right)\\ & =
  g_M \left( X, {\mathcal{H}}\nabla_{Z} J_{1}Y + {\mathcal{T}}_{Z}J_{1}Y\right) + g_M \left( \omega_{3}(X)\phi_{J_{2}}Z , Y \right) + g_M \left(\omega_{3}(X)\omega_{J_{2}}Z , Y \right)\\ &- g_M \left( \omega_{2}(X)  \phi_{J_{3}}Z,Y \right) -g_M \left(\omega_{2}(X) \omega_{J_{3}}Z,Y \right)\\ & =
  g_M \left( X, {\mathcal{H}}\nabla_{Z} J_{1}Y \right) + g_M \left({X, \mathcal{T}}_{Z}J_{1}Y\right)\\ &
 =g_M\left(X, {\mathcal{T}}_Z J_{1} Y\right),
\end{aligned}
$$
which means our claim.

Given $ Y,Z \in \Gamma \left( {\mathcal{D}}_{2}\right)$ and $X \in \Gamma \left( {\mathcal{D}}_{1}\right)$, we obtain
$$ g_M\left( [Y,Z],X\right) = g_M\left( \nabla_{Y}Z-\nabla_{Z}Y,X\right) = g_M \left( {\mathcal{T}}_{Y}J_{1}Z-{\mathcal{T}}_{Z}J_{1}Y, J_{1}X\right)=0 $$
Similarly we can find for $\alpha={2}$ $ \&$ $ \alpha ={3}$,
which implies (i).\newline
\noindent
For (ii), Given $ X,Y\in \Gamma \left( {\mathcal{D}}_{1}\right)$ and $Z \in \Gamma \left( {\mathcal{D}}_{2}\right)$, and $ \alpha \in \{1,2,3\}$, we have
$$
\begin{aligned}
 g_M\left( [X,Y],Z\right) &= \frac{1}{\lambda^{2}}g_N\left( F_{\ast}\nabla_{X}J_{1}Y-F_{\ast}\nabla_{Y}J_{1}X, F_{\ast}J_{1}Z\right) \\ &= \frac{1}{\lambda^{2}} g_N\left( \left(\nabla F_{\ast}\right) \left( Y,J_{1}X \right) - \left( \nabla F_{\ast}\right) \left( X,J_{1}Y\right), F_{\ast}J_{1}Z \right)
\end{aligned}
$$
 since we are given that distribution ${\mathcal{D}}_{1}$ is integrable which means $g_M\left([X,Y],Z\right)=0$ \newline
 so that we get (a) $\Leftrightarrow $ (b), (a) $\Leftrightarrow $ (c), (a) $\Leftrightarrow $(d).
 Therefore, we get the result.
 
\begin{theorem}
\label{th-3} Let $F$ be an almost $h$-conformal semi-invariant submersion from a 
quaternionic K\"{a}hler manifold $(M,J_{1},J_{2},J_{3},g_{M})$ onto a 
Riemannian manifold $(N,g_{N})$ such that $\{J_{1},J_{2},J_{3}\}$ is an almost $h$-conformal semi-invariant basis. Then the following conditions are 
equivalent: \newline
{\rm (a)} the distribution $(\operatorname{ker}F_{\ast })^{\perp }$ is integrable,\newline
{\rm (b)} ${\mathcal{A}}_{Y}\omega_{J_{1}}B_{J_{1}}X-{\mathcal{A}}_{X}\omega_{J_{1}}B_{J_{1}}Y + \phi_{J_{1}}\left( {\mathcal{A}}_{Y}C_{J_{1}}X- {\mathcal{A}}_{X}C_{J_{1}}Y\right) \in \Gamma \left({\mathcal{D}}_{2}^{J_{1}} \right) $ and 
\begin{eqnarray*}
\lambda ^{-2}g_{N}\left( \nabla _{Y}^{N}F_{\ast }C_{J_{1}}X-\nabla
_{X}^{N}F_{\ast }C_{J_{1}}Y,F_{\ast }J_{1}V\right) &=&g_{M}\left( {\cal A}
_{X}{B}_{J_{1}}Y-{\cal A}_{Y}{B}_{{J}_{1}}X,J_{1}V\right) \\
&-&g_{M}\left( {\cal H}grad\ln \lambda ,{C}_{J_{1}}Y\right)
g_{M}(X,J_{1}V) \\
&+&g_{M}\left( {\cal H}grad\ln \lambda ,{C}_{J_{1}}X\right)
g_{M}(Y,J_{1}V) \\
&-&2g_{M}\left( { C}_{J_{1}}X,Y\right) g_{M}({\cal H}grad\ln \lambda
,J_{1}V) \\
&+&\omega _{3}(Y)g_{M}\left( {C}_{J_{2}}X,J_{1}V\right) -\omega
_{2}(Y)g_{M}\left( {C}_{J_{3}}X,J_{1}V\right) \\
&-&\omega _{3}(X)g_{M}\left( { C}_{J_{2}}Y,J_{1}V\right) +\omega
_{2}(X)g_{M}\left( {C}_{J_{3}}Y,J_{1}V\right) ,
\end{eqnarray*}
for $X,Y\in \Gamma ((\operatorname{ker}F_{\ast })^{\perp })$ and $V\in \Gamma \left({\mathcal{D}}_{2}^{J_{1}}\right).$\\
{\rm (c)} ${\mathcal{A}}_{Y}\omega_{J_{2}}B_{J_{2}}X-{\mathcal{A}}_{X}\omega_{J_{2}}B_{J_{2}}Y + \phi_{J_{2}}\left( {\mathcal{A}}_{Y}C_{J_{2}}X- {\mathcal{A}}_{X}C_{J_{2}}Y\right) \in \Gamma \left({\mathcal{D}}_{2}^{J_{2}} \right) $ and
\begin{eqnarray*}
\lambda ^{-2}g_{N}\left( \nabla _{Y}^{N}F_{\ast }C_{J_{2}}X-\nabla
_{X}^{N}F_{\ast }C_{J_{2}}Y,F_{\ast }J_{2}V\right) 
&=&g_{M}\left( {\cal A}_{X}{B}_{J_{2}}Y-{\cal A}_{Y}{B}
_{J_{2}}X,J_{2}V\right) \\
&-&g_{M}\left( {\cal H}grad\ln \lambda ,{C}_{J_{2}}Y\right)
g_{M}(X,J_{2}V) \\
&+&g_{M}\left( {\cal H}grad\ln \lambda ,{C}_{J_{2}}X\right)
g_{M}(Y,J_{2}V) \\
&-&2g_{M}\left( { C}_{J_{2}}X,Y\right) g_{M}({\cal H}grad\ln \lambda
,J_{2}V) \\
&+&\omega _{1}(Y)g_{M}\left( { C}_{J_{3}}X,J_{2}V\right) -\omega
_{3}(Y)g_{M}\left( {C}_{J_{1}}X,J_{2}V\right) \\
&-&\omega _{1}(X)g_{M}\left( {C}_{J_{3}}Y,J_{2}V\right) +\omega
_{3}(X)g_{M}\left( {C}_{J_{1}}Y,J_{2}V\right) ,
\end{eqnarray*}for $X,Y\in \Gamma ((\operatorname{ker}F_{\ast })^{\perp })$ and $V\in \Gamma \left({\mathcal{D}}_{2}^{J_{2}}\right).$\\
{\rm (d)} ${\mathcal{A}}_{Y}\omega_{J_{3}}B_{J_{3}}X-{\mathcal{A}}_{X}\omega_{J_{3}}B_{J_{3}}Y + \phi_{J_{3}}\left( {\mathcal{A}}_{Y}C_{J_{3}}X- {\mathcal{A}}_{X}C_{J_{3}}Y\right) \in \Gamma \left({\mathcal{D}}_{2}^{J_{3}} \right) $  and \\
\begin{eqnarray*}
\lambda ^{-2}g_{N}\left( \nabla _{Y}^{N}F_{\ast }C_{J_{3}}X-\nabla
_{X}^{N}F_{\ast }C_{J_{3}}Y,F_{\ast }J_{3}V\right) &=&g_{M}\left( {\cal A}
_{X}{ B}_{J_{3}}Y-{\cal A}_{Y}{ B} _{J_{3}}X,J_{3}V\right) \\
&-&g_{M}\left( {\cal H} grad \ln \lambda ,{ C}_{J_{3}}Y\right)
g_{M}(X,J_{3}V) \\
&+&g_{M}\left( {\cal H} grad \ln \lambda ,{C}_{J_{3}}X\right)
g_{M}(Y,J_{3}V) \\
&-&2g_{M}\left( {C}_{J_{3}}X,Y\right) g_{M}({\cal H} grad \ln \lambda
,J_{3}V) \\
&+&\omega _{2}(Y)g_{M}\left( {C}_{J_{1}}X,J_{3}V\right) -\omega
_{1}(Y)g_{M}\left( {C}_{J_{2}}X,J_{3}V\right) \\
&-&\omega _{2}(X)g_{M}\left( { C}_{J_{1}}Y,J_{3}V\right) +\omega_{1}(X)g_{M}\left( { C}_{J_{2}}Y,J_{3}V\right).
\end{eqnarray*}
for $X,Y\in \Gamma ((\operatorname{ker}F_{\ast })^{\perp })$ and $V\in \Gamma \left({\mathcal{D}}_{2}^{J_{3}}\right)$.
\end{theorem}
\noindent {\bf Proof:} For $X,Y\in \Gamma ((\operatorname{ker}F_{\ast })^{\perp })$ and $
W\in \Gamma ({\mathcal{D}}_{1}^{J_{\alpha}})$
\begin{eqnarray*}
g_{M}\left( [X,Y],W\right) &=&g_{M}\left( J_{\alpha}[X,Y],J_{\alpha}W\right) \\
&=&g_{M}\left( J_{\alpha}\nabla _{X}Y-J_{\alpha}\nabla _{Y}X,J_{\alpha}W\right) \\
&=&g_{M}\left( \nabla _{X}J_{\alpha}Y-\nabla _{Y}J_{\alpha}X,J_{\alpha}W\right) +g_{M}\left(
\left( \nabla _{Y}J_{\alpha}\right) X-\left( \nabla _{X}J_{\alpha}\right) Y,J_{\alpha}W\right) .
\end{eqnarray*}
Now 
\begin{eqnarray*}
g_{M}\left( \nabla _{X}J_{\alpha}Y-\nabla _{Y}J_{\alpha}X,J_{\alpha}W\right) &=&g_{M}\left(
\nabla _{X}B_{J_{\alpha}}Y,J_{\alpha}W\right) +g_{M}\left( \nabla _{X}C_{J_{\alpha}}Y,J_{\alpha}W\right)
\\
&-&g_{M}\left( \nabla _{Y}B_{J_{\alpha}}X,J_{\alpha}W\right) -g_{M}\left( \nabla
_{Y}C_{J_{\alpha}}X,J_{\alpha}W\right) \\
&=&-g_{M} \left( \nabla_{X}J_{\alpha}B_{J_{\alpha}}Y,W \right)+ g_{M}\left( \nabla_{Y}J_{\alpha}B_{J_{\alpha}}X,W\right)\\
&+& g_{M}\left( {\mathcal H} \nabla_{X}C_{J_{\alpha}}Y + {\mathcal{A}}_{X}C_{J_{\alpha}}Y, J_{\alpha}W \right)\\
&-& g_{M}\left({\mathcal{H}}\nabla_{Y}C_{J_{\alpha}}X + {\mathcal{A}}_{Y}C_{J_{\alpha}}X,J_{\alpha}W \right)\\
&=& -g_{M}\left( \nabla_{X} \left( \phi_{J_{\alpha}}B_{J_{\alpha}}Y +\omega_{J_{\alpha}}B_{J_{\alpha}}Y \right), W \right) \\
&+& g_{M}\left( \nabla_{Y} \left( \phi_{J_{\alpha}}B_{J_{\alpha}}X +\omega_{J_{\alpha}}B_{J_{\alpha}}X \right), W \right) \\
&+& g_{M} \left( {\mathcal{A}}_{X}C_{J_{\alpha}}Y, J_{\alpha}W\right) -g_{M}\left({\mathcal{A}}_{Y}C_{J_{\alpha}}X,J_{\alpha}W \right)\\
&=& -g_{M}\left( \nabla_{X} \omega_{J_{\alpha}}B_{J_{\alpha}}Y,W \right) + g_{M}\left( \nabla_{Y}\omega_{J_{\alpha}}B_{J_{\alpha}}X,W \right)\\
&-& g_{M}\left(J_{\alpha}{\mathcal{A}}_{X}C_{J_{\alpha}}Y,W\right)+g_{M}\left(J_{\alpha}{\mathcal{A}}_{Y}C_{J_{\alpha}}X,W\right)\\
&=& -g_{M}\left( {\mathcal{H}}\nabla_{X} \omega_{J_{\alpha}}B_{J_{\alpha}}Y +{\mathcal{A}}_{X}\omega_{J_{\alpha}}B_{J_{\alpha}}Y ,W \right)\\
&+& g_{M}\left( {\mathcal{H}} \nabla_{Y}\omega_{J_{\alpha}}B_{J_{\alpha}}X + {\mathcal{A}}_{Y}\omega_{J_{\alpha}}B_{J_{\alpha}}X,W \right)\\
&-& g_{M}\left(\phi_{J_{\alpha}}{\mathcal{A}}_{X}C_{J_{\alpha}}Y +\omega_{J_{\alpha}}{\mathcal{A}}_{X}C_{J_{\alpha}}Y ,W\right) \\
&+& g_{M}\left(\phi_{J_{\alpha}}{\mathcal{A}}_{Y}C_{J_{\alpha}}X + \omega_{J_{\alpha}}{\mathcal{A}}_{Y}C_{J_{\alpha}}X ,W\right)\\
&=& -g_{M}\left( {\mathcal{A}}_{X}\omega_{J_{\alpha}}B_{J_{\alpha}}Y,W\right) + g_{M}\left( {\mathcal{A}}_{Y}\omega_{J_{\alpha}}B_{J_{\alpha}}X,W \right)\\
&-& g_{M}\left( \phi_{J_{\alpha}}{\mathcal{A}}_{X}C_{J_{\alpha}}Y,W\right)+ g_{M}\left( \phi_{J_{\alpha}}{\mathcal{A}}_{Y}C_{J_{\alpha}}X,W \right)\\
&=& g_{M} \left( {\mathcal{A}}_{Y}\omega_{J_{\alpha}}B_{J_{\alpha}}X - {\mathcal{A}}_{X}\omega_{J_{\alpha}}B_{J_{\alpha}}Y+ \phi_{J_{\alpha}} \left({\mathcal{A}}_{Y}C_{J_{\alpha}}X - {\mathcal{A}}_{X}C_{J_{\alpha}}Y\right),W\right).
\end{eqnarray*}
So, 
\begin{eqnarray*}
g_{M}\left( [X,Y],W\right) &=&g_{M} \left( {\mathcal{A}}_{Y}\omega_{J_{\alpha}}B_{J_{\alpha}}X - {\mathcal{A}}_{X}\omega_{J_{\alpha}}B_{J_{\alpha}}Y + \phi_{J_{\alpha}} \left({\mathcal{A}}_{Y}C_{J_{\alpha}}X - {\mathcal{A}}_{X}C_{J_{\alpha}}Y\right),W\right) \\
&+& g_{M}\left(\left( \nabla _{Y}J_{\alpha}\right) X,J_{\alpha}W\right)-g_{M}\left( \left( \nabla _{X}J_{\alpha}\right) Y,J_{\alpha}W\right).   
\end{eqnarray*}
For $J_{1}$, we have
\begin{eqnarray*}
 g_{M}\left( [X,Y],W\right) &=&g_{M} \left( {\mathcal{A}}_{Y}\omega_{J_{1}}B_{J_{1}}X - {\mathcal{A}}_{X}\omega_{J_{1}}B_{J_{1}}Y + \phi_{J_{1}} \left({\mathcal{A}}_{Y}C_{J_{1}}X - {\mathcal{A}}_{X}C_{J_{1}}Y\right),W\right) \\
&+& g_{M}\left(\left( \nabla _{Y}J_{1}\right) X,J_{1}W\right)-g_{M}\left( \left( \nabla _{X}J_{1}\right) Y,J_{1}W\right)\\
&=& g_{M} \left( {\mathcal{A}}_{Y}\omega_{J_{1}}B_{J_{1}}X - {\mathcal{A}}_{X}\omega_{J_{1}}B_{J_{1}}Y + \phi_{J_{1}} \left({\mathcal{A}}_{Y}C_{J_{1}}X - {\mathcal{A}}_{X}C_{J_{1}}Y\right),W\right)\\
&+& g_{M} \left( \omega_{3}(Y)J_{2}X, J_{1}W\right)- g_{M}\left(\omega_{2}(Y)J_{3}X,J_{1}W\right)\\
&-&g_{M}\left( \omega_{3}(X) J_{2}Y, J_{1}W\right)+ g_{M}\left( \omega_{2}(X)J_{3}Y, J_{1}W\right)\\
&=& g_{M} \left( {\mathcal{A}}_{Y}\omega_{J_{1}}B_{J_{1}}X - {\mathcal{A}}_{X}\omega_{J_{1}}B_{J_{1}}Y + \phi_{J_{1}} \left({\mathcal{A}}_{Y}C_{J_{1}}X - {\mathcal{A}}_{X}C_{J_{1}}Y\right),W\right)\\ 
&+& g_{M} \left( \omega_{3}(Y)\left(B_{J_{2}}X+C_{J_{2}}X\right), J_{1}W\right)- g_{M}\left(\omega_{2}(Y) \left(B_{J_{3}}+C_{J_{3}}X\right),J_{1}W\right)\\
&-&g_{M}\left( \omega_{3}(X) \left(B_{J_{2}}Y+ C_{J_{2}}Y\right), J_{1}W\right)+ g_{M}\left( \omega_{2}(X) \left(B_{J_{3}}Y+C_{J_{3}}Y\right), J_{1}W\right) \\
&=& g_{M} \left( {\mathcal{A}}_{Y}\omega_{J_{1}}B_{J_{1}}X - {\mathcal{A}}_{X}\omega_{J_{1}}B_{J_{1}}Y + \phi_{J_{1}} \left({\mathcal{A}}_{Y}C_{J_{1}}X - {\mathcal{A}}_{X}C_{J_{1}}Y\right),W\right).
\end{eqnarray*}
Since $ (\operatorname{ker}F_{\ast })^{\perp } $ is integrable so, $ g_{M} \left( [X,Y],W \right)=0$ for $ W \in \Gamma({D_{1}}^{J_{1}})$, so\\
$$  {\mathcal{A}}_{Y}\omega_{J_{1}}B_{J_{1}}X - {\mathcal{A}}_{X}\omega_{J_{1}}B_{J_{1}}Y + \phi_{J_{1}} \left({\mathcal{A}}_{Y}C_{J_{1}}X - {\mathcal{A}}_{X}C_{J_{1}}Y\right) \in \Gamma ({D_{2}}^{J_{1}}). $$
\newline
\noindent
Now, given that $ V \in \Gamma({D_{2}}^{J_{1}})$
\begin{eqnarray*}
    g_{M}\left( [X,Y],V \right) &=&  g_{M}\left(J_{\alpha}[X,Y],J_{\alpha}V\right)\\
    &=& g_{M}\left(J_{\alpha}\nabla_{X}Y-J_{\alpha}\nabla_{Y}X,J_{\alpha}V\right) \\
    &=& g_{M}\left(\nabla_{X}J_{\alpha}Y-(\nabla_{X}J_{\alpha})Y-\nabla_{Y}J_{\alpha}X + (\nabla_{Y}J_{\alpha})X,J_{\alpha}V\right)\\
    &=& g_{M}\left( \nabla_{X}J_{\alpha}Y-\nabla_{Y}J_{\alpha}X ,J_{\alpha}V\right)+ g_{M}\left( (\nabla_{Y}J_{\alpha})X-(\nabla_{X}J_{\alpha})Y ,J_{\alpha}V\right)\\
&=&g_{M}\left(\nabla_{X}B_{J_{\alpha}}Y,J_{\alpha}V \right)+ g_{M}\left(\nabla_{X}C_{J_{\alpha}}Y,J_{\alpha}V\right)
- g_{M}\left( \nabla_{Y}B_{J_{\alpha}}X,J_{\alpha}V\right)\\
&-&g_{M}\left(\nabla_{Y}C_{J_{\alpha}}X,J_{\alpha}V\right)+g_{M}\left( (\nabla_{Y}J_{\alpha})X-(\nabla_{X}J_{\alpha})Y ,J_{\alpha}V\right)\\
&=&g_{M}\left({\mathcal{A}}_{X}B_{J_{\alpha}}Y + {\mathcal{V}} \nabla_{X}B_{J_{\alpha}}Y,J_{\alpha}V \right) -g_{M}\left( {\mathcal{A}}_{Y}B_{J_{\alpha}}X + {\mathcal{V}}\nabla_{Y} B_{J_{\alpha}}X,J_{\alpha}V\right)\\
&+&g_{M}\left(\nabla_{X}C_{J_{\alpha}}Y,J_{\alpha}V\right) -g_{M}\left(\nabla_{Y}C_{J_{\alpha}}X,J_{\alpha}V\right)+g_{M}\left( (\nabla_{Y}J_{\alpha})X-(\nabla_{X}J_{\alpha})Y ,J_{\alpha}V\right)\\
&=& g_{M}\left({\mathcal{A}}_{X}B_{J_{\alpha}}Y -{\mathcal{A}}_{Y}B_{J_{\alpha}}X,J_{\alpha}V\right)+\lambda
^{-2}g_{N}\left( F_{\ast }\left( \nabla _{X}C_{J_{\alpha}}Y\right) ,F_{\ast
}J_{\alpha}V\right) \\
&-&\lambda ^{-2}g_{N}\left( F_{\ast }\left( \nabla _{Y}C_{J_{\alpha}}X\right)
,F_{\ast }J_{\alpha}V\right) +g_{M}\left( \left( \nabla _{Y}J_{\alpha}\right) X-\left( \nabla
_{X}J_{\alpha}\right) Y,J_{\alpha}V\right), 
\end{eqnarray*}
\noindent
using (\ref{eq-bb}) in above equation, we have
\begin{eqnarray*}
g_{M}\left( [X,Y],V\right) &=&g_{M}\left( {\cal A}_{X}B_{J_{\alpha}}Y-{\cal A}
_{Y}B_{J_{\alpha}}X,J_{\alpha}V\right) +\lambda ^{-2}g_{N}\left( \nabla _{X}F_{\ast
}C_{J_{\alpha}}Y-\left( \nabla F_{\ast }\right) \left( X,C_{J_{\alpha}}Y\right) ,F_{\ast
}J_{\alpha}V\right) \\
&&-\lambda ^{-2}g_{N}\left( \nabla _{Y}F_{\ast }C_{J_{\alpha}}X-\left( \nabla
F_{\ast }\right) \left( Y,C_{J_{\alpha}}X\right) ,F_{\ast }J_{\alpha}V\right) \\
&&+g_{M}\left( \left( \nabla _{Y}J_{\alpha}\right) X-\left( \nabla _{X}J_{\alpha}\right)
Y,J_{\alpha}V\right).
\end{eqnarray*} 
By using (\ref{eq-cc}) in the above equation, we get

\begin{eqnarray*}
g_{M}\left( [X,Y],V\right) 
&=&g_{M}\left( {\cal A}_{X}B_{J_{\alpha}}Y-{\cal A}_{Y}B_{J_{\alpha}}X,J_{\alpha}V\right) +\lambda
^{-2}g_{N}\left( \nabla _{X}F_{\ast }C_{J_{\alpha}}Y-\nabla _{Y}F_{\ast
}C_{J_{\alpha}}X,F_{\ast }J_{\alpha}V\right) \\
&+&\lambda ^{-2}g_{N}\left( Y\left( \ln \lambda \right) F_{\ast
}C_{J_{\alpha}}X+C_{J_{\alpha}} {X}\left( \ln \lambda \right) F_{\ast }Y-g_{M}\left(
Y,C_{J_{\alpha}}X\right) F_{\ast }\left( grad\ln \lambda \right) ,F_{\ast }J_{\alpha}V\right)
\\
&-&\lambda ^{-2}g_{N}\left( X\left( \ln \lambda \right) F_{\ast
}C_{J_{\alpha}}Y+C_{J_{\alpha}} {Y}\left( \ln \lambda \right) F_{\ast }X-g_{M}\left(
X,C_{J_{\alpha}}Y\right) F_{\ast }\left( grad \ln \lambda \right) ,F_{\ast }J_{\alpha}V\right)
\\
&+&g_{M}\left( \left( \nabla _{Y}J_{\alpha}\right) X-\left( \nabla _{X}J_{\alpha}\right)
Y,J_{\alpha}V\right) \\&=&g_{M}\left( {\cal A}_{X}B_{J_{\alpha}}Y-{\cal A}
_{Y}B_{J_{\alpha}}X,J_{\alpha}V\right) +\lambda ^{-2}g_{N}\left( \nabla _{X}F_{\ast
}C_{J_{\alpha}}Y-\nabla _{Y}F_{\ast }C_{J_{\alpha}}X,F_{\ast }J_{\alpha}V\right) \\
&+&g_{M}\left( {\cal H} grad \ln \lambda ,C_{J_{\alpha}}X\right) g_{M}\left(
Y,J_{\alpha}V\right) -g_{M}\left( {\cal H} grad \ln \lambda ,C_{J_{\alpha}}Y\right)
g_{M}\left( X,J_{\alpha}V\right) \\
&-&2g_{M}\left( C_{J_{\alpha}}X,Y\right) g_{M}\left( {\cal H} grad \ln \lambda
,J_{\alpha}V\right) +g_{M}\left( \left( \nabla _{Y}J_{\alpha}\right) X-\left( \nabla
_{X}J_{\alpha}\right) Y,J_{\alpha}V\right) .
\end{eqnarray*}
$\left( \operatorname{ker}F_{\ast }\right) ^{\perp }$ is integrable which means that
$ g_{M}\left( [X,Y],V\right) =0$, therefore  
\begin{eqnarray*}
0 &=&g_{M}\left( {\cal A}_{X}B_{J_{\alpha}}Y-{\cal A}_{Y}B_{J_{\alpha}}X,J_{\alpha}V\right) +\lambda
^{-2}g_{N}\left( \nabla _{X}F_{\ast }C_{J_{\alpha}}Y-\nabla _{Y}F_{\ast
}C_{J_{\alpha}}X,F_{\ast }J_{\alpha}V\right) \\
&+&g_{M}\left( {\cal H} grad \ln \lambda ,C_{J_{\alpha}}X\right) g_{M}\left(
Y,J_{\alpha}V\right) -g_{M}\left( {\cal H} grad \ln \lambda ,C_{J_{\alpha}}Y\right)
g_{M}\left( X,J_{\alpha}V\right) \\
&-&2g_{M}\left( C_{J_{\alpha}}X,Y\right) g_{M}\left( {\cal H} grad \ln \lambda
,J_{\alpha}V\right) +g_{M}\left( \left( \nabla _{Y}J_{\alpha}\right) X-\left( \nabla
_{X}J_{\alpha}\right) Y,J_{\alpha}V\right),
\end{eqnarray*}
\begin{eqnarray*}
-\lambda ^{-2}g_{N}\left( \nabla _{X}F_{\ast }C_{J_{\alpha}}Y-\nabla
_{Y}F_{\ast }C_{J_{\alpha}}X,F_{\ast }J_{\alpha}V\right) &=&g_{M}\left( {\cal A}
_{X}B_{J_{\alpha}}Y- {\cal A}_{Y}B_{J_{\alpha}}X,J_{\alpha}V\right) \\
&+&g_{M}\left( {\cal H} grad \ln \lambda ,C_{J_{\alpha}}X\right) g_{M}\left(
Y,J_{\alpha}V\right) \\
&-&g_{M}\left( {\cal H} grad \ln \lambda ,C_{J_{\alpha}}Y\right) g_{M}\left(
X,J_{\alpha}V\right) \\
&-&2g_{M}\left( C_{J_{\alpha}}X,Y\right) g_{M}\left( {\cal H} grad \ln \lambda
,J_{\alpha}V\right) \\
&+&g_{M}\left( \left( \nabla _{Y}J_{\alpha}\right) X-\left( \nabla _{X}J_{\alpha}\right)
Y,J_{\alpha}V\right) .
\end{eqnarray*}
\newline
\noindent
For ${\alpha}={1}$  
\begin{eqnarray}
\lambda ^{-2}g_{N}\left( \nabla _{Y}F_{\ast }C_{J_{1}}X-\nabla
_{X}F_{\ast }C_{J_{1}}Y,F_{\ast }J_{1}V\right) &=&g_{M}\left( {\cal A}
_{X}B_{J_{1}}Y-{\cal A}_{Y}B_{J_{1}}X,J_{1}V\right)  \nonumber \\
&+&g_{M}\left( {\cal H} grad \ln \lambda ,C_{J_{1}}X\right) g_{M}\left(
Y,J_{1}V\right)  \nonumber \\
&-&g_{M}\left( {\cal H} grad \ln \lambda ,C_{J_{1}}Y\right) g_{M}\left(
X,J_{1}V\right)  \nonumber \\
&-&2g_{M}\left( C_{J_{1}}X,Y\right) g_{M}\left( {\cal H} grad \ln \lambda
,J_{1}V\right)  \nonumber \\
&+&g_{M}\left( \left( \nabla _{Y}J_{1}\right) X-\left( \nabla
_{X}J_{1}\right) Y,J_{1}V\right).  \label{eq-dd}
\end{eqnarray}
Now,  
\begin{eqnarray}
g_{M}\left( \left( \nabla _{Y}J_{1}\right) X-\left( \nabla _{X}J_{1}\right)
Y,J_{1}V\right) &=&g_{M}\left( \omega _{3}(Y)J_{2}X-\omega
_{2}(Y)J_{3}X-\omega _{3}(X)J_{2}Y+\omega _{2}(X)J_{3}Y,J_{1}V\right) 
\nonumber \\
&=&\omega _{3}(Y)g_{M}\left( J_{2}X,J_{1}V\right) -\omega _{2}(Y)g_{M}\left(
J_{3}X,J_{1}V\right)  \nonumber \\
&-&\omega _{3}(X)g_{M}\left( J_{2}Y,J_{1}V\right) +\omega _{2}(X)g_{M}\left(
J_{3}Y,J_{1}V\right)  \nonumber \\
&=&\omega _{3}(Y)g_{M}\left( B_{J_{2}}X + C_{J_{2}}X,J_{1}V\right) - \omega _{2}(Y)g_{M}\left( B_{J_{3}}X+C_{J_{3}}X,J_{1}V\right) \nonumber \\
&-& \omega
_{3}(X)g_{M}\left( B_{J_{2}}Y+C_{J_{2}}Y,J_{1}V\right)  +\omega _{2}(X)g_{M}\left( B_{J_{3}}Y+C_{J_{3}}Y,J_{1}V\right)\nonumber \\
&=&\omega _{3}(Y)g_{M}\left( C_{J_{2}}X,J_{1}V\right) -\omega
_{2}(Y)g_{M}\left( C_{J_{3}}X,J_{1}V\right)  \nonumber \\
&-&\omega _{3}(X)g_{M}\left( C_{J_{2}}Y,J_{1}V\right) +\omega
_{2}(X)g_{M}\left( C_{J_{3}}Y,J_{1}V\right).  \label{eq-ddd}
\end{eqnarray}
Using (\ref{eq-dd}) and (\ref{eq-ddd}), we get the result. Similarly, we can
find the other results.

\begin{theorem}
Let $F$ be an almost $h$-conformal semi-invariant submersion from a quaternionic 
K\"ahler manifold $\left( M,J_{1},J_{2},J_{3},g_{M}\right) $ onto a 
Riemannian manifold $\left( N,g_{N}\right) $ such that $\{J_{1},J_{2},J_{3} 
\} $ is an almost $h$-conformal semi-invariant basis. Also, assume that the 
distribution $\left( {\operatorname{ker}}F_{\ast }\right) ^{\perp }$ is integrable. Then 
the following conditions are equivalent:\newline
{\rm (a)} the map $F$ is $h$-homothetic,\newline
{\rm (b)} for $X,Y\in \Gamma \left( \left( {\operatorname{ker}}F_{\ast }\right) ^{\perp
}\right) $ and $V\in \Gamma \left( {\mathcal{D}}_{2}^{J_{1}}\right) $  
\begin{eqnarray*}
\lambda ^{-2}g_{N}\left( \nabla _{Y}^{N}F_{\ast }C_{J_{1}}X-\nabla
_{X}^{N}F_{\ast }C_{J_{1}}Y,F_{\ast }J_{1}V\right) &=&g_{M}\left( {\cal A}
_{X}{\cal B}_{J_{1}}Y-{\cal A}_{Y}{\cal B}_{{J}_{1}}X,J_{1}V\right) \\
&+&\omega _{3}(Y)g_{M}\left( {\cal C}_{J_{2}}X,J_{1}V\right) -\omega
_{2}(Y)g_{M}\left( {\cal C}_{J_{3}}X,J_{1}V\right) \\
&-&\omega _{3}(X)g_{M}\left( {\cal C}_{J_{2}}Y,J_{1}V\right) +\omega
_{2}(X)g_{M}\left( {\cal C}_{J_{3}}Y,J_{1}V\right) ,
\end{eqnarray*}
{\rm (c)} for $X,Y\in \Gamma \left( \left( {\operatorname{ker}}F_{\ast }\right) ^{\perp
}\right) $ and $V\in \Gamma \left( {\mathcal{D}}_{2}^{J_{2}}\right) $  
\begin{eqnarray*}
\lambda ^{-2}g_{N}\left( \nabla _{Y}^{N}F_{\ast }C_{J_{2}}X-\nabla
_{X}^{N}F_{\ast }C_{J_{2}}Y,F_{\ast }J_{2}V\right) 
&=&g_{M}\left( {\cal A}_{X}{\cal B}_{J_{2}}Y-{\cal A}_{Y}{\cal B}
_{J_{2}}X,J_{2}V\right) \\
&+&\omega _{1}(Y)g_{M}\left( {\cal C}_{J_{3}}X,J_{2}V\right) -\omega
_{3}(Y)g_{M}\left( {\cal C}_{J_{1}}X,J_{2}V\right) \\
&-&\omega _{1}(X)g_{M}\left( {\cal C}_{J_{3}}Y,J_{2}V\right) +\omega
_{3}(X)g_{M}\left( {\cal C}_{J_{1}}Y,J_{2}V\right) ,
\end{eqnarray*}
{\rm (d)} for $X,Y\in \Gamma \left( \left( {\operatorname{ker}}F_{\ast }\right) ^{\perp
}\right) $ and $V\in \Gamma \left( {\mathcal{D}}_{2}^{J_{3}}\right) $  
\begin{eqnarray*}
\lambda ^{-2}g_{N}\left( \nabla _{Y}^{N}F_{\ast }C_{J_{3}}X-\nabla
_{X}^{N}F_{\ast }C_{J_{3}}Y,F_{\ast }J_{3}V\right) &=&g_{M}\left( {\cal A}
_{X}{\cal B}_{J_{3}}Y-{\cal A}_{Y}{\cal B} _{J_{3}}X,J_{3}V\right) \\
&+&\omega _{2}(Y)g_{M}\left( {\cal C}_{J_{1}}X,J_{3}V\right) -\omega
_{1}(Y)g_{M}\left( {\cal C}_{J_{2}}X,J_{3}V\right) \\
&-&\omega _{2}(X)g_{M}\left( {\cal C}_{J_{1}}Y,J_{3}V\right) +\omega_{1}(X)g_{M}\left( {\cal C}_{J_{2}}Y,J_{3}V\right).
\end{eqnarray*}
\end{theorem}

\noindent {\bf Proof:} Given $X,Y\in \Gamma \left( \left( {\operatorname{ker}}F_{\ast
}\right) ^{\perp }\right) $, $V\in \Gamma \left( {\mathcal{D}}_{2}^{R}\right) $ 
and $R\in \{J_{1},J_{2},J_{3}\}$. The map is $h$-homothetic, by using the Theorem \ref{th-3}, we obtain the 
results.

\noindent Conversely, from theorem \ref{th-3} (b), [(c) or (d)], we obtain  
\begin{eqnarray}
0 &=&g_{M}\left( {\cal H} grad \ln \lambda ,C_{R}X\right) g_{M}\left(
Y,RV\right) -g_{M}\left( {\cal H} grad \ln \lambda ,C_{R}Y\right)
g_{M}\left( X,RV\right)  \nonumber \\
&&-2g_{M}\left( C_{R}X,Y\right) g_{M}\left( {\cal H} grad \ln \lambda
,RV\right) .  \label{eq-e}
\end{eqnarray}
Replace $Y$ by $RV$ in (\ref{eq-e}) and then using $g_{M}\left( 
C_{R}X,RV\right) =0$, we get  
\[
g_{M}\left( {\cal H} grad \ln \lambda ,C_{R}X\right) g_{M}\left( 
RV,RV\right) =0,  
\]
\[
g_{M}\left( {\cal H} grad \ln \lambda ,C_{R}X\right) =0,  
\]
\begin{equation}
g_{M}\left( \nabla (\lambda ),X\right) =0,\quad {\rm for}\quad X\in \Gamma
\left( \mu ^{R}\right) .  \label{eq-ee}
\end{equation}
Replacing $Y$ by $C_{R}X$ for $X\in \Gamma \left( \mu ^{R}\right) $ in (\ref%
{eq-e}), we get  
\[
-2g_{M}\left( C_{R}X,C_{R}X\right) g_{M}\left( {\cal H} grad \ln \lambda 
,RV\right) =0,  
\]
\[
g_{M}\left( {\cal H} grad \ln \lambda ,RV\right) =0,  
\]
\begin{equation}
g_{M}\left( \nabla (\lambda ),RV\right) =0,\quad {\rm for}\quad V\in \Gamma
\left( {\cal D}_{2}^{R}\right) .  \label{eq-eee}
\end{equation}
By (\ref{eq-ee}) and (\ref{eq-eee}), we can say that $\lambda $ is constant.
\newline
Similarly, we get the other implications.

\begin{cor}
Let $F$ be an almost $h$-conformal anti-holomorphic semi-invariant submersion from a quaternionic K\"{a}%
hler manifold $\left( M,J_{1},J_{2},J_{3},g_{M}\right) $ onto a Riemannian 
manifold $\left( N,g_{N}\right) $ such that $\{J_{1},J_{2},J_{3}\}$ is an almost $h$-conformal anti-holomorphic semi-invariant  basis. Then the following conditions are equivalent:  
\newline
{\rm (a)} The distribution $\left( {\operatorname{ker}}F_{\ast }\right) ^{\perp }$ is 
integrable. \newline
{\rm (b)} ${\cal A}_{J_{1}V_{1}}J_{1}V_{2}= {\cal A}_{J_{1}V_{2}}J_{1}V_{1}$ for $ V_{1},V_{2} \in \Gamma({\mathcal{D}}_{2}^{J_{1}})$. \newline
{\rm (c)} ${\cal A}_{J_{3}V_{1}}J_{3}V_{2}= {\cal A}_{J_{3}V_{2}}J_{3}V_{1}$ for $ V_{1},V_{2} \in \Gamma({\mathcal{D}}_{2}^{J_{3}})$.
\end{cor}

\noindent {\bf Proof:} For anti-holomorphic semi-invariant submersion, we have $B_{R}=R$ and $C_{R}=0$ on $\left( {\operatorname{ker}}F_{\ast
}\right) ^{\perp }$ and $ \omega_{R}=R$ on $ {\mathcal{D}}_{2}^{R}$ for $R\in \{J_{1},J_{3}\}$. 
Applying $ X=RV_{1}$  and $ Y= RV_{2}, V_{1}, V_{2} \in \Gamma({\mathcal{D}}_{2}^R), $ in Theorem \ref{th-3}(Case (b)).
For $R=J_{1}$, we obtain  
$$
{\mathcal{A}}_{J_{1} {V_1}} J_{1} V_2-{\mathcal{A}}_{J_{1}{V_2}} J_{1} V_1 \in \Gamma\left({\mathcal{D}}_2^{J_{1}}\right)
$$
and
$$
0=g_M\left({\mathcal{A}}_{J_{1} V_2} J_{1} V_1-{\mathcal{A}}_{J_{1} V_1} J_{1} V_2, V\right) \text { for } V \in \Gamma\left({\mathcal{D}}_2^{J_{1}}\right),
$$
which are equivalent to
$$
{\mathcal{A}}_{J_{1} V_1} J_{1} V_2={\mathcal{A}}_{J_{1} V_2} J_{1} V_1 \quad \text { for } V_1, V_2 \in \Gamma\left({\mathcal{D}}_2^{J_{1}}\right) .
$$
Hence, we get $(a) \Leftrightarrow(b).$ 
Similarly, we can obtain $(a) \Leftrightarrow(c)$.\\

\begin{theorem}
\label{th-4} Let $F$ be an almost $h$-conformal semi-invariant submersion from a 
quaternionic K\"{a}hler manifold $\left( M,J_{1},J_{2},J_{3},g_{M}\right) $ 
onto a Riemannian manifold $\left( N,g_{N}\right) $ such that $
(J_{1},J_{2},J_{3})$ is an almost $h$-conformal semi-invariant basis. Then the 
following conditions are equivalent:\newline
{\rm (a)} The distribution $\left( {\operatorname{ker}}F_{\ast }\right) ^{\perp }$ defines 
a totally geodesic foliation on $M$. \newline
{\rm (b)} $ {\mathcal{A}}_{X}C_{J_{1}}Y + {\mathcal{V}}\nabla_{X}B_{J_{1}}Y \in \Gamma \left( {\mathcal{D}}_{2}^{J_{1}}\right) $  and 
\begin{eqnarray*}
g_{N}\left( F_{\ast }C_{J_{1}}Y, \nabla _{X}F_{\ast
}J_{1}V\right) &=& \lambda^{2}g_{M}\left({\mathcal{A}}_{X}B_{J_{1}}Y + 
 g_{M}\left( C_{J_{1}}Y,X\right) grad\left(ln \lambda \right) -g_{M} \left( C_{J_{1}}Y,ln \lambda \right)X, J_{1}V  \right)\\
 &-& \lambda^{2}\left( g_{M} \left( \omega_{3}(X)C_{J_{2}}Y,J_{1}V\right) + g_{M} \left( \omega_{2} (X)C_{J_{3}}Y,J_{1}V\right)\right),
\end{eqnarray*}
for $X,Y\in \Gamma \left( \left( {\operatorname{ker}}F_{\ast }\right) ^{\perp
}\right) $ and $V\in \Gamma \left( {\mathcal{D}}_{2}^{J_{1}}\right).$ \\{\rm (c)} $ {\mathcal{A}}_{X}C_{J_{2}}Y + {\mathcal{V}}\nabla_{X}B_{J_{2}}Y \in \Gamma \left( {\mathcal{D}}_{2}^{J_{2}}\right)$  and 
\begin{eqnarray*}
g_{N}\left( F_{\ast }C_{J_{2}}Y, \nabla _{X}F_{\ast
}J_{2}V\right) &=& \lambda^{2}g_{M}\left({\mathcal{A}}_{X}B_{J_{2}}Y + 
 g_{M}\left( C_{J_{2}}Y,X\right) grad\left(ln \lambda \right) -g_{M} \left( C_{J_{2}}Y,ln \lambda \right)X, J_{2}V  \right)\\
 &-& \lambda^{2}\left( g_{M} \left( \omega_{3}(X)C_{J_{3}}Y,J_{2}V\right) + g_{M} \left( \omega_{2} (X)C_{J_{1}}Y,J_{2}V\right)\right),
\end{eqnarray*}
for $X,Y\in \Gamma \left( \left( {\operatorname{ker}}F_{\ast }\right) ^{\perp
}\right) $ and $V\in \Gamma \left( {\mathcal{D}}_{2}^{J_{2}}\right).$\\
{\rm (d)} $ {\mathcal{A}}_{X}C_{J_{3}}Y + {\mathcal{V}}\nabla_{X}B_{J_{3}}Y \in \Gamma \left( {\mathcal{D}}_{2}^{J_{3}}\right)$  and 
\begin{eqnarray*}
g_{N}\left( F_{\ast }C_{J_{3}}Y, \nabla _{X}F_{\ast
}J_{3}V\right) &=& \lambda^{2}g_{M}\left({\mathcal{A}}_{X}B_{J_{3}}Y + 
 g_{M}\left( C_{J_{3}}Y,X\right) grad\left(ln \lambda \right) -g_{M} \left( C_{J_{3}}Y,ln \lambda \right)X, J_{3}V  \right)\\
 &-& \lambda^{2}\left( g_{M} \left( \omega_{3}(X)C_{J_{1}}Y,J_{3}V\right) + g_{M} \left( \omega_{2} (X)C_{J_{2}}Y,J_{3}V\right)\right),
\end{eqnarray*}
for $X,Y\in \Gamma \left( \left( {\operatorname{ker}}F_{\ast }\right) ^{\perp
}\right) $ and $V\in \Gamma \left( {\mathcal{D}}_{2}^{J_{3}}\right).$
\end{theorem}

\noindent {\bf Proof:}  For $X,Y\in \Gamma \left( \left( {\operatorname{ker}}F_{\ast }\right) ^{\perp }\right) $ 
and $W\in \Gamma \left( {\mathcal{D}}_{1}^{R}\right) $ 
  \begin{eqnarray*}
g_{M}\left( \nabla _{X}Y,W\right) &=&g_{M}\left( R\left( \nabla _{X}Y\right)
,RW\right) \\
&=&g_{M}\left( \nabla _{X}RY-\left( \nabla _{X}R\right) Y,RW\right) \\
&=&g_{M}\left( \nabla _{X}RY,RW\right) -g_{M}\left( \left( \nabla
_{X}R\right) Y,RW\right) \\
&=&g_{M}\left( \nabla_{X}(B_{R}Y+ C_{R}Y),RW\right) -g_{M}\left((\nabla_{X}R)Y,RW\right) \\
&=&g_{M}\left( \nabla_{X}B_{R}Y,RW\right) +g_{M}\left( \nabla_{X}C_{R}Y,RW\right) -g_{M}\left( \left( \nabla _{X}R\right) Y,RV\right) \\
&=&g_{M}\left( {\cal A}_{X}B_{R}Y+  {\mathcal{V}} \nabla_{X}B_{R}Y,RW \right)+g_{M}\left( {\cal H}\nabla _{X}C_{R}Y + {\mathcal{A}}_{X}C_{R}Y,RW\right) \\
&& -g_{M}\left( \left(\nabla _{X}R\right) Y,RW\right) \\
&=&g_{M}\left( {\mathcal{V}}{\nabla}_{X}B_{R}Y,RW\right) +g_{M}\left( {\cal A}_{X}C_{R}Y,RW\right) -g_{M}\left( \left( \nabla _{X}R\right) Y,RW\right) \\
&=& -g_{M} \left( R {\mathcal{V}}\nabla_{X}B_{R}Y + R {\mathcal{A}}_{X}C_{R}Y,W\right)-g_{M}\left( \left( \nabla _{X}R\right) Y,RW\right)\\
&=& -g_{M}\left( \phi_{R} {\mathcal{V}}\nabla_{X}B_{R}Y + \omega_{R}{\mathcal{V}}\nabla_{X}B_{R}Y + \phi_{R}{\mathcal{A}}_{X}C_{R}Y + \omega_{R}{\mathcal{A}}_{X}C_{R}Y, W\right)\\
&-& g_{M}\left( \left( \nabla _{X}R\right) Y,RW\right)\\
&=& -g_{M}\left( \phi_{R} \left({\mathcal{A}}_{X}C_{R}Y +{\mathcal{V}}\nabla_{X}B_{R}Y \right), W\right)-g_{M}\left( \left( \nabla _{X}R\right) Y,RW\right).
\end{eqnarray*}

For $R=J_{1}$, we have 

\begin{eqnarray*}
    g_{M}\left( \nabla_{X}Y,W\right) &=& -g_{M}\left( \phi_{J_{1}} \left({\mathcal{A}}_{X}C_{J_{1}}Y +{\mathcal{V}}\nabla_{X}B_{J_{1}}Y \right), W\right)-g_{M}\left( \left( \nabla _{X}J_{1}\right) Y,J_{1}W\right)\\
    &=& -g_{M}\left( \phi_{J_{1}} \left({\mathcal{A}}_{X}C_{J_{1}}Y +{\mathcal{V}}\nabla_{X}B_{J_{1}}Y \right), W\right)- g_{M}\left(\omega_{3}(X)J_{2}Y-\omega_{2}(X)J_{3}Y,J_{1}W \right)\\
    &=& -g_{M}\left( \phi_{J_{1}} \left({\mathcal{A}}_{X}C_{J_{1}}Y +{\mathcal{V}}\nabla_{X}B_{J_{1}}Y \right), W\right)- g_{M} \left(\omega_{3}(X)J_{2}Y,J_{1}W \right)\\
    &+& g_{M}\left( \omega_{2}(X) J_{3}Y,J_{1}W\right)\\
    &=&-g_{M}\left( \phi_{J_{1}} \left({\mathcal{A}}_{X}C_{J_{1}}Y +{\mathcal{V}}\nabla_{X}B_{J_{1}}Y \right), W\right)-g_{M}\left( \omega_{3}(X)(B_{J_{2}}Y + C_{{J}_{2}}Y),J_{1}W \right)\\
    &+& g_{M}\left( \omega_{2}(X)(B_{{J}_{3}}Y+ C_{J_{3}}Y),J_{1}W \right)\\
    &=&-g_{M}\left( \phi_{J_{1}} \left({\mathcal{A}}_{X}C_{J_{1}}Y +{\mathcal{V}}\nabla_{X}B_{J_{1}}Y \right), W\right).
\end{eqnarray*}

$g_{M}\left( \nabla_{X}Y,W\right)=0$ which gives $ g_{M}\left( \phi_{J_{1}} \left({\mathcal{A}}_{X}C_{J_{1}}Y +{\mathcal{V}}\nabla_{X}B_{J_{1}}Y \right), W\right)= 0.$ \\

\noindent 
Given that $ V\in \Gamma ({\mathcal{D}}_{2}^{R})$
using 
\begin{eqnarray*}
g_{M}(\nabla_{X}Y,V) &=& g_M(R \nabla_{X}Y,RV)\\
&=& g_{M}(\nabla_{X}RY-(\nabla_{X}R)Y,RV)\\
&=& g_{M}(\nabla_{X}RY,RV)- g_{M}((\nabla_{X}R)Y,RV).
\end{eqnarray*}
Solving $g_{M}(\nabla_{X}RY,RV)$
\begin{eqnarray*}
g_{M}(\nabla_{X}RY,RV)&=& g_{M}\left( {\mathcal{A}}_{X}B_{R}Y + {\mathcal{V}}\nabla_{X}B_{R}Y,RV\right) + g_{M}(\nabla_{X}C_{R}Y,RV)\\
&=& g_{M}\left( {\mathcal{A}}_{X}B_{R}Y ,RV \right)+  g_{M}(\nabla_{X}C_{R}Y,RV)\\
&=& g_{M}\left( {\mathcal{A}}_{X}B_{R}Y ,RV \right)-  g_{M}(C_{R}Y,\nabla_{X}RV)\\
&=& g_{M}\left( {\mathcal{A}}_{X}B_{R}Y ,RV \right)+\lambda^{-2}g_{N}( F_{\ast}C_{R}Y,RV(\ln \lambda) F_{\ast}X\\
&-& g_{M}(X,RV)F_{\ast}\nabla (\ln \lambda)-{\nabla}_{X}^{F} F_{\ast}RV) )\\
&=& g_{M}\left( {\mathcal{A}}_{X}B_{R}Y +g_{M}(C_{R}Y,X) \nabla(\ln \lambda)- C_{R}Y(\ln \lambda)X,RV \right) \\
&-& \lambda^{-2}g_{N}(F_{\ast}C_{R}Y, {\nabla}_{X}^{F}F_{\ast}RV).
\end{eqnarray*}
Now solving \( g_M \left( \left( \nabla _{X}R \right) Y, RV \right) \).

\begin{align*}
-g_M \left( \left( \nabla _{X}R \right) Y, RV \right) &= -g_M\left( \left( \nabla _{X}J_{1} \right) Y, J_{1}V \right) \\
&= -g_M \left( \omega _{3}(X) J_{2}Y - \omega _{2}(X) J_{3}Y, J_{1}V \right) \\
&= -g_M \left( \omega _{3}(X) J_{2}Y, J_{1}V \right) + g_M \left( \omega _{2}(X) J_{3}Y, J_{1}V \right)\\
&= -g_M\left( \omega _{3}(X)\left( B_{J_{2}}Y + C_{J_{2}}Y \right), J_{1}V \right) 
+ g_M\left( \omega _{2}(X)\left( B_{J_{3}}Y + C_{J_{3}}Y \right), J_{1}V \right) \\
&= -g_M\left( \omega _{3}(X)B_{J_{2}}Y, J_{1}V \right) 
- g_M\left( \omega _{3}(X) C_{J_{2}}Y, J_{1}V \right) \\
&\quad+ g_M\left( \omega _{2}(X) B_{J_{3}}Y, J_{1}V \right) 
+ g_M\left( \omega _{2}(X) C_{J_{3}}Y, J_{1}V \right).
\end{align*}

Thus

\begin{align*}
-g_M\left( \left( \nabla _{X}J_{1} \right) Y, J_{1}V \right) 
&= -g_M\left( \omega _{3}(X) C_{J_{2}}Y, J_{1}V \right) + g_M\left( \omega _{2}(X) C_{J_{3}}Y, J_{1}V \right).
\end{align*}

Now

\begin{align}
g_{M}(\nabla _{X}Y,V) 
&= g_{M}\left({\cal{A}}_{X}B_{J_{1}}Y + g_{M}(C_{J_{1}}Y,X)grad (\ln \lambda) - g_{M}(C_{J_{1}}Y,\ln \lambda )X, J_{1}V \right) \nonumber \\
&\quad -\frac{1}{\lambda^2}g_{N}(F_{\ast }C_{J_{1}}Y,\nabla _{X}^{F}F_{\ast }J_{1}V) -g_M\left( \omega _{3}(X) C_{J_{2}}Y, J_{1}V \right) + g_M\left( \omega _{2}(X) C_{J_{3}}Y, J_{1}V \right). \label{eq-eeee}
\end{align}
Now, for \( g_M(\nabla_XY, V) = 0 \), we have

\begin{align*}
 g_N(F_{\ast }C_{J_{1}}Y,\nabla _{X}^{F}F_{\ast }J_{1}V) 
&= \lambda^2 \Big(g_M({\cal{A}}_{X}B_{J_{1}}Y + g_{M}(C_{J_{1}}Y,X)grad (\ln \lambda )  - g_M(C_{J_{1}}Y,\ln \lambda) X, J_{1}V) \Big) \\
&\quad - \lambda^2  \omega _{3}(X) g_M\left( C_{J_{2}}Y, J_{1}V \right) + \lambda^2  \omega _{2}(X) g_M\left( C_{J_{3}}Y, J_{1}V \right).
\end{align*}
Similarly, we can find the other results.

\begin{lem} Let $F$ be an almost h-conformal semi-invariant submersion from a quaternionic K\"ahler manifold $(M,J_{1},J_{2},J_{3},g_{M})$ onto a Riemannian manifold $(N,g_{N})$ such that $(J_{1},J_{2},J_{3} )$ is an almost h-conformal semi-invariant basis. Assume that the distribution $ {\cal{D}}_{2}^{R} $ is parallel along $\left( \operatorname{ker} F_{\ast} \right) ^{\perp}$ for $ R \in \{ J_{1},J_{2},J_{3} \}.$ Then the following conditions are equivalent: \newline
{\rm (a)} The map $F$ is h-homothetic. \newline
{\rm (b)} For $X,Y\in \Gamma \left( \left( {\operatorname{ker}}F_{\ast }\right) ^{\perp
}\right) $ and $V\in \Gamma \left( {\cal{D}}_{2}^{J_{1}} \right) $
\[  \frac{1}{\lambda^{2}}g_{N} \left( {F_{\ast}C_{J_{1}}Y,\nabla}_{X}^{F}F_{\ast}J_{1}V\right)   =  g_{M} \left( {\cal{A}}_{X}B_{J_{1}}Y,J_{1}V\right) -g_{M}(\omega
_{3}(X)C_{J_{2}}Y,J_{1}V)+g_{M}(\omega _{2}(X)C_{J_{3}}Y,J_{1}V).    \]
{\rm (c)} For $X,Y\in \Gamma \left( \left( {\operatorname{ker}}F_{\ast }\right) ^{\perp
}\right) $ and $V\in \Gamma \left( {\cal{D}}_{2}^{J_{2}} \right) $
\[  \frac{1}{\lambda^{2}}g_{N} \left(F_{\ast}C_{J_{2}}Y, {\nabla}_{X}^{F}F_{\ast}J_{2}V\right) =    g_{M} \left( {\cal{A}}_{X}B_{J_{2}}Y,J_{2}V\right) -g_M (\omega
_{1}(X)C_{J_{3}}Y,J_{2}V)+g_{M}(\omega _{3}(X)C_{J_{1}}Y,J_{2}V). \] 
{\rm (d)} For $X,Y\in \Gamma \left( \left( {\operatorname{ker}}F_{\ast }\right) ^{\perp
}\right) $ and $V\in \Gamma \left( {\cal{D}}_{2}^{J_{3}} \right) $
\[  \frac{1}{\lambda ^{2}}g_{N} \left( F_{\ast}C_{J_{3}}Y,{\nabla}_{X}^{F}F_{\ast}J_{3}V\right)=g_{M} \left( {\cal{A}}_{X}B_{J_{3}}Y,J_{3}V\right)  -g_M (\omega
_{2}(X)C_{J_{1}}Y,J_{3}V)+g_{M}(\omega _{1}(X)C_{J_{2}}Y,J_{3}V).     \]
\end{lem}
\noindent {\bf Proof: }
Given $X,Y\in \Gamma \left((\operatorname{ker} F_{\ast }\right)^{\perp}),V \in \Gamma ({\cal{D}}_{2}^{R}),$
and $R\in \left \{ J_{1},J_{2},J_{3}\right\}.$ \newline
By equation (\ref{eq-eeee}), we have
\begin{align*}
  g_{M}(\nabla _{X}Y,V) & =  g_{M}({\cal{A}}_{X}B_{J_{1}}Y+g_{M}(C_{J_{1}}Y,X)grad (\ln \lambda)-g_{M}(C_{J_{1}}Y,\ln \lambda )X,J_{1}V) \\
& -\frac{1}{\lambda ^{2}}g_{N}(F_{\ast
}C_{J_{1}}Y,\nabla _{X}^{F}F_{\ast }J_{1}V)
 -g_{M}(\omega
_{3}(X)C_{J_{2}}Y,J_{1}V)+g_{M}(\omega _{2}(X)C_{J_{3}}Y,J_{1}V).
\end{align*}
Since $g_{M}(\nabla _{X}Y,V)=-g_{M}(Y,\nabla _{X}V)=0$, therefore from the above equation, we get $(a)\Rightarrow (b),(a)\Rightarrow
(c),(a)\Rightarrow (d).$

For converse part, by use of \ref{th-4}, we obtain 
\begin{equation}
-g_{M}(C_{R}Y,grad (\ln \lambda ))g_{M}(X,RV)+g_{M}(X,C_{R}Y)g_{M}(grad
(\ln \lambda ),RV)=0. \label{eq-f}
\end{equation}
Applying $X=C_{R}Y$ in (\ref{eq-f}) and using $g_{M}\left( C_{R}Y,RV\right) 
=0,$ we have

\[ 
g_{M}(C_{R}Y,C_{R}Y)g_{M}(grad (\ln \lambda ),RV)=0,
\]
so
\begin{align}
g_{M}(grad (\ln \lambda ),RV)=0 \hspace{0.15cm} {\text {for}}\hspace{0.15cm} V\in \Gamma ({\cal{D}}_{2}^{R}),\label{eq-f'}
\end{align}  
since $||C_{R}Y||\neq 0.$ 

Applying $X=RV$ in (\ref{eq-f}) and again using $g_{M}\left( 
C_{R}Y,RV\right) =0,$ we get
\[
-g_{M}(C_{R}Y,grad (\ln \lambda ))g_{M}(RV,RV)=0,
\]
so
\begin{align}
 g_{M}(X,grad (\ln \lambda ))=0 \hspace{0.15cm} {\text{for}} \hspace{0.15cm} X\in \Gamma (\mu ^{R}).\label{eq-f''}
 \end{align} 
By (\ref{eq-f'}) and (\ref{eq-f''}), the map $F$ is horizontally-homothetic.\newline 
Hence, we have (b) $\Rightarrow $ (a), (c) $\Rightarrow $ (a), (d) $\Rightarrow $ (a).

\begin{lem}
    Let $F$ be an almost $h$-conformal anti-holomorphic semi-invariant submersion from a quaternionic K\"ahler manifold $(M,J_{1},J_{2},J_{3},g_{M})$ onto a Riemannian manifold $(N,g_{N})$ such that $\{J_{1}, J_{2}, J_{3}\}$ is an almost $h$-conformal anti-holomorphic semi-invariant basis. Then the following conditions are equivalent: \newline
    {\rm (a)} The distribution $\left( \operatorname{ker} F_{\ast}\right)^{\perp}$ defines a totally geodesic foliation on $M.$ \newline
    {\rm (b)} The distribution ${\cal{D}}_{2}^{J_{1}}$ is parallel along $\left( \operatorname{ker} F_{\ast}\right)^{\perp}.$ \newline
    {\rm (c)} The distribution ${\cal{D}}_{2}^{J_{3}}$ is parallel along $\left( \operatorname{ker} F_{\ast}\right)^{\perp}.$
\end{lem}
\noindent {\bf Proof:} 
Since $B_{R}=R$ and $C_{R}=0$ on $(\operatorname{ker} F_{\ast })^{\perp}$ for $R\in \left\{J_{1},J_{3}\right\} .$ \newline
Given $X,Y\in \Gamma ((\operatorname{ker} F_{\ast })^{\perp})$ and $V\in \Gamma ({\cal{D}}_{2}^{R})$,
from Theorem \ref{th-4}, we get
\begin{align*}
 g_{N}(F_{\ast }C_{J_{1}}Y,\nabla _{X}^{F}F_{\ast }J_{1}V) &= \lambda
^{2}g_{M}({\cal{A}}_{X}B_{J_{1}}Y+g_{M}(C_{J_{1}}Y,X)grad (\ln \lambda )-g_{M}(C_{J_{1}}Y,\ln
\lambda )X,J_{1}V)\\
& -g_{M}(\omega _{3}(X)C_{J_{2}}Y,J_{1}V)+g_{M}(\omega _{2}(X)C_{J_{3}}Y,J_{1}V) 
\end{align*} 
for $R=J_{1}$. 
Putting $C_{J_{1}}=0$ and $B_{J_{1}}=J_{1}.$

$(a)\Leftrightarrow {\mathcal{V}} \nabla _{X}J_{1}Y\in \Gamma ({\cal{D}}_{2}^{J_{1}})$ and $%
g_{M}({\cal{A}}_{X}J_{1}Y,J_{1}V)=0.$

$(a)\Leftrightarrow \nabla _{X}J_{1}Y\in \Gamma ({\cal{D}}_{2}^{J_{1}}).$\\
Similarly, we can find $ (a)\Leftrightarrow (c)$ 
Therefore, we obtain the result.

\begin{theorem}
\label{th-7}    Let $F$ be an almost $h$-conformal semi-invariant submersion from a quaternionic K\"ahler manifold $(M,J_{1},J_{2},J_{3},g_{M})$ onto a Riemannian manifold $(N,g_{N})$ such that $\{J_{1}, J_{2}, J_{3}\}$ is an almost $h$-conformal semi-invariant basis. Then the following conditions are equivalent: \newline
    {\rm (a)} The distribution $\operatorname{ker} F_{\ast}$ defines a totally geodesic foliation on $M.$ \newline
    {\rm (b)}  ${\cal{T}}_{V}\omega_{J_{1}}U + \widehat{\nabla}_{V}\phi_{J_{1}}U \in \Gamma \left( {\cal{D}}_{1}^{J_{1}}\right) $ and \newline
    $ g_{N}\left( {\nabla}_{\omega_{J_{1}}V}^{F} F_{\ast}X, F_{\ast} \omega_{J_{1}}U\right) = \lambda^{2} g_{M} \left( C_{J_{1}}{\cal{T}}_{U}\phi_{J_{1}}V + {\cal{A}}_{\omega_{J_{1}}V}\phi_{J_{1}}U + g_{M}\left( \omega_{J_{1}}V,\omega_{J_{1}}U\right) grad (\ln \lambda),X \right)$  \newline
    for $U,V \in \Gamma(\operatorname{ker} F_{\ast})$ and $X \in \Gamma (\mu ^{J_{1}}).$ \newline
    {\rm (c)} ${\cal{T}}_{V}\omega_{J_{2}}U + \widehat{\nabla}_{V}\phi_{J_{2}}U \in \Gamma \left( {\cal{D}}_{1}^{J_{2}}\right) $ and \newline
    $ g_{N}\left( {\nabla}_{\omega_{J_{2}}V}^{F} F_{\ast}X, F_{\ast} \omega_{J_{2}}U\right) = \lambda^{2} g_{M} \left( C_{J_{2}}{\cal{T}}_{U}\phi_{J_{2}}V + {\cal{A}}_{\omega_{J_{2}}V}\phi_{J_{2}}U + g_{M}\left( \omega_{J_{2}}V,\omega_{J_{2}}U\right) grad (\ln \lambda),X \right)$  \newline
    for $U,V \in \Gamma(\operatorname{ker} F_{\ast})$ and $X \in \Gamma (\mu ^{J_{2}}).$ \newline
    {\rm (d)} ${\cal{T}}_{V}\omega_{J_{3}}U + \widehat{\nabla}_{V}\phi_{J_{3}}U \in \Gamma \left( {\cal{D}}_{1}^{J_{3}}\right) $ and \newline
    $ g_{N}\left( {\nabla}_{\omega_{J_{3}}V}^{F} F_{\ast}X, F_{\ast} \omega_{J_{3}}U\right) = \lambda^{2} g_{M} \left( C_{J_{3}}{\cal{T}}_{U}\phi_{J_{3}}V + {\cal{A}}_{\omega_{J_{3}}V}\phi_{J_{3}}U + g_{M}\left( \omega_{J_{3}}V,\omega_{J_{3}}U\right) grad (\ln \lambda),X \right)$  \newline
    for $U,V \in \Gamma(\operatorname{ker} F_{\ast})$ and $X \in \Gamma (\mu ^{J_{3}}).$ \newline
\end{theorem}

\noindent {\bf Proof:}
Given $U,V\in \Gamma (\operatorname{ker} F_{\ast }), W\in \Gamma ({\cal{D}}_{2}^{R})$ and $R\in
\left\{ J_{1},J_{2},J_{3}\right\} $, we have
\begin{align*}
    g_{M}(\nabla _{V}U,RW) & =-g_{M}(R\nabla _{V}U,W)\\
& =-g_{M}(\nabla _{V}RU-(\nabla _{V}R)U,W)\\
& =-g_{M}(\nabla _{V}RU,W)+g_{M}((\nabla _{V}R)U,W).
\end{align*}
Solving the first term of R.H.S. $-g_{M}(\nabla _{V}RU,W),$ we have 
\begin{align*}
    -g_{M}(\nabla _{V}RU,W) & =-g_{M}(\nabla _{V}(\phi _{R}U+\omega _{R}U),W)\\
& =-g_{M}(\nabla _{V}\phi _{R}U,W)-g_{M}(\nabla _{V}\omega _{R}U,W)\\
& =-g_{M}(\widehat{\nabla}_{V}\phi _{R}U+{\mathcal{T}}_{V}\phi _{R}U,W)-g_{M}({\mathcal{H}}\nabla
_{V}\omega _{R}U+{\mathcal{T}}_{V}\omega _{R}U,W)\\
& =-g_{M}(\widehat{\nabla} _{V}\phi _{R}U,W)-g_{M}({\mathcal{T}}_{V}\phi
_{R}U,W)-g_{M}({\mathcal{H}}\nabla _{V}\omega _{R}U,W)-g_{M}({\mathcal{T}}_{V}\omega _{R}U,W)\\
&=-g_{M}(\widehat{\nabla} _{V}\phi _{R}U,W)-g_{M}({\mathcal{T}}_{V}\omega _{R}U,W)\\
&=-g_{M}((\widehat{\nabla} _{V}\phi _{R}U+{\mathcal{T}}_{V}\omega _{R}U),W)\\
&=-g_{M}(R(\widehat{\nabla} _{V}\phi _{R}U+{\mathcal{T}}_{V}\omega _{R}U),RW)\\
& =-g_{M}(\phi _{R}(\widehat{\nabla} _{V}\phi _{R}U+{\mathcal{T}}_{V}\omega _{R}U)+\omega
_{R}(\widehat{\nabla} _{V}\phi _{R}U+{\mathcal{T}}_{V}\omega _{R}U),RW)\\
& =-g_{M}(\omega _{R}(\widehat{\nabla} _{V}\phi _{R}U+{\mathcal{T}}_{V}\omega _{R}U),RW).
\end{align*}
So for $R= J_{1},$ we have 
$$-g_{M}(\nabla _{V}J_{1}U,W)=-g_{M}(\omega _{J_{1}}(\widehat{\nabla}
_{V}\phi _{J_{1}}U+{\mathcal{T}}_{V}\omega _{J_{1}}U),J_{1}W).$$
Now, solving $g_{M}((\nabla _{V}R)U,W)$ for $R=J_{1}$, we have 
\begin{align*}
g_{M}((\nabla _{V}J_{1})U,W)& =g_{M}((\omega _{3}(V)J_{2}-\omega _{2}(V)J_{3})U,W)\\
& =g_{M}(\omega _{3}(V)J_{2}U,W)-g_{M}(\omega _{2}(V)J_{3}U,W)\\
& =g_{M}(\omega _{3}(V)(\phi _{J_{2}}U+\omega _{J_{2}}U),W)-g_{M}(\omega _{2}(V)(\phi_{J_{3}}U+\omega _{J_{3}}U),W)\\
& =g_{M}(\omega _{3}(V)\phi _{J_{2}}U,W)+g_{M}(\omega _{3}(V)\omega
_{J_{2}}U,W)-g_{M}(\omega _{2}(V)\phi _{J_{3}}U,W)\\
&-g_{M}(\omega _{2}(V)\omega_{J_{3}}U,W)\\
&=0.
\end{align*}
Combining both, we get
$-g_{M}(\nabla _{V}U,J_{1}W)=-g_{M}(\omega _{J_{1}}(\widehat{\nabla} _{V}\phi_{J_{1}}U+T_{V}\omega _{J_{1}}U),J_{1}W).$\newline
Since $\operatorname{ker} F_{\ast }$ defines the totally geodesic foliation, so
$$g_{M}(\nabla _{V}U,J_{1}W)=0.$$
$$\Leftrightarrow (\widehat{\nabla} _{V}\phi _{J_{1}}U+T_{V}\omega _{J_{1}}U)\in \Gamma
({\cal{D}}_{1}^{J_{1}}).$$ 
Let $X\in \Gamma (\mu ^{R})$
\begin{align*}
g_{M}(\nabla _{U}V,X) &= g_{M}(R\nabla _{U}V,RX)\\
&= g_{M}(\nabla _{U}RV-(\nabla _{U}R)V,RX)\\
&= g_{M}(\nabla _{U}RV,RX)-g_{M}((\nabla _{U}R)V,RX)\\
&= g_{M}(\nabla _{U}(\phi _{R}V+\omega _{R}V),RX)-g_{M}((\nabla _{U}R)V,RX)\\
&= g_{M}(\nabla _{U}\phi _{R}V,RX)+g_{M}(\nabla _{U}\omega
_{R}V,RX)-g_{M}((\nabla _{U}R)V,RX)\\
&= g_{M}(\nabla _{U}\phi _{R}V,RX)+g_{M}(\phi _{R}U,\nabla _{\omega_{R}V}X)+g_{M}(\omega _{R}U,\nabla _{\omega _{R}V}X)\\
&+g_{M}(\nabla _{\omega_{R}V}R)U,X)-g_{M}((\nabla _{U}R)V,RX)\\
&= g_{M}({\cal{T}}_{U}\phi _{R}V+ \widehat{\nabla} _{U}\phi _{R}V,RX)+g_{M}(\phi_{R}U,\nabla _{\omega _{R}V}X)+g_{M}(\omega _{R}U,\nabla _{\omega_{R}V}X)\\
&+ g_{M}(\nabla _{\omega _{R}V}R)U,X)-g_{M}((\nabla _{U}R)V,RX)\\
&= g_{M}({\cal{T}}_{U}\phi _{R}V,RX)+g_{M}(\phi _{R}U,{\cal{A}}_{\omega _{R}V}X+{\cal{H}}\nabla_{\omega _{R}V}X)\\
&+g_{M}(\omega _{R}U,\nabla _{\omega _{R}V}X)+g_{M}(\nabla
_{\omega _{R}V}R)U,X)-g_{M}((\nabla _{U}R)V,RX)\\
&=g_{M}({\cal{T}}_{U}\phi _{R}V,RX)+g_{M}(\phi _{R}U,{\cal{A}}_{\omega _{R}V}X)+g_{M}(\omega_{R}U,\nabla _{\omega _{R}V}X)\\
&+g_{M}(\nabla _{\omega_{R}V}R)U,X)-g_{M}((\nabla _{U}R)V,RX)\\
&=g_{M}({\cal{T}}_{U}\phi _{R}V,RX)+g_{M}(\phi _{R}U,{\cal{A}}_{\omega _{R}V}X)-\frac{1}{\lambda ^{2}}g_{M}(grad (\ln \lambda ),X)g_{N}(F_{\ast \omega
_{R}}V,F_{\ast \omega _{R}}U)\\
&+\frac{1}{\lambda ^{2}}g_{N}(\nabla _{\omega
_{R}V}^{F}F_{\ast }X,F_{\ast \omega _{R}}U)+g_{M}(\nabla _{\omega_{R}V}R)U,X)-g_{M}((\nabla _{U}R)V,RX)\\
&=-g_{M}(R{\cal{T}}_{U}\phi _{R}V,X)+g_{M}(\phi _{R}U,{\cal{A}}_{\omega_{R}V}X)-g_{M}((\omega _{R}V,\omega _{R}U)grad (\ln \lambda ),X)\\
&+\frac{1}{%
\lambda ^{2}}g_{N}(\nabla _{\omega _{R}V}^{F}F_{\ast }X,F_{\ast \omega
_{R}}U)+g_{M}(\nabla _{\omega _{R}V}R)U,X)-g_{M}((\nabla _{U}R)V,RX)\\
&=-g_{M}(B_{R}{\cal{T}}_{U}\phi _{R}V+C_{R}{\cal{T}}_{U}\phi _{R}V,X)+g_{M}(\phi
_{R}U,{\cal{A}}_{\omega _{R}V}X)-g_{M}((\omega _{R}V,\omega _{R}U)grad (\ln \lambda ),X)\\
&+\frac{1}{\lambda ^{2}}g_{N}(\nabla _{\omega _{R}V}^{F}F_{\ast}X,F_{\ast \omega _{R}}U)+g_{M}(\nabla _{\omega _{R}V}R)U,X)-g_{M}((\nabla_{U}R)V,RX)\\
&=g_{M}(-C_{R}{\cal{T}}_{U}\phi _{R}V,X)+g_{M}(\phi _{R}U,{\cal{A}}_{\omega_{R}V}X)-g_{M}((\omega _{R}V,\omega _{R}U)grad (\ln \lambda ),X)\\
&+\frac{1}{%
\lambda ^{2}}g_{N}(\nabla _{\omega _{R}V}^{F}F_{\ast }X,F_{\ast \omega_{R}}U)+g_{M}(\nabla _{\omega _{R}V}R)U,X)-g_{M}((\nabla _{U}R)V,RX)\\
&=g_{M}(-C_{R}T_{U}\phi _{R}V-{\cal{A}}_{\omega _{R}V}\phi _{R}U-g_{M}(\omega
_{R}V,\omega _{R}U)grad (\ln \lambda ),X)\\
&+\frac{1}{\lambda ^{2}}%
g_{N}(\nabla _{\omega _{R}V}^{F}F_{\ast }X,F_{\ast \omega
_{R}}U)+g_{M}(\nabla _{\omega _{R}V}R)U,X)-g_{M}((\nabla _{U}R)V,RX).
\end{align*}
Solving for $R=J_{1}$
\begin{eqnarray*}
g_{M}(\nabla _{U}V,X) &=& g_{M}(-C_{J_{1}}{\cal{T}}_{U}\phi _{J_{1}}V-{\cal{A}}_{\omega _{J_{1}}V}\phi
_{J_{1}}U-g_{M}(\omega _{J_{1}}V,\omega _{J_{1}}U)grad (\ln \lambda ),X)\\
&+& \frac{1}{%
\lambda ^{2}}g_{N}(\nabla _{\omega _{J_{1}}V}^{F}F_{\ast }X,F_{\ast \omega
_{J_{1}}}U)+g_{M}(\nabla _{\omega _{J_{1}}V}J_{1})U,X)-g_{M}((\nabla _{U}J_{1})V,J_{1}X)\\
&=& g_{M}(-C_{J_{1}}{\cal{T}}_{U}\phi _{J_{1}}V-{\cal{A}}_{\omega _{J_{1}}V}\phi _{J_{1}}U-g_{M}(\omega
_{J_{1}}V,\omega _{J_{1}}U)grad (\ln \lambda ),X)\\
&+& \frac{1}{\lambda ^{2}}%
g_{N}(\nabla _{\omega _{J_{1}}V}^{F}F_{\ast }X,F_{\ast \omega_{J_{1}}}U) +g_{M}((\omega _{3}\left( \omega _{J_{1}}V\right) J_{2}-\omega _{2}\left(
\omega _{J_{1}}U\right) J_{3})U,X)\\
&-&g_{M}((\omega _{3}(U)J_{2}-\omega _{2}(U)J_{3})V,J_{1}X)\\
&=& g_{M}(-C_{J_{1}}{\cal{T}}_{U}\phi _{J_{1}}V-{\cal{A}}_{\omega _{J_{1}}V}\phi _{J_{1}}U-g_{M}(\omega
_{J_{1}}V,\omega _{J_{1}}U)grad (\ln \lambda ),X)\\
&+& \frac{1}{\lambda ^{2}}%
g_{N}(\nabla _{\omega _{J_{1}}V}^{F}F_{\ast }X,F_{\ast \omega
_{J_{1}}}U)
+ g_{M}(\omega _{3}\left( \omega _{J_{1}}V\right) (\phi _{J_{2}}U+\omega
_{J_{2}}U),X)\\
&-& g_{M}(\omega _{2}\left( \omega _{J_{1}}U\right) (\phi _{J_{3}}U+\omega
_{J_{3}}U),X)\\
&-& g_{M}(\omega _{3}(U)(\phi _{J_{2}}V+\omega _{J_{2}}V),J_{1}X)+g_{M}(\omega
_{2}(U)(\phi _{J_{3}}V+\omega _{J_{3}}V),J_{1}X)\\
&=&g_{M}(-C_{J_{1}}{\cal{T}}_{U}\phi _{J_{1}}V-{\cal{A}}_{\omega _{J_{1}}V}\phi _{J_{1}}U-g_{M}(\omega
_{J_{1}}V,\omega _{J_{1}}U)grad (\ln \lambda ),X)\\
&+&\frac{1}{\lambda ^{2}}%
g_{N}(\nabla _{\omega _{J_{1}}V}^{F}F_{\ast }X,F_{\ast \omega _{J_{1}}}U).
\end{eqnarray*}
since $\operatorname{ker} F_{\ast}$ defines totally geodesic foliation on M, so $%
g_{M}(\nabla _{U}V,X)=0$

$\Leftrightarrow g_{N}(\nabla _{\omega _{J_{1}}V}^{F}F_{\ast }X,F_{\ast \omega
_{J_{1}}}U)=\lambda ^{2}g_{M}(C_{J_{1}}{\cal{T}}_{U}\phi _{J_{1}}V+{\cal{A}}_{\omega _{J_{1}}V}\phi
_{J_{1}}U+g_{M}(\omega _{J_{1}}V,\omega _{J_{1}}U)grad (\ln \lambda ),X).$
Therefore $(a)\Leftrightarrow
(b)$.\newline
Similarly, we find $(a)\Leftrightarrow
(c),(a)\Leftrightarrow (d).$
    
\begin{lem}
  Let $F$ be an almost $h$-conformal semi-invariant submersion from a quaternionic K\"ahler manifold $(M,J_{1},J_{2},J_{3},g_{M})$ onto a Riemannian manifold $(N,g_{N})$ such that $(J_{3},J_{2},J_{3})$ is an almost $h$-conformal semi-invariant basis. Assume that the distribution $\mu^{R}$ is parallel along $\operatorname{ker} F_{\ast}$ for any $R \in \{J_{3},J_{2},J_{3}\}.$   Then the following conditions are equivalent: \newline
  {\rm (a)} Dilation $\lambda$ is constant on $\mu^{R}.$ \newline
  {\rm (b)}  $ g_{N} \left( \nabla_{\omega_{R}V}^{F}F_{\ast}X,F_{\ast}\omega_{R}U\right) = \lambda^{2}g_{M} \left( C_{R}{\cal{T}}_{U}\phi_{R}V + {\cal{A}}_{\omega_{R}V}\phi_{R}U,X \right)$ 
  for $X \in \Gamma(\mu^{R})$ and \newline $U,V \in \Gamma \left( \operatorname{ker} F_{\ast}\right).$
\end{lem}

\noindent {\bf Proof:} 
Given $X\in \Gamma (\mu ^{R})$ and $U,V\in \Gamma (\operatorname{ker} F_{\ast }),$ 
by using the proof of \ref{th-7}, we have
\begin{eqnarray*}
g_{M}\left( \nabla _{U}V,X\right) &=& g_{M}(-C_{R}{\cal T}_{U}\phi _{R}V-{\cal A}_{{\omega}
_{RV}}\phi _{R}U-g_{M}(\omega _{R}V,\omega _{R}U)grad(\ln \lambda ),X)\\
&+& \frac{1}{\lambda ^{2}}g_{N}(\nabla _{\omega _{R}V}^{F}F_{\ast }X,F_{\ast \omega_{R}}U)\\
&=& g_{M}(-C_{R}{\cal{T}}_{U}\phi _{R}V-{\cal{A}}_{\omega _{R}V}\phi _{R}U,X)+\frac{1}{%
\lambda ^{2}}g_{N}(\nabla _{\omega _{R}V}^{F}F_{\ast }X,F_{\ast \omega
_{R}}U),
\end{eqnarray*}
since $g_{M}(\nabla _{U}V,X)=-g_{M}(V,\nabla _{U}X)=0$. Therefore, we obtain $(a)\Leftrightarrow (b).$\\

Using Theorem \ref{th-4} and Theorem \ref{th-7}, we have the following:

\begin{theorem}
    Let $F$ be an almost $h$-conformal semi-invariant submersion from a quaternionic K\"ahler manifold $(M,J_{1},J_{2},J_{3},g_{M})$ onto a Riemannian manifold  $(N,g_{N})$  such that $\{J_{1}, J_{2}, J_{3}\}$ is an almost $h$-conformal semi-invariant basis. Then the following conditions are equivalent: \newline
    {\rm (a)} $M$ is locally a product Riemannian manifold $M_{\operatorname{ker} F_{\ast}} \times M_{(\operatorname{ker} F_{\ast})^{\perp}}.$ \newline
    {\rm (b)}  $ {\mathcal{A}}_{X}C_{J_{1}}Y + {\mathcal{V}}\nabla_{X}B_{J_{1}}Y \in \Gamma \left( {\mathcal{D}}_{2}^{J_{1}}\right) $  and 
    \begin{eqnarray*}
g_{N}\left( F_{\ast }C_{J_{1}}Y, \nabla _{X}^{N}F_{\ast
}J_{1}V\right) &=& \lambda^{2}g_{M}\left({\mathcal{A}}_{X}B_{J_{1}}Y + 
 g_{M}\left( C_{J_{1}}Y,X\right) grad\left(\ln \lambda \right) -g_{M} \left( C_{J_{1}}Y,\ln \lambda \right)X, J_{1}V  \right)\\
 &-& \lambda^{2}\left( g_{M} \left( \omega_{3}(X)C_{J_{2}}Y,J_{1}V\right) + g_{M} \left( \omega_{2} (X)C_{J_{3}}Y,J_{1}V\right)\right),
\end{eqnarray*}
for $X,Y\in \Gamma \left( \left( {\operatorname{ker}}F_{\ast }\right) ^{\perp
}\right) $ and $V\in \Gamma \left( {\mathcal{D}}_{2}^{J_{1}}\right), $ \newline
    ${\cal{T}}_{V}\omega_{J_{1}}U + \widehat{\nabla}_{V}\phi_{J_{1}}U \in \Gamma \left( {\cal{D}}_{1}^{J_{1}}\right) $ and \newline
    $ g_{N}\left( {\nabla}_{\omega_{J_{1}}V}^{F} F_{\ast}X, F_{\ast} \omega_{J_{1}}U\right) = \lambda^{2} g_{M} \left( C_{J_{1}}{\cal{T}}_{U}\phi_{J_{1}}V + {\cal{A}}_{\omega_{J_{1}}V}\phi_{J_{1}}U + g_{M}\left( \omega_{J_{1}}V,\omega_{J_{1}}U\right) grad (\ln \lambda),X \right)$  \newline
    for $U,V \in \Gamma(\operatorname{ker} F_{\ast})$ and $X \in \Gamma (\mu ^{J_{1}}).$ \newline            
    {\rm (c)} $ {\mathcal{A}}_{X}C_{J_{2}}Y + {\mathcal{V}}\nabla_{X}B_{J_{2}}Y \in \Gamma \left( {\mathcal{D}}_{2}^{J_{2}}\right)$  and 
\begin{eqnarray*}
g_{N}\left( F_{\ast }C_{J_{2}}Y, \nabla _{X}^{N}F_{\ast
}J_{2}V\right) &=& \lambda^{2}g_{M}\left({\mathcal{A}}_{X}B_{J_{2}}Y + 
 g_{M}\left( C_{J_{2}}Y,X\right) grad\left(\ln \lambda \right) -g_{M} \left( C_{J_{2}}Y,\ln \lambda \right)X, J_{2}V  \right)\\
 &-& \lambda^{2}\left( g_{M} \left( \omega_{3}(X)C_{J_{3}}Y,J_{2}V\right) + g_{M} \left( \omega_{2} (X)C_{J_{1}}Y,J_{2}V\right)\right),
\end{eqnarray*}
for $X,Y\in \Gamma \left( \left( {\operatorname{ker}}F_{\ast }\right) ^{\perp
}\right) $ and $V\in \Gamma \left( {\mathcal{D}}_{2}^{J_{2}}\right), $ \newline
    ${\cal{T}}_{V}\omega_{J_{2}}U + \widehat{\nabla}_{V}\phi_{J_{2}}U \in \Gamma \left( {\cal{D}}_{1}^{J_{2}}\right) $ and \newline
    $ g_{N}\left( {\nabla}_{\omega_{J_{2}}V}^{F} F_{\ast}X, F_{\ast} \omega_{J_{2}}U\right) = \lambda^{2} g_{M} \left( C_{J_{2}}{\cal{T}}_{U}\phi_{J_{2}}V + {\cal{A}}_{\omega_{J_{2}}V}\phi_{J_{2}}U + g_{M}\left( \omega_{J_{2}}V,\omega_{J_{2}}U\right) grad (\ln \lambda),X \right)$  \newline
    for $U,V \in \Gamma(\operatorname{ker} F_{\ast})$ and $X \in \Gamma (\mu ^{J_{2}}).$ \newline              
    {\rm (d)} $ {\mathcal{A}}_{X}C_{J_{3}}Y + {\mathcal{V}}\nabla_{X}B_{J_{3}}Y \in \Gamma \left( {\mathcal{D}}_{2}^{J_{3}}\right)$  and 
\begin{eqnarray*}
g_{N}\left( F_{\ast }C_{J_{3}}Y, \nabla _{X}^{N}F_{\ast
}J_{3}V\right) &=& \lambda^{2}g_{M}\left({\mathcal{A}}_{X}B_{J_{3}}Y + 
 g_{M}\left( C_{J_{3}}Y,X\right) grad\left(\ln \lambda \right) -g_{M} \left( C_{J_{3}}Y,\ln \lambda \right)X, J_{3}V  \right)\\
 &-& \lambda^{2}\left( g_{M} \left( \omega_{3}(X)C_{J_{1}}Y,J_{3}V\right) + g_{M} \left( \omega_{2} (X)C_{J_{2}}Y,J_{3}V\right)\right),
\end{eqnarray*}
for $X,Y\in \Gamma \left( \left( {\operatorname{ker}}F_{\ast }\right) ^{\perp
}\right) $ and $V\in \Gamma \left( {\mathcal{D}}_{2}^{J_{3}}\right),$ \newline
    ${\cal{T}}_{V}\omega_{J_{3}}U + \widehat{\nabla}_{V}\phi_{J_{3}}U \in \Gamma \left( {\cal{D}}_{1}^{J_{3}}\right) $ and \newline
    $ g_{N}\left( {\nabla}_{\omega_{J_{3}}V}^{F} F_{\ast}X, F_{\ast} \omega_{J_{3}}U\right) = \lambda^{2} g_{M} \left( C_{J_{3}}{\cal{T}}_{U}\phi_{J_{3}}V + {\cal{A}}_{\omega_{J_{3}}V}\phi_{J_{3}}U + g_{M}\left( \omega_{J_{3}}V,\omega_{J_{3}}U\right) grad (\ln \lambda),X \right)$  \newline
    for $U,V \in \Gamma(\operatorname{ker} F_{\ast})$ and $X \in \Gamma (\mu ^{J_{3}}).$ \newline             
\end{theorem}

\begin{theorem}
\label{th-8} Let $F$ be an $h$-conformal semi-invariant submersion from a quaternionic- K\"ahler manifold $(M,J_{1},J_{2},J_{3},g_{M})$ onto a Riemannian manifold $(N,g_{N})$ such that $\{J_{1}, J_{2}, J_{3}\}$ is an $h$-conformal semi-invariant basis. Then the following conditions are equivalent: \newline
{\rm (a)} The distribution ${\cal{D}}_{1}$\ defines a totally geodesic foliation on $M.$ \newline
{\rm (b)} $(\nabla F_{\ast })(V,J_{1}W)\in \Gamma (F_{\ast }\mu ^{J_{1}}),$
\begin{eqnarray*}
g_{N}((\nabla F_{\ast })(V,J_{1}W),F_{\ast }C_{J_{1}}X) &=&\lambda
^{2}g_{M}(U,{\cal{T}}_{V}\omega _{J_{1}}B_{J_{1}}X)\\ 
&-& \lambda ^{2}(g_{M}(\omega _{3}(V)\phi _{J_{2}}U,J_{1}X)-g_{M}(\omega _{3}(V)\omega_{J_{2}}U,J_{1}X))\\
&+& \lambda^{2}(g_{M}(\omega _{2}(V)\phi _{J_{3}}U,J_{1}X)+g_{M}(\omega _{2}(V)\omega_{J_{3}}U,J_{1}X))   
\end{eqnarray*} 
for $ U,V \in \Gamma ({\cal{D}}_{1}), W \in \Gamma ({\cal{D}}_{2})$ and $ X \in \Gamma \left( \left( \operatorname{ker} F_{\ast} \right)^{\perp}\right).$ \newline
{\rm (c)} $(\nabla F_{\ast })(V,J_{2}W)\in \Gamma (F_{\ast }\mu ^{J_{2}}),$
\begin{eqnarray*}
g_{N}((\nabla F_{\ast })(V,J_{2}W),F_{\ast }C_{J_{2}}X) &=& \lambda
^{2}g_{M}(U,{\cal{T}}_{V}\omega _{J_{2}}B_{J_{2}}X)\\
&-& \lambda ^{2}(g_{M}(\omega _{3}(V)\phi _{J_{3}}U,J_{2}X)-g_{M}(\omega _{3}(V)\omega
_{J_{3}}U,J_{2}X))\\
&+& \lambda{2}(g_{M}(\omega _{2}(V)\phi _{J_{1}}U,J_{2}X)+g_{M}(\omega _{2}(V)\omega
_{J_{1}}U,J_{2}X))   
\end{eqnarray*}
for $ U,V \in \Gamma ({\cal{D}}_{1}), W \in \Gamma ({\cal{D}}_{2})$ and $ X \in \Gamma \left( \left( \operatorname{ker} F_{\ast} \right)^{\perp}\right).$ \newline
{\rm (d)} $(\nabla F_{\ast })(V,J_{3}W)\in \Gamma (F_{\ast }\mu ^{J_{3}}),$
\begin{eqnarray*}
g_{N}((\nabla F_{\ast })(V,J_{3}W),F_{\ast }C_{J_{3}}X) &=& \lambda
^{2}g_{M}(U,{\cal{T}}_{V}\omega _{J_{3}}B_{J_{3}}X)\\
&-& \lambda ^{2}(g_{M}(\omega _{3}(V)\phi _{J_{1}}U,J_{3}X)-g_{M}(\omega _{3}(V)\omega
_{J_{1}}U,J_{3}X)) \\
&+& \lambda^{2}(g_{M}(\omega _{2}(V)\phi _{J_{2}}U,J_{3}X)+g_{M}(\omega _{2}(V)\omega_{J_{2}}U,J_{3}X))
\end{eqnarray*}
for $ U,V \in \Gamma ({\cal{D}}_{1}), W \in \Gamma ({\cal{D}}_{2})$ and $ X \in \Gamma \left( \left( \operatorname{ker} F_{\ast} \right)^{\perp}\right).$
\end{theorem}

\noindent {\bf Proof:}
Given $U,V\in \Gamma \left( {\cal{D}}_{1}\right) ,W\in \Gamma \left({\cal{D}}_{2}\right) $ and $R\in \left\{ J_{1},J_{2},J_{3}\right\}$,
we get
\begin{align*}
g_{M}\left( \nabla _{V}U,W\right) &=g_{M}(R\nabla _{V}U,RW)\\
&=g_{M}(\nabla _{V}RU,RW)-g_{M}((\nabla _{V}R)U,RW)\\
&=g_{M}({\cal{H}}\nabla _{V}RU,RW)-g_{M}((\nabla _{V}R)U,RW)\\
&=\frac{1}{\lambda ^{2}}g_{N}((\nabla _{V}^{F}F_{\ast }RU-(\nabla F_{\ast})(V,RU)),F_{\ast }RW)-g_{M}((\nabla _{V}R)U,RW)\\
&=\frac{1}{\lambda ^{2}}g_{N}(\nabla _{V}^{F}F_{\ast }RU,F_{\ast }RW)-\frac{1%
}{\lambda ^{2}}g_{N}((\nabla F_{\ast })(V,RU),F_{\ast }RW)-g_{M}((\nabla_{V}R)U,RW)\\
&=-\frac{1}{\lambda ^{2}}g_{N}((\nabla F_{\ast })(V,RU),F_{\ast
}RW)-g_{M}((\nabla _{V}R)U,RW),
\end{align*}
solving for $R=J_{1}$, we have
\begin{align*}
g_{M}\left( \nabla _{V}U,W\right) &=-\frac{1}{\lambda ^{2}}g_{N}((\nabla F_{\ast })(V,J_{1}U),F_{\ast }J_{1}W)-g_{M}((\nabla _{V}J_{1})U,J_{1}W)\\
& =-\frac{1}{\lambda ^{2}}g_{N}((\nabla F_{\ast })(V,J_{1}U),F_{\ast
}J_{1}W)-g_{M}(\omega _{3}(V)J_{2}U-\omega _{2}(V)J_{3}U,J_{1}W)\\
&=-\frac{1}{\lambda ^{2}}g_{N}((\nabla F_{\ast })(V,J_{1}U),F_{\ast
}J_{1}W)-g_{M}(\omega _{3}(V)(\phi _{J_{2}}U+\omega _{J_{2}}U),J_{1}W)\\
&+g_{M}(\omega_{2}(V)(\phi _{J_{3}}U+\omega _{J_{3}}U),J_{1}W)\\
&=-\frac{1}{\lambda ^{2}}g_{N}((\nabla F_{\ast })(V,J_{1}U),F_{\ast }J_{1}W).  
\end{align*}
Now since distribution ${\cal{D}}_{1}$ defines a totally geodesic foliation, we have \newline 
$$g_{M}(\nabla _{V}U,W)=0\Leftrightarrow ((\nabla F_{\ast })(V,J_{1}U))\in \Gamma
\left( F_{\ast }\mu ^{J_{1}}\right).$$
Given $X\in \Gamma \left( \left(\operatorname{ker} F_{\ast}\right)^{\perp}\right) ,$ 
we obtain
\begin{eqnarray*}
g_{M}(\nabla _{V}U,X) &=&g_{M}(R\nabla _{V}U,RX)\\
&=& g_{M}(\nabla _{V}RU-(\nabla _{V}R)U,RX)\\
&=& g_{M}(\nabla _{V}RU,RX)-g_{M}((\nabla _{V}R)U,RX)\\
&=& g_{M}(\nabla _{V}RU,B_{R}X+C_{R}X)-g_{M}((\nabla _{V}R)U,RX)\\
&=& g_{M}(\nabla _{V}RU,B_{R}X)+g_{M}(\nabla _{V}RU,C_{R}X)-g_{M}((\nabla_{V}R)U,RX)\\
&=& g_{M}(U,\nabla _{V}RB_{R}X)+g_{M}(\nabla _{V}RU,C_{R}X)-g_{M}((\nabla_{V}R)U,RX)\\
&-&g_{M}(U,(\nabla_{V}R)B_{R}X)\\
&=& g_{M}(U,{\cal{T}}_{V}\omega _{R}B_{R}X)-\frac{1}{\lambda ^{2}}g_{N}((\nabla F_{\ast })(V,RU),F_{\ast }C_{R}X)-g_{M}((\nabla _{V}R)U,RX).
\end{eqnarray*}Solving for $R=J_{1}$, we have 
\begin{eqnarray*}g_{M}(\nabla _{V}U,X)&=& g_{M}(U,{\cal{T}}_{V}\omega _{J_{1}}B_{J_{1}}X)-\frac{1}{\lambda ^{2}}%
g_{N}((\nabla F_{\ast })(V,J_{1}U),F_{\ast }C_{J_{1}}X)-g_{M}((\nabla _{V}J_{1})U,J_{1}X)\\
&=& g_{M}(U,{\cal{T}}_{V}\omega _{J_{1}}B_{J_{1}}X)-\frac{1}{\lambda ^{2}}g_{N}((\nabla F_{\ast })(V,J_{1}U),F_{\ast }C_{J_{1}}X)\\
&-& g_{M}(\omega _{3}(V)J_{2}U-\omega _{2}(V)J_{3}U,J_{1}X)\\
&=& g_{M}(U,{\cal{T}}_{V}\omega _{J_{1}}B_{J_{1}}X)-\frac{1}{\lambda ^{2}}g_{N}((\nabla F_{\ast })(V,J_{1}U),F_{\ast }C_{J_{1}}X)\\
&-&g_{M}(\omega _{3}(V)(\phi _{J_{2}}U+\omega_{J_{2}}U)-\omega _{2}(V)(\phi _{J_{3}}U+\omega _{J_{3}}U),J_{1}X)\\
&=& g_{M}(U,{\cal{T}}_{V}\omega _{J_{1}}B_{J_{1}}X)-\frac{1}{\lambda ^{2}}g_{N}((\nabla F_{\ast })(V,J_{1}U),F_{\ast }C_{J_{1}}X)
- g_{M}(\omega _{3}(V)\omega _{J_{2}}U,J_{1}X)\\
&-&g_{M}(\omega _{3}(V)\phi _{J_{2}}U,J_{1}X)+g_{M}(\omega _{2}(V)\omega _{J_{3}}U,J_{1}X)+g_{M}(\omega _{2}(V)\phi_{J_{3}}U,J_{1}X),
\end{eqnarray*}
since $g_{M}(\nabla _{V}U,X)=0$
\begin{align*}
\Leftrightarrow g_{N}((\nabla F_{\ast })(V,J_{1}U),F_{\ast }C_{J_{1}}X)
& =\lambda^{2}g_{M}(U,{\cal{T}}_{V}\omega _{J_{1}}B_{J_{1}}X) \\
&+ \lambda ^{2}(g_{M}(\omega _{2}(V)\omega _{J_{3}}U,J_{1}X) -g_{M}(\omega _{3}(V)\omega _{J_{2}}U,J_{1}X) ) \\
& + \lambda^{2} \left( g_{M} (\omega_{2}(V) \phi _{J_{3}}U,J_{1}X) - g_{M}(\omega_{3}(V) \phi _{J_{2}}U,J_{1}X)\right).
\end{align*}
Hence we have $ (a) \Leftrightarrow (b).$ \newline
Similarly, we can obtain $ (a) \Leftrightarrow (c),  (a) \Leftrightarrow (d).$

\begin{theorem}
\label{th-9} Let $F$ be an $h$-conformal semi-invariant submersion from a quaternionic- K\"ahler manifold $(M,J_{1},J_{2},J_{3},g_{M})$ onto a Riemannian manifold $(N,g_{N})$ such that $\{J_{1}, J_{2}, J_{3}\}$ is an $h$-conformal semi-invariant basis. Then the following conditions are equivalent:  \newline
{\rm (a)} The distribution ${\cal{D}}_{2}$ defines a totally geodesic foliation on $M.$ \newline
{\rm (b)} $\left( \nabla F_{\ast }\right) (V,J_{1}W)\in \Gamma (F_{\ast }\mu ^{J_{1}}),$
\begin{eqnarray*}
\frac{1}{\lambda ^{2}}g_{N}(\nabla _{J_{1}V}^{F}F_{\ast }J_{1}U,F_{\ast
}J_{1}C_{J_{1}}X) &=& g_{M}(V,B_{J_{1}}{\cal{T}}_{U}B_{J_{1}}X)+g_{M}(U,V)g_{M}({\cal{H}}grad (\ln \lambda),J_{1}C_{J_{1}}X)\\
&-& g_{M}(\omega _{3}(U)B_{J_{2}}V,B_{J_{1}}X) - g_{M}(\omega _{3}(U)C_{J_{2}}V,C_{J_{1}}X)\\
&+& g_{M}(\omega_{2}(U)B_{J_{3}}V,B_{J_{1}}X)+g_{M}(\omega _{2}(U)C_{J_{3}}V,C_{J_{1}}X)
\end{eqnarray*}
for $ U,V \in \Gamma \left( {\cal{D}}_{2}\right), W\in \Gamma \left( {\cal{D}}_{1}\right) $ and $ X \in \Gamma \left( \left( \operatorname{ker} F_{\ast}\right)^{\perp}\right).$ \newline 
{\rm (c)}  $\left( \nabla F_{\ast }\right) (V,J_{2}W)\in \Gamma (F_{\ast }\mu ^{J_{2}}),$
\begin{eqnarray*}
\frac{1}{\lambda ^{2}}g_{N}(\nabla _{J_{2}V}^{F}F_{\ast }J_{2}U,F_{\ast
}J_{2}C_{J_{2}}X)&=& g_{M}(V,B_{J_{2}}{\cal{T}}_{U}B_{J_{2}}X)+g_{M}(U,V)g_{M}({\cal{H}}grad (\ln \lambda),J_{2}C_{J_{2}}X)\\
&-& g_{M}(\omega _{3}(U)B_{J_{3}}V,B_{J_{2}}X)-g_{M}(\omega _{3}(U)C_{J_{3}}V,C_{J_{2}}X)\\
&+& g_{M}(\omega_{2}(U)B_{J_{1}}V,B_{J_{2}}X)+g_{M}(\omega _{2}(U)C_{J_{1}}V,C_{J_{2}}X)   
\end{eqnarray*}
for $ U,V \in \Gamma \left( {\cal{D}}_{2}\right), W\in \Gamma \left( {\cal{D}}_{1}\right) $ and $ X \in \Gamma \left( \left( \operatorname{ker} F_{\ast}\right)^{\perp}\right).$ \newline 
{\rm (d)} $\left( \nabla F_{\ast }\right) (V,J_{3}W)\in \Gamma (F_{\ast }\mu ^{J_{3}}),$
\begin{eqnarray*}
 \frac{1}{\lambda ^{2}}g_{N}(\nabla _{J_{3}V}^{F}F_{\ast }J_{3}U,F_{\ast
}J_{3}C_{J_{3}}X)&=& g_{M}(V,B_{J_{3}}{\cal{T}}_{U}B_{J_{3}}X)+g_{M}(U,V)g_{M}({\cal{H}}grad(\ln \lambda
),J_{3}C_{J_{3}}X)\\
&-& g_{M}(\omega _{3}(U)B_{J_{2}}V,B_{J_{3}}X)-g_{M}(\omega _{3}(U)C_{J_{2}}V,C_{J_{1}}X)\\
&+& g_{M}(\omega_{2}(U)B_{J_{3}}V,B_{J_{1}}X)+g_{M}(\omega _{2}(U)C_{J_{3}}V,C_{J_{1}}X)
\end{eqnarray*}
for $ U,V \in \Gamma \left( {\cal{D}}_{2}\right), W\in \Gamma \left( {\cal{D}}_{1}\right) $ and $ X \in \Gamma \left( \left( \operatorname{ker} F_{\ast}\right)^{\perp}\right).$
\end{theorem}

\noindent {\bf Proof:}
Given $U, V \in \Gamma({\cal{D}}_2), W \in \Gamma({\cal{D}}_1), R \in \{J_{1}, J_{2}, J_{3}\},$ we get

\begin{align*}
g_M\left( \nabla_U V, W \right) &= g_M(R \nabla_U V, RW) \\
&= g_M\left( \nabla_U RV - (\nabla_U R)V, RW \right) \\
&= -g_{M}\left( RV, \nabla_{U}RW \right) -g_{M}\left( (\nabla_{U}R)V,RW\right)\\
&= -\frac{1}{\lambda_{2}}g_{N}\left(F_*RV,F_*\nabla_{U}RW  \right)-g_{M}\left( (\nabla_{U}R)V,RW\right)\\
&= -\frac{1}{\lambda_{2}}g_{N}\left(F_*RV, \nabla_{U}F_*RW-(\nabla F_*) (U,RW) \right)-g_{M}\left( (\nabla_{U}R)V,RW\right)\\
&= \frac{1}{\lambda_{2}}g_{N}\left( F_*RV,( \nabla F_*) (U,RW) \right) - g_{M}\left( (\nabla_{U}R)V,RW\right).
\end{align*}

Solving for $R = J_{1}$, we have

\begin{align*}
g_M\left( \nabla_U V, W \right) &= \frac{1}{\lambda^2} g_N\left( (\nabla F_*)(U, J_{1}W), F_* J_{1}V \right) - g_M\left( (\nabla_U J_{1})V, J_{1}W \right) \\
&= \frac{1}{\lambda^2} g_N\left( (\nabla F_*)(U, J_{1}W), F_* J_{1}V \right) - g_M\left( \omega_3(U) J_{2}V - \omega_2(U) J_{3}V, J_{1}W \right) \\
&= \frac{1}{\lambda^2} g_N\left( (\nabla F_*)(U, J_{1}W), F_* J_{1}V \right) - g_M\left( \omega_3(U) {J_{2}}V, J_{1}W \right) \\
&\quad - g_M\left( \omega_2(U) J_{3}V, J_{1}W \right) \\
&= \frac{1}{\lambda^2} g_N\left( (\nabla F_*)(U, J_{1}W), F_* J_{1}V \right).
\end{align*}

Since distribution ${\cal{D}}_2$ defines a totally geodesic foliation, we have  $ g_M(\nabla_U V, W) = 0.$ 
$\Leftrightarrow (\nabla F_*)(U, J_{1}W) \in \Gamma \left( F_* \mu^J_{1} \right).$ \newline
Given $X \in \Gamma \left( (\operatorname{ker} F_*)^\perp \right),$ we obtain
\begin{align*}
g_M(\nabla_U V, X) &= g_M(R \nabla_U V, RX) \\
&= g_M(\nabla_U RV - (\nabla_U R)V, RX) \\
&= g_M(\nabla_U RV, B_R X + C_R X) - g_M((\nabla_U R)V, RX) \\
&= -g_M(RV, \nabla_U B_R X) - g_M(\nabla_{RV} U, C_R X) - g_M((\nabla_U R)V, RX) \\
&= -g_M(RV, {\cal{T}}_U B_R X + \widehat{\nabla}_U B_R X) - g_M(\nabla_{RV} U, C_R X) - g_M((\nabla_U R)V, RX) \\
&= -g_M(RV, {\cal{T}}_U B_R X) - g_M(\nabla_{RV} U, C_R X) - g_M((\nabla_U R)V, RX) \\
&= g_M(V, R {\cal{T}}_U B_R X) + g_M(U, V) g_M({\cal{H}} grad (\ln \lambda), R C_R X) \\
&\quad - \frac{1}{\lambda^2} g_N(\nabla_{RV}^{N} F_* RU, F_* R C_R X) - g_M((\nabla_U R)V, RX) \\
&= g_M(V, B_R {\cal{T}}_U B_R X) + g_M(U, V) g_M({\cal{H}} grad (\ln \lambda), R C_R X) \\
&\quad - \frac{1}{\lambda^2} g_N(\nabla_{RV}^{F} F_* RU, F_* R C_R X) - g_M((\nabla_U R)V, RX).
\end{align*}

Solving for $R = J_{1}$
\begin{align*}
g_M(\nabla_U V, X) &= g_M(V, B_{J_{1}} {\cal{T}}_U B_{J_{1}} X) + g_M(U, V) g_M({\cal{H}}grad(\ln \lambda), {J_{1}} C_{J_{1}} X) \\
&\quad - \frac{1}{\lambda^2} g_N(\nabla_{J_{1}V}^{F} F_* {J_{1}}U, F_* {J_{1}} C_{J_{1}} X) - g_M((\nabla_U J_{1}) V, {J_{1}}X) \\
&= g_M(V, B_{J_{1}} {\cal{T}}_U B_{J_{1}} X) + g_M(U, V) g_M({\cal{H}} grad(\ln \lambda), {J_{1}} C_{J_{1}} X) \\
&\quad - \frac{1}{\lambda^2} g_N(\nabla_{J_{1}V}^{F} F_* {J_{1}}U, F_* {J_{1}} C_{J_{1}} X) - g_M(\omega_3(U) {J_{2}}V - \omega_2(U) {J_{3}}V, {J_{1}}X) \\
&= g_M(V, B_{J_{1}} {\cal{T}}_U B_{J_{1}} X) + g_M(U, V) g_M({\cal{H}} grad(\ln \lambda), {J_{1}} C_{J_{1}} X)\\ 
&\quad - \frac{1}{\lambda^2} g_N(\nabla_{J_{1}V}^{F} F_* {J_{1}}U, F_* {J_{1}} C_{J_{1}} X)\\
&\quad  - g_M(\omega_3(U) (B_{J_{2}} V + C_{J_{2}} V) - \omega_2(U) (B_{J_{3}} V + C_{J_{3}} V), {J_{1}}X) \\
&= g_M(V, B_{J_{1}} {\cal{T}}_U B_{J_{1}} X) + g_M(U, V) g_M({\cal{H}}grad(\ln \lambda), {J_{1}} C_{J_{1}} X) \\
&\quad - \frac{1}{\lambda^2} g_N(\nabla_{J_{1}V}^{F} F_* {J_{1}}U, F_* {J_{1}} C_{J_{1}} X) - g_M(\omega_3(U)(B_{J_{2}} V + C_{J_{2}} V), {J_{1}}X) \\
&\quad + g_M(\omega_2(U)(B_{J_{3}} V + C_{J_{3}} V), {J_{1}}X) \\
&= g_M(V, B_{J_{1}} {\cal{T}}_U B_{J_{1}} X) + g_M(U, V) g_M({\cal{H}} grad(\ln \lambda), {J_{1}} C_{J_{1}} X) \\
&\quad - \frac{1}{\lambda^2} g_N(\nabla_{J_{1}V}^{F} F_* {J_{1}}U, F_* {J_{1}} C_{J_{1}} X) \\
&\quad - g_M(\omega_3(U) B_{J_{2}} V, B_{J_{1}} X) - g_M(\omega_3(U) C_{J_{2}} V, C_{J_{1}} X) \\
&\quad + g_M(\omega_2(U) B_{J_{3}} V, B_{J_{1}} X) + g_M(\omega_2(U) C_{J_{3}} V, C_{J_{1}} X).
\end{align*}

Since distribution ${\cal{D}}_2$ defines a totally geodesic foliation, we have $g_M(\nabla_U V, X) = 0$\begin{align*}
\Leftrightarrow  \frac{1}{\lambda^2} g_N(\nabla_{{J_{1}}V}^{F} F_* {J_{1}}U, F_* {J_{1}} C_{J_{1}} X) &= g_M(V, B_{J_{1}} {\cal{T}}_U B_{J_{1}} X) + g_M(U, V) g_M({\cal{H}} grad(\ln \lambda), {J_{1}} C_{J_{1}} X) \\
&\quad - g_M(\omega_3(U) B_{J_{2}} V, B_{J_{1}} X) - g_M(\omega_3(U) C_{J_{2}} V, C_{J_{1}} X) \\
&\quad + g_M(\omega_2(U) B_{J_{3}} V, B_{J_{1}} X) + g_M(\omega_2(U) C_{J_{3}} V, C_{J_{1}} X).
\end{align*}

Therefore, $(a) \Leftrightarrow (b)$.

Similarly, we can find $(a) \Leftrightarrow (c)$ and $(a) \Leftrightarrow (d)$.
\medskip

Using Theorem \ref{th-8} and Theorem \ref{th-9}, we have  

\begin{theorem}
Let $F$ be an $h$-conformal semi-invariant submersion from a quaternionic- K\"ahler manifold $(M,{J_{1}},{J_{2}},{J_{3}},g_{M})$ onto a Riemannian manifold $(N,g_{N})$ such that $({J_{1}},{J_{2}},{J_{3}})$ is an $h$-conformal semi-invariant basis. Then the following conditions are equivalent: \newline 
{\rm (a)} The fibres of $F$ are locally product Riemannian manifolds $M_{{\cal{D}}_{1}}\times M_{{\cal{D}}_{2}}.$ \newline
{\rm (b)} $(\nabla F_{\ast })(V,{J_{1}}W)\in \Gamma (F_{\ast }\mu ^{{J_{1}}}),$
\begin{eqnarray*}
g_{N}((\nabla F_{\ast })(V,{J_{1}}W),F_{\ast }C_{{J_{1}}}X) &=& \lambda
^{2}g_{M}(U,{\cal{T}}_{V}\omega _{{J_{1}}}B_{{J_{1}}}X)\\
&-&\lambda ^{2}(g_{M}(\omega _{3}(V)\phi _{{J_{2}}}U,{J_{1}}X)-g_{M}(\omega _{3}(V)\omega_{{J_{2}}}U,{J_{1}}X))\\
&+& \lambda^{2}(g_{M}(\omega _{2}(V)\phi _{{J_{3}}}U,{J_{1}}X)+g_{M}(\omega _{2}(V)\omega_{{J_{3}}}U,{J_{1}}X))
\end{eqnarray*}
for $ V,W \in \Gamma ({\cal{D}}_{1}),$  and $ X \in \Gamma \left( \left( \operatorname{ker} F_{\ast} \right)^{\perp}\right).$ \newline

$\left( \nabla F_{\ast }\right) (V,{J_{1}}W)\in \Gamma (F_{\ast }\mu ^{{J_{1}}}),$
\begin{eqnarray*}
-\frac{1}{\lambda ^{2}}g_{N}(\nabla _{{J_{1}}V}^{F}F_{\ast }{J_{1}}U,F_{\ast
}{J_{1}}C_{{J_{1}}}X)&=& g_{M}(V,B_{{J_{1}}}{\cal{T}}_{U}B_{{J_{1}}}X)+g_{M}(U,V)g_{M}({\cal{H}}grad (\ln \lambda),{J_{1}}C_{{J_{1}}}X)\\
&-& g_{M}(\omega _{3}(U)B_{{J_{2}}}V,B_{{J_{1}}}X)
- g_{M}(\omega _{3}(U)C_{{J_{2}}}V,C_{{J_{1}}}X)\\
&+& g_{M}(\omega_{2}(U)B_{{J_{3}}}V,B_{{J_{1}}}X)+g_{M}(\omega _{2}(U)C_{{J_{3}}}V,C_{{J_{1}}}X)    
\end{eqnarray*}
for $U,V \in \Gamma ({\cal{D}}_{2}), W \in \Gamma ({\cal{D}}_{1}) $ and $X \in \Gamma \left( \left( \operatorname{ker} F_{\ast} \right)^{\perp}\right)$ \newline
{\rm (c)} $(\nabla F_{\ast })(V,{J_{2}}W)\in \Gamma (F_{\ast }\mu ^{{J_{2}}}),$
\begin{eqnarray*}
g_{N}((\nabla F_{\ast })(V,{J_{2}}W),F_{\ast }C_{{J_{2}}}X)&=& \lambda
^{2}g_{M}(U,{\cal{T}}_{V}\omega _{{J_{2}}}B_{{J_{2}}}X)\\
&-& \lambda ^{2}(g_{M}(\omega _{3}(V)\phi _{{J_{3}}}U,{J_{2}}X)-g_{M}(\omega _{3}(V)\omega_{{J_{3}}}U,{J_{2}}X))\\
&+& \lambda^{2}(g_{M}(\omega _{2}(V)\phi _{{J_{1}}}U,{J_{2}}X)+g_{M}(\omega _{2}(V)\omega_{{J_{1}}}U,{J_{2}}X))  
\end{eqnarray*}
for $ V,W \in \Gamma ({\cal{D}}_{1}),$  and $ X \in \Gamma \left( \left( \operatorname{ker} F_{\ast} \right)^{\perp}\right).$ \newline
$\left( \nabla F_{\ast }\right) (V,{J_{2}}W)\in \Gamma (F_{\ast }\mu ^{{J_{2}}}),$
\begin{eqnarray*}
-\frac{1}{\lambda ^{2}}g_{N}(\nabla _{{J_{2}}V}^{F}F_{\ast }{J_{2}}U,F_{\ast
}{J_{2}}C_{{J_{2}}}X)&=& g_{M}(V,B_{{J_{2}}}{\cal{T}}_{U}B_{{J_{2}}}X)+g_{M}(U,V)g_{M}({\cal{H}}grad (\ln \lambda),{J_{2}}C_{{J_{2}}}X)\\
&-& g_{M}(\omega _{3}(U)B_{{J_{3}}}V,B_{{J_{2}}}X)-g_{M}(\omega _{3}(U)C_{{J_{3}}}V,C_{{J_{2}}}X)\\
&+& g_{M}(\omega_{2}(U)B_{{J_{1}}}V,B_{{J_{2}}}X)+g_{M}(\omega _{2}(U)C_{{J_{1}}}V,C_{{J_{2}}}X) 
\end{eqnarray*}
for $U,V \in \Gamma ({\cal{D}}_{2}), W \in \Gamma ({\cal{D}}_{1}) $ and $X \in \Gamma \left( \left( \operatorname{ker} F_{\ast} \right)^{\perp}\right)$ \newline
{\rm (d)} $(\nabla F_{\ast })(V,{J_{3}}W)\in \Gamma (F_{\ast }\mu ^{{J_{3}}}),$
\begin{eqnarray*}
g_{N}((\nabla F_{\ast })(V,{J_{3}}W),F_{\ast }C_{{J_{3}}}X)&=& \lambda
^{2}g_{M}(U,{\cal{T}}_{V}\omega _{{J_{3}}}B_{{J_{3}}}X)\\
&-& \lambda ^{2}(g_{M}(\omega _{3}(V)\phi _{{J_{1}}}U,{J_{3}}X)-g_{M}(\omega _{3}(V)\omega_{{J_{1}}}U,{J_{3}}X))\\
&+& \lambda^{2}(g_{M}(\omega _{2}(V)\phi _{{J_{2}}}U,{J_{3}}X)+g_{M}(\omega _{2}(V)\omega_{{J_{2}}}U,{J_{3}}X))   
\end{eqnarray*} 
for $ V,W \in \Gamma ({\cal{D}}_{1}),$  and $ X \in \Gamma \left( \left( \operatorname{ker} F_{\ast} \right)^{\perp}\right).$ \newline
$\left( \nabla F_{\ast }\right) (V,{J_{3}}W)\in \Gamma (F_{\ast }\mu ^{{J_{3}}}),$
\begin{eqnarray*}
-\frac{1}{\lambda ^{2}}g_{N}(\nabla _{{J_{3}}V}^{F}F_{\ast }{J_{3}}U,F_{\ast
}{J_{3}}C_{{J_{3}}}X)&=& g_{M}(V,B_{{J_{3}}}{\cal{T}}_{U}B_{{J_{3}}}X)+g_{M}(U,V)g_{M}({\cal{H}}grad (\ln \lambda),{J_{3}}C_{{J_{3}}}X)\\
&-&g_{M}(\omega _{3}(U)B_{{J_{2}}}V,B_{{J_{3}}}X)-g_{M}(\omega _{3}(U)C_{{J_{2}}}V,C_{{J_{1}}}X)\\
&+& g_{M}(\omega_{2}(U)B_{{J_{3}}}V,B_{{J_{1}}}X)+g_{M}(\omega _{2}(U)C_{{J_{3}}}V,C_{{J_{1}}}X)   
\end{eqnarray*}
for $U,V \in \Gamma ({\cal{D}}_{2}), W \in \Gamma ({\cal{D}}_{1}) $ and $X \in \Gamma \left( \left( \operatorname{ker} F_{\ast} \right)^{\perp}\right).$ \newline
\end{theorem}
Now, we state the following result given by Baird and Wood \cite{Baird}:

\begin{lem}
Let $F$ be an $h$-conformal submersion from a Riemannian manifold $\left( 
M,g_{M}\right) $ onto a Riemannian manifold $\left( N,g_{N}\right) $ with 
dilation $\lambda $. Let $H$ be the mean curvature vector field of the 
distribution $\operatorname{ker}F_{\ast }$ and $m=\dim \operatorname{ker}F_{\ast }$, $n=\dim N$. Then the 
tension field $\tau (F)$ of $F$ is given by  
\[
\tau (F)=-mF_{\ast }H+(2-n)F_{\ast }( grad \ln \lambda ).  
\]
\end{lem}

Using the above result, we can say that

\begin{cor}
Let $F$ be an almost $h$-conformal semi-invariant submersion from a quaternionic K%
\"{a}hler manifold $\left( M,J_{1},J_{2},J_{3},g_{M}\right) $ onto a 
Riemannian manifold $\left( N,g_{N}\right) $ such that $\{J_{1},J_{2},J_{3} 
\} $ is an almost h-conformal semi-invariant basis. Also, assume that $F$ is 
harmonic with $\dim \operatorname{ker}F_{\ast }>0$ and $\dim N>2$. Then the following 
conditions are equivalent:\newline
{\rm (a)} all the fibers of $F$ are minimal.\newline
{\rm (b)} the map $F$ is $h$-homothetic. 
\end{cor}

\begin{cor}
Let $F$ be an almost $h$-conformal semi-invariant submersion from a quaternionic K%
\"{a}hler manifold $\left( M,J_{1},J_{2},J_{3},g_{M}\right) $ onto a 
Riemannian manifold $\left( N,g_{N}\right) $ such that $\{J_{1},J_{2},J_{3} 
\} $ is an almost $h$-conformal semi-invariant basis. Also, assume that $\dim 
\operatorname{ker}F_{\ast }>0$ and $\dim N=2$. Then the map $F$ is harmonic if and only if \newline
{\rm (a)} all the fibers of $F$ are minimal. \newline
{\rm (b)} The map $F$ is harmonic.
\end{cor}

\begin{theorem}
Let $F$ be an almost $h$-conformal semi-invariant submersion from a quaternionic K%
\"{a}hler manifold $\left( M,J_{1},J_{2},J_{3},g_{M}\right) $ onto a 
Riemannian manifold $\left( N,g_{N}\right) $ such that $\{J_{1},J_{2},J_{3}
\} $ is an almost h-conformal semi-invariant basis. Then the following conditions 
are equivalent:\newline
{\rm (a)} The map $F$ is $h$-homothetic.\newline
{\rm (b)} The map $F$ is $\left( J_{1}{\cal{D}}_{2}^{J_{1}},\mu ^{J_{1}}\right) $-totally geodesic.\newline
({\rm c)} The map $F$ is $\left( J_{2}{\cal{D}}_{2}^{J_{2}},\mu ^{J_{2}}\right) $-totally geodesic.\newline
{\rm (d)} The map $F$ is $\left( J_{3}{\cal{D}}_{2}^{J_{3}},\mu ^{J_{3}}\right) $-totally geodesic. 
\end{theorem}

\noindent {\bf Proof:} $ U \in \Gamma({\mathcal{D}}^{R}_{2}), X \in \Gamma(\mu ^{R})$ and $R \in \{J_{1},J_{2},J_{3}\}$, we have 
\begin{align} \label{eq-i}
    (\nabla F_{*})(RU,X) &= RU(\ln \lambda) F_{*}X +X(\ln \lambda)F_{*}RU-g_{M}(RU,X)F_{*}\nabla(\ln \lambda) \nonumber \\
    &= RU(\ln \lambda) F_{*}X +X(\ln \lambda)F_{*}RU.
\end{align}
Since the map $F$ is $h$-homothetic so we will get (a) $\Rightarrow $ (b), (a) $\Rightarrow $ (c), (a) $\Rightarrow $ (d). \\
Conversely, from equation \ref{eq-i}, we obtain
$$ RU(\ln \lambda)F_{*}X + X (\ln \lambda)F_{*}RU =0.$$
Since $\{F_{*}X,F_{*}RU\} $ is linearly independent for nonzero $X,U$, we have $RU(\ln \lambda) =0$ and $X(\ln \lambda) =0,$ which means (a) $ \Leftarrow$ (b), (a) $ \Leftarrow$ (c), (a) $ \Leftarrow$ (d).\\

\begin{theorem}
 Let $F$ be an almost $h$-conformal semi-invariant submersion from a quaternionic K\"ahler manifold $(M,J_{1},J_{2},J_{3},g_{M})$ onto a Riemannian manifold $(N,g_{N})$ such that $\{J_{1}, J_{2}, J_{3}\}$ is an almost $h$-conformal semi-invariant basis. Then the following conditions are equivalent: \newline 
 \begin{description}
     \item[(a)] The map $F$ is totally geodesic. 
     \item[(b)] For $U \in \Gamma \left( \operatorname{ker} F_{\ast}\right)$ and $W \in \Gamma \left( {\cal{D}}_{2}^{J_{1}}\right)$
\begin{description}
         \item[(i)] $C_{J_{1}}{\cal{T}}{U}J_{1}V + \omega_{J_{1}} \widehat{\nabla}_{U} J_{1}V =\omega_{3}(U)J_{3}V + \omega_{2}(U)J_{2}$ for $U,V \in \Gamma \left({\cal{D}}_{1}^{J_{1}}\right),$
         \item[(ii)] $C_{J_{1}}{\cal{H}}\nabla_{U}J_{1}W + \omega_{J_{1}} {\cal{T}}_{U}J_{1}W= \omega_{3}(U)J_{3}W+ \omega_{2}(U)J_{2}W.$ 
     \end{description}
 \item[(c)] 
For $U \in \Gamma \left( \operatorname{ker} F_{\ast}\right)$ and $W \in \Gamma \left( {\cal{D}}_{2}^{J_{2}}\right)$
 \begin{description}
     \item[(i)] $C_{J_{2}}{\cal{T}}{U}J_{2}V + \omega_{J_{2}} \widehat{\nabla}_{U} J_{2}V =\omega_{3}(U)J_{1}V + \omega_{2}(U)J_{3}$ for $U,V \in \Gamma \left({\cal{D}}_{1}^{J_{2}}\right),$
     \item[(ii)]  $C_{J_{2}}{\cal{H}}\nabla_{U}J_{2}W + \omega_{J_{2}} {\cal{T}}_{U}J_{2}W= \omega_{3}(U)J_{1}W+ \omega_{2}(U)J_{3}W.$ 
 \end{description}
 \item[(d)]
For $U \in \Gamma \left( \operatorname{ker} F_{\ast}\right)$ and $W \in \Gamma \left( {\cal{D}}_{2}^{J_{3}}\right)$ 
 \begin{description}
     \item[(i)] $C_{J_{3}}{\cal{T}}{U}J_{3}V + \omega_{J_{3}} \widehat{\nabla}_{U} J_{3}V =\omega_{3}(U)J_{2}V + \omega_{2}(U)J_{1}$ for $U,V \in \Gamma \left({\cal{D}}_{1}^{J_{3}}\right),$
     \item[(ii)] $C_{J_{3}}{\cal{H}}\nabla_{U}J_{3}W + \omega_{J_{3}} {\cal{T}}_{U}J_{3}W= \omega_{3}(U)J_{2}W+ \omega_{2}(U)J_{1}W.$ 
 \end{description}
 \end{description}
\end{theorem}

\noindent {\bf Proof:}
Given $U,V \in \Gamma \left( {\cal{D}}_{1}^{R}\right) $ and $R\in \left\{ 
J_{1}, J_{2}, J_{3} \right\} $, we have

\begin{align*}
    \left( \nabla F_{\ast} \right) (U,V) 
    &= -F_{\ast }(\nabla _{U}V) \\
    &= F_{\ast }(R^{2}\nabla _{U}V) \\
    &= F_{\ast }(R(\nabla _{U}RV - (\nabla _{U}R)V))\\
    &= F_{\ast }\left(R({\cal{T}}_{U}RV + \widehat{\nabla} _{U}RV) - R(\nabla _{U}R)V\right) \\
    &= F_{\ast }\left(B_{R}{\cal{T}}_{U}RV + C_{R}{\cal{T}}_{U}RV + \phi _{R}\widehat{\nabla} _{U}RV + \omega
    _{R}\widehat{\nabla} _{U}RV - R(\nabla _{U}R)V\right) \\
    &= F_{\ast }\left(C_{R}{\cal{T}}_{U}RV + \omega _{R}\widehat{\nabla} _{U}RV - R(\nabla _{U}R)V\right).
\end{align*}

Now, solve for $R = J_{1}$

\begin{align*}
    \left( \nabla F_{\ast} \right) (U,V) 
    &= F_{\ast }\left(C_{J_{1}}{\cal{T}}_{U}J_{1}V + \omega _{J_{1}}\widehat{\nabla} _{U}J_{1}V - J_{1}(\nabla _{U}J_{1})V\right) \\
    &= F_{\ast }\left(C_{J_{1}}{\cal{T}}_{U}J_{1}V + \omega _{J_{1}}\widehat{\nabla} _{U}J_{1}V - J_{1}\left(\omega _{3}(U)J_{2}V - \omega _{2}(U)J_{3}V\right)\right).
\end{align*}
Thus
\[
    \left( \nabla F_{\ast} \right) (U,V) = 0 \Leftrightarrow
    C_{J_{1}}{\cal{T}}_{U}J_{1}V + \omega _{J_{1}}\widehat{\nabla} _{U}J_{1}V = \omega _{3}(U)J_{3}V + \omega _{2}(U)J_{2}V.
\]
Given  $ U \in \Gamma \left( \operatorname{ker} F_{\ast} \right)$   and $ W \in \Gamma \left( {\cal{D}}_{2}^{R} \right),$   we get
\begin{align*}
\left( \nabla F_{\ast} \right)(U,W) &= -F_{\ast}(\nabla_{U} W) \\
&= -F_{\ast}(R \nabla_{U} W) \\
&= -F_{\ast}(\nabla_{U} RW - (\nabla_{U} R)W) \\
&= F_{\ast}(R(\nabla_{U} RW) - R(\nabla_{U} R)W) \\
&= F_{\ast}(R({\cal{T}}_{U} RW + {\cal{H}} \nabla_{U} RW) - R(\nabla_{U} R)W) \\
&= F_{\ast}(\phi_{R} {\cal{T}}_{U} RW + \omega_{R} {\cal{T}}_{U} RW + B_{R} {\cal{H}} \nabla_{U} RW + C_{R} {\cal{H}} \nabla_{U} RW - R(\nabla_{U} R)W) \\
&= F_{\ast}(\omega_{R} {\cal{T}}_{U} RW + C_{R} {\cal{H}} \nabla_{U} RW - R(\nabla_{U} R)W).
\end{align*}

Solving for $  R = J_{1} $, we have 
\begin{align*}
\left( \nabla F_{\ast} \right)(U,W)&=F_{\ast}\left(\omega_{J_{1}} {\cal{T}}_{U} J_{1} W + C_{J_{1}} {\cal{H}} \nabla_{U} J_{1} W - J_{1}(\nabla_{U} J_{1})W\right) \\
&= F_{\ast} \left(\omega_{J_{1}} {\cal{T}}_{U} J_{1} W + C_{J_{1}} {\cal{H}} \nabla_{U} J_{1} W- J_{1}\left(\omega_{3}(U) J_{2} W - \omega_{2}(U) J_{3} W\right)\right).
\end{align*}
So that

$\left( \nabla F_{\ast} \right)(U,W) = 0 \Leftrightarrow \omega_{J_{1}} {\cal{T}}_{U} J_{1} W + C_{J_{1}} {\cal{H}} \nabla_{U} J_{1} W 
= \omega_{3}(U) J_{3} W + \omega_{2}(U) J_{2} W.$

Similarly, we can solve for $  R = J_{2}$   and $  R = J_{3}.$

\begin{lem}
\label{lem-10}Let $F$ be an almost $h$-conformal semi- invariant submersion with totally umbilical fibers from a quaternionic K\"ahler manifold $(M,J_{1},J_{2},J_{3},g_{M})$ onto a Riemannian manifold $(N,g_{N})$ such that $\{J_{1}, J_{2}, J_{3}\}$ is an almost $h$-conformal semi-invariant basis. Then for $R\in \left\{ 
J_{1},J_{2},J_{3}\right\} $,   mean curvature vector field of the distribution $\operatorname{ker} F_{\ast}$, $H\in \Gamma \left( R{\cal D}_2^R\right)$.
\end{lem}
\noindent {\bf Proof:} $ X,Y \in \Gamma \left( {\cal D}_1\right) , W \in \Gamma \left( \mu\right)$ and $ R \in \{J_{1},J_{2},J_{3}\} $\newline
\begin{equation*}
 \nabla_{X}RY = {\mathcal{T}}_{X}RY+ \widehat{\nabla}_{X}RY  \end{equation*}
 \begin{equation*}
\left( \nabla_{X}R\right)Y + R \nabla_{X}Y = {\mathcal{T}}_{X}RY+ \widehat{\nabla}_{X}RY \end{equation*}
 \begin{equation*}
\left( \nabla_{X}R\right)Y + R{\mathcal{T}}_{X}Y+R\widehat{\nabla}_{X}Y={\mathcal{T}}_{X}RY+ \widehat{\nabla}_{X}RY\end{equation*}
 \begin{equation*}
      \left( \nabla_{X}R\right)Y + B_{R}{\mathcal{T}}_{X}Y+C_{R}{\mathcal{T}}_{X}Y+\phi_{R}\widehat{\nabla}_{X}Y+\omega_{R}\widehat{\nabla}_{X}Y={\mathcal{T}}_{X}RY+ \widehat{\nabla}_{X}RY 
  \end{equation*}
  So,
\begin{align*}
 g_{M}\left( \nabla_{X}RY,W\right)&= g_{M}\left( {\mathcal{T}}_{X}RY,W\right)\\
 &=g_{M}\left((\nabla_{X}R)Y,W \right)+g_{M}\left( C_{R}{\mathcal{T}}_{X}Y,W\right)\\
 &= g_{M}\left((\nabla_{X}R)Y,W \right)-g_{M}\left( {\mathcal{T}}_{X}Y,RW\right).
\end{align*}
Since ${\mathcal{T}}_{X}Y = g_{M}(X,Y)H$.  
Therefore
\begin{equation} \label{eq-h}
 g_{M}(X,RY)g_{M}(H,W)= g_{M}\left((\nabla_{X}R)Y,W \right)-g_{M}\left( {\mathcal{T}}_{X}Y,RW\right).  
\end{equation} 
Interchanging the role of $X$ and $Y,$ we get
\begin{equation} \label{eq-hh}
 g_{M}(Y,RX)g_{M}(H,W)= g_{M}\left((\nabla_{Y}R)X,W \right)-g_{M}\left( {\mathcal{T}}_{Y}X,RW\right).   
\end{equation}
Combining \ref{eq-h} and \ref{eq-hh}, we get
$$g_{M}\left( (\nabla_{X}R)Y + (\nabla_{Y}R)X,W \right)= 2g_{M} (X,Y)g_{M}(H,RW).$$
Therefore $$ 2g_{M}(X,Y)g_{M}(H,RW)=0$$ 
which implies $H \in \Gamma (R{\mathcal{D}}^{R}_{2}).$\\
For $R=J_{1}$ \newline
\begin{align*}
  g_{M}\left( (\nabla_{X}J_{1})Y + (\nabla_{Y}J_{1})X,W\right)&= g_{M}\left( \omega_{3}(X)J_{2}Y - \omega_{2}(X)J_{3}Y + \omega_{3}(Y)J_{2}X - \omega_{2}(Y)J_{3}X,W\right). \\
  &=0.
\end{align*}
Similarly for $R=J_{2}$ and $J_{3}.$
\begin{theorem}
 Let $F$ be an almost $h$-conformal semi-invariant submersion with totally umbilical fibers from a quaternionic K\"ahler manifold $(M,J_{1},J_{2},J_{3},g_{M})$ onto a Riemannian manifold $(N,g_{N})$ such that $\{J_{1}, J_{2}, J_{3}\}$ is an almost $h$-conformal semi-invariant basis. Then all the fibers of $F$ are totally geodesic.   
\end{theorem}

\noindent {\bf Proof:} By Lemma \ref{lem-10}, we have
$ H \in \Gamma (R {\mathcal{D}}_{2})$ for $ R \in \{J_{1},J_{2},J_{3}\}$ \newline
\noindent so that $$\{J_{1}H,J_{2}H,J_{3}H\} \in \Gamma( {\mathcal{D}}_{2}).$$
But $$ J_{3}H=J_{1}J_{2}H=J_{1}(J_{2}H)\in \Gamma({\mathcal{D}}_{2}).$$
Since $J_{1}{\mathcal{D}}_{2} \subset \left( \operatorname{ker} F_{\ast}\right)^{\perp}, $ we must have $ H=0.$

\section{Examples}

\begin{ex-new}
Every $h$-anti-invariant Riemannian submersion is an $h$-conformal 
anti-invariant submersion with dilation $\lambda=id$ {\rm \cite{Park-17}}. So it is $h$-conformal 
semi-invariant submersion with dilation $\lambda=id$ and distribution ${\cal D}_1=0$.
\end{ex-new}

\begin{ex-new}
Let $M$ be a $4n$-dimensional almost quaternionic Hermitian manifold and $N$
be a Riemannian manifold such that $dim \ N=dim \ M -1$. Let $F : 
M\rightarrow N$ be an $h$-conformal submersion with dilation $\lambda$.  Then
map $F$ is an $h$-conformal anti-invariant submersion with dilation $ \lambda
$. So it is an $h$-conformal 
semi-invariant submersion with dilation $\lambda$ and distribution ${\cal D}_1=0$.
\end{ex-new}


\begin{ex-new}
Consider the $4$-dimensional manifold $M=\{(x_{1},x_{2},x_{3},x_{4})\in  
\mathbb{R}^{4},x_{4}\not=0\}$, where $(x_{1},x_{2},x_{3},x_{4})$ are the 
standard coordinates of ${\mathbb R}^{4}$. The vector fields  
\[
e_{1}=x_{3}\frac{\partial }{\partial x_{1}},\quad e_{2}=x_{3}\frac{\partial 
}{\partial x_{2}},\quad e_{3}=-x_{3}\frac{\partial }{\partial x_{3}},\quad 
e_{4}=x_{3}\frac{\partial }{\partial x_{4}}  
\]
are linearly independent at each point of $M$. Let $g$ be the standard 
Riemannian metric. Let $J_{1},J_{2},J_{3}$ be the $(1,1)$-tensor fields 
defined by  
\[ J_{1}\left( \frac{\partial }{\partial x_{1}}\right) =\frac{\partial }{
\partial x_{2}},\ \ \ \ \ \ \ \ \ \ J_{1}\left( \frac{\partial }{\partial
x_{2}}\right) =-\frac{\partial }{\partial x_{1}},  
\ \ \ \ \ \ \ \ \ \ 
J_{1}\left( \frac{\partial }{\partial x_{3}}\right) =\frac{\partial }{
\partial x_{4}},\ \ \ \ \ \ \ \ \ \ J_{1}\left( \frac{\partial }{\partial
x_{4}}\right) =-\frac{\partial }{\partial x_{3}}, \]
\[
J_{2}\left( \frac{\partial }{\partial x_{1}}\right) =\frac{\partial }{
\partial x_{3}},\ \ \ \ \ \ \ \ \ \ J_{2}\left( \frac{\partial }{\partial
x_{2}}\right) =-\frac{\partial }{\partial x_{4}},  
\ \ \ \ \ \ \ \ \ 
J_{2}\left( \frac{\partial }{\partial x_{3}}\right) =-\frac{\partial }{
\partial x_{1}},\ \ \ \ \ \ \ \ J_{2}\left( \frac{\partial }{\partial x_{4}}
\right) =\frac{\partial }{\partial x_{2}},  
\]
\[
J_{3}\left( \frac{\partial }{\partial x_{1}}\right) =\frac{\partial }{
\partial x_{4}},\ \ \ \ \ \ \ \ \ \ J_{3}\left( \frac{\partial }{\partial
x_{2}}\right) =\frac{\partial }{\partial x_{3}},  
\ \ \ \ \ \ \ \ \ 
J_{3}\left( \frac{\partial }{\partial x_{3}}\right) =-\frac{\partial }{
\partial x_{2}},\ \ \ \ \ \ \ \ J_{3}\left( \frac{\partial }{\partial x_{4}}
\right) =-\frac{\partial }{\partial x_{1}}.  
\]
By using Koszul's formula, we have  
\[
\left( \nabla _{e_{1}}J_{1}\right) e_{1}=e_{4},\ \ \ \ \ \ \ \ \left( \nabla
_{e_{2}}J_{1}\right) e_{1}=-e_{3},\ \ \ \ \ \ \ \ \left( \nabla 
_{e_{3}}J_{1}\right) e_{1}=0,\ \ \ \ \ \ \ \ \ \ \left( \nabla 
_{e_{4}}J_{1}\right) e_{1}=0,  
\]
\[
\left( \nabla _{e_{1}}J_{1}\right) e_{2}=e_{3},\ \ \ \ \ \ \ \ \left( \nabla
_{e_{2}}J_{1}\right) e_{2}=e_{4},\ \ \ \ \ \ \ \ \ \ \left( \nabla 
_{e_{3}}J_{1}\right) e_{2}=-e_{1},\ \ \ \ \ \ \ \left( \nabla 
_{e_{4}}J_{1}\right) e_{2}=0,  
\]
\[
\left( \nabla _{e_{1}}J_{1}\right) e_{3}=-e_{2},\ \ \ \ \ \ \left( \nabla 
_{e_{2}}J_{1}\right) e_{3}=e_{1},\ \ \ \ \ \ \ \ \ \ \left( \nabla 
_{e_{3}}J_{1}\right) e_{3}=0,\ \ \ \ \ \ \ \ \ \ \left( \nabla 
_{e_{4}}J_{1}\right) e_{3}=0,  
\]
\[
\left( \nabla _{e_{1}}J_{1}\right) e_{4}=-e_{1},\ \ \ \ \ \ \left( \nabla 
_{e_{2}}J_{1}\right) e_{4}=-e_{2},\ \ \ \ \ \ \ \ \left( \nabla 
_{e_{3}}J_{1}\right) e_{4}=0,\ \ \ \ \ \ \ \ \ \ \left( \nabla 
_{e_{4}}J_{1}\right) e_{4}=0,  
\]
\[
\left( \nabla _{e_{1}}J_{2}\right) e_{1}=0,\ \ \ \ \ \ \ \ \ \left( \nabla 
_{e_{2}}J_{2}\right) e_{1}=e_{2},\ \ \ \ \ \ \ \ \ \ \left( \nabla 
_{e_{3}}J_{2}\right) e_{1}=0,\ \ \ \ \ \ \ \ \ \ \left( \nabla 
_{e_{4}}J_{2}\right) e_{1}=e_{4},  
\]
\[
\left( \nabla _{e_{1}}J_{2}\right) e_{2}=0,\ \ \ \ \ \ \ \ \ \left( \nabla 
_{e_{2}}J_{2}\right) e_{2}=-e_{1},\ \ \ \ \ \ \ \ \left( \nabla 
_{e_{3}}J_{2}\right) e_{2}=0,\ \ \ \ \ \ \ \ \ \left( \nabla 
_{e_{4}}J_{2}\right) e_{2}=e_{3},  
\]
\[
\left( \nabla _{e_{1}}J_{2}\right) e_{3}=0,\ \ \ \ \ \ \ \ \left( \nabla 
_{e_{2}}J_{2}\right) e_{3}=e_{4},\ \ \ \ \ \ \ \ \ \ \left( \nabla 
_{e_{3}}J_{2}\right) e_{3}=0,\ \ \ \ \ \ \ \ \ \ \ \left( \nabla 
_{e_{4}}J_{2}\right) e_{3}=-e_{2},  
\]
\[
\left( \nabla _{e_{1}}J_{2}\right) e_{4}=0,\ \ \ \ \ \ \ \ \ \left( \nabla 
_{e_{2}}J_{2}\right) e_{4}=-e_{3},\ \ \ \ \ \ \ \ \left( \nabla 
_{e_{3}}J_{2}\right) e_{4}=0,\ \ \ \ \ \ \ \ \ \left( \nabla 
_{e_{4}}J_{2}\right) e_{4}=-e_{1},  
\]
\[
\left( \nabla _{e_{1}}J_{3}\right) e_{1}=-e_{2},\ \ \ \ \ \left( \nabla 
_{e_{2}}J_{3}\right) e_{1}=0,\ \ \ \ \ \ \ \ \ \ \left( \nabla 
_{e_{3}}J_{3}\right) e_{1}=0,\ \ \ \ \ \ \ \ \ \left( \nabla 
_{e_{4}}J_{3}\right) e_{1}=-e_{3},  
\]
\[
\left( \nabla _{e_{1}}J_{3}\right) e_{2}=e_{1},\ \ \ \ \ \ \ \left( \nabla 
_{e_{2}}J_{3}\right) e_{2}=0,\ \ \ \ \ \ \ \ \ \ \left( \nabla 
_{e_{3}}J_{3}\right) e_{2}=0,\ \ \ \ \ \ \ \ \ \left( \nabla 
_{e_{4}}J_{3}\right) e_{2}=e_{4},  
\]
\[
\left( \nabla _{e_{1}}J_{3}\right) e_{3}=-e_{4},\ \ \ \ \ \ \left( \nabla 
_{e_{2}}J_{3}\right) e_{3}=0,\ \ \ \ \ \ \ \ \ \left( \nabla 
_{e_{3}}J_{3}\right) e_{3}=0,\ \ \ \ \ \ \ \ \ \ \left( \nabla 
_{e_{4}}J_{3}\right) e_{3}=e_{1},  
\]
\[
\left( \nabla _{e_{1}}J_{3}\right) e_{4}=e_{3},\ \ \ \ \ \ \ \left( \nabla 
_{e_{2}}J_{3}\right) e_{4}=0,\ \ \ \ \ \ \ \ \ \left( \nabla 
_{e_{3}}J_{3}\right) e_{4}=0,\ \ \ \ \ \ \ \ \ \left( \nabla 
_{e_{4}}J_{3}\right) e_{4}=-e_{2}.  
\]
Let $\omega _{1},\omega _{2},\omega _{3},\omega _{4}$ be the $1$-forms 
defined by  
\[
\omega _{1}=\left(  
\begin{array}{cccc}
-1 & 0 & 0 & 0 \\ 
0 & 0 & 0 & 0 \\ 
0 & 0 & 0 & 0 \\ 
0 & 0 & 0 & 0%
\end{array}
\right) ,\qquad \omega _{2}=\left(  
\begin{array}{cccc}
0 & 0 & 0 & 0 \\ 
0 & 0 & 0 & 0 \\ 
0 & 0 & 0 & 1 \\ 
0 & 0 & 0 & 0%
\end{array}
\right) ,  
\]
\[ \omega _{3}=\left(  
\begin{array}{cccc}
0 & 0 & 0 & 0 \\ 
0 & -1 & 0 & 0 \\ 
0 & 0 & 0 & 0 \\ 
0 & 0 & 0 & 0%
\end{array}
\right) ,\qquad \omega _{4}=\left(  
\begin{array}{cccc}
0 & 0 & 0 & 0 \\ 
0 & 0 & 0 & 0 \\ 
0 & 0 & 0 & 0 \\ 
0 & 0 & 0 & 1%
\end{array}
\right) .  \]
From the above results, it is easy to check that $M$ is quaternionic K\"{a}hler manifold \cite{G2}.
\newline
\noindent {\rm (a)} Let $F : M \rightarrow \mathbb {R}^3$ be an $h$-conformal 
submersion with dilation $\lambda$. The map $F$ then is a $h$-conformal antiinvariant submersion with dilation $\lambda$.

\noindent {\rm (b)} Consider a Riemannian submersion $F : M \rightarrow  
\mathbb{R}^2$ by  
\[
F(x_{1},x_{2},x_{3},x_{4})=(e^{x_{3}}sin \ x_{4},e^{x_{3}}cos \ x_{4} )  
\]
Then $F$ is an $h$-conformal semi-invariant submersion, where distribution $ {\mathcal{D}}_{1} = \{0\},$  with dilation $
e^{x_{3}}$.

\noindent {\rm (c)} Consider a Riemannian submersion $ F  : M \rightarrow \mathbb{R}$ by 
\[ F(x_{1},x_{2},x_{3},x_{4})= e^{x_{3}}\]
Then $F$ is an $h$-conformal semi-invariant submersion with dilation $e^{x_{3}}.$
\end{ex-new}
\noindent
\textbf{Declaration:}  We declare that no conflicts of interest are associated with this publication and there has been no significant financial support for this work.  We certify that the submission is original work. 

\end{document}